\documentclass[a4paper,UKenglish,cleveref, autoref, thm-restate]{lipics-v2021}
\usepackage{subfiles}
\usepackage[mode=buildnew]{standalone}

\usepackage{tcolorbox}
\tcbuselibrary{theorems}
\usepackage[autolanguage]{numprint}
\npthousandsep{\,}
\usepackage{ifthen}
\usepackage{pgfmath}
\usepackage{graphicx}
\usepackage{datatool}
\usepackage{dataplot}
\usepackage[group-separator={,}]{siunitx}
\usepackage{xint}
\usepackage{xintexpr}
\usepackage{xfp}
\usepackage{xinttools}
\usepackage{xstring}
\usepackage{makecell}
\usepackage{fp}
\usepackage{listings}
\usepackage{caption}
\usepackage{complexity}
\usepackage{longtable}
\usepackage{pdflscape}
\usepackage{multicol}
\usepackage{multirow}
\usepackage{booktabs}
\usepackage{tabularx,ragged2e}
\usepackage{mathtools}
\usepackage{multicol}
\usepackage{float}
\usepackage{algorithm2e}
\hideLIPIcs 
\LinesNumbered
\graphicspath{ {./figures/} }
\newcommand{\ie}{i.e.\ }
\newcommand{\etal}{et~al.}

\newcommand{\I}{\mathcal{I}}
\newcommand{\tI}{\mathcal{\tilde{I}}}
\newcommand{\w}{\omega}
\newcommand{\aw}{\alpha_\omega}

\newcommand{\baseline}{\textsc{BASELINE}}
\newcommand{\pils}{\textsc{CHILS}}
\newcommand{\chils}{\textsc{CHILS}}
\newcommand{\pilscore}{\textsc{D-Core}}
\newcommand{\chilscore}{\textsc{D-Core}}
\newcommand{\lnr}{\textsc{LearnAndReduce}}

\newcommand{\hils}{\textsc{HILS}}
\newcommand{\metamis}{\textsc{METAMIS}}
\newcommand{\ilp}{\textsc{BSA}}
\newcommand{\htwis}{\textsc{HtWIS}}
\newcommand{\mmwiss}{\textsc{m$^2$wis + s}}

\newcommand{\reductiondetails}[3]{
\begin{tabular}{p{2.4cm}p{10.5cm}}
\\[-.5em]
    {Reduced Graph} & #1 \\
    {Offset} & #2 \\
    {Reconstruction} & #3 \\
\end{tabular}\\
}

\newcommand{\convertDecimal}[1]{%
  \StrSubstitute{#1}{,}{.}[\temp]%
}

\newcommand{\nptwo}[1]{%
\nprounddigits{2}%
\numprint{#1}%
}
\newcommand{\checkSmallTime}[1]{
    \convertDecimal{#1}%
    \pgfmathsetmacro{\result}{\temp} 
    \pgfmathparse{\result<0.01 ? 1 : 0} 
    \ifnum\pgfmathresult=1%
        $<\,$\numprint{0.01}%
    \else%
        \nptwo{\result}%
    \fi%
}

\newcommand{\niceGraphName}[1]{%
    \begingroup%
    \edef\tempgraph{#1}%
    \StrSubstitute{\tempgraph}{-uniform}{}[\tempgraph]%
    \StrSubstitute{\tempgraph}{AM}{}[\tempgraph]%
    \StrSubstitute{\tempgraph}{mesh-}{}[\tempgraph]%
    \StrSubstitute{\tempgraph}{snap-}{}[\tempgraph]%
    \StrSubstitute{\tempgraph}{ssmc-}{}[\tempgraph]%
    \StrSubstitute{\tempgraph}{osm-}{}[\tempgraph]%
    \StrSubstitute{\tempgraph}{fe-}{}[\tempgraph]%
    \StrSubstitute{\tempgraph}{LiveJournal1}{LiveJ.}[\tempgraph]%
    \StrSubstitute{\tempgraph}{de-island}{de-i.}[\tempgraph]%
    \StrSubstitute{\tempgraph}{district-of-columbia}{d.-of-c.}[\tempgraph]%
    \StrSubstitute{\tempgraph}{pokec-relationships}{p.-rel}[\tempgraph]%
    \tempgraph%
    \endgroup%
}

\newcommand{\CheckMarkGraphsParameter}[1]{%
    \def\result{false}
    \edef\tempgraph{#1}
    \IfSubStr{\tempgraph}{pwt}{\def\result{true}}{}%
    \IfSubStr{\tempgraph}{rotor}{\def\result{true}}{}%
    \IfSubStr{\tempgraph}{hawaii-AM3}{\def\result{true}}{}%
    \IfSubStr{\tempgraph}{vermont-AM3}{\def\result{true}}{}%
    \IfSubStr{\tempgraph}{as-skitter}{\def\result{true}}{}%
    \IfSubStr{\tempgraph}{LiveJournal}{\def\result{true}}{}%
    \IfSubStr{\tempgraph}{CR-S-L-1}{\def\result{true}}{}%
    \IfSubStr{\tempgraph}{CW-T-C-2}{\def\result{true}}{}%
    \IfSubStr{\tempgraph}{CW-S-L-4}{\def\result{true}}{}%
    \IfSubStr{\tempgraph}{MR-D-FN}{\def\result{true}}{}%
    \IfSubStr{\tempgraph}{MW-D-40}{\def\result{true}}{}%
    \IfSubStr{\tempgraph}{MW-W-05}{\def\result{true}}{}%
}

\newcommand{\CheckMarkGraphsReducedEasily}[1]{%
    \def\result{false}
    \edef\tempgraph{#1}
    \IfSubStr{\tempgraph}{venus}{\def\result{true}}{}%
    \IfSubStr{\tempgraph}{california-AM1}{\def\result{true}}{}%
    \IfSubStr{\tempgraph}{canada-AM2}{\def\result{true}}{}%
    \IfSubStr{\tempgraph}{connecticut-AM3}{\def\result{true}}{}%
    \IfSubStr{\tempgraph}{florida-AM2}{\def\result{true}}{}%
    \IfSubStr{\tempgraph}{greenland-AM1}{\def\result{true}}{}%
    \IfSubStr{\tempgraph}{hawaii-AM1}{\def\result{true}}{}%
    \IfSubStr{\tempgraph}{kansas-AM2}{\def\result{true}}{}%
    \IfSubStr{\tempgraph}{cow}{\def\result{true}}{}%
    \IfSubStr{\tempgraph}{face}{\def\result{true}}{}%
    \IfSubStr{\tempgraph}{kansas-AM1}{\def\result{true}}{}%
    \IfSubStr{\tempgraph}{kentucky-AM1}{\def\result{true}}{}%
    \IfSubStr{\tempgraph}{louisiana-AM2}{\def\result{true}}{}%
    \IfSubStr{\tempgraph}{massachusetts-AM1}{\def\result{true}}{}%
    \IfSubStr{\tempgraph}{mexico-AM1}{\def\result{true}}{}%
    \IfSubStr{\tempgraph}{michigan}{\def\result{true}}{}%
    \IfSubStr{\tempgraph}{hampshire-AM2}{\def\result{true}}{}%
    \IfSubStr{\tempgraph}{new-york-AM1}{\def\result{true}}{}%
    \IfSubStr{\tempgraph}{north-carolina-AM2}{\def\result{true}}{}%
    \IfSubStr{\tempgraph}{oregon-AM2}{\def\result{true}}{}%
    \IfSubStr{\tempgraph}{rhode-island-AM1}{\def\result{true}}{}%
    \IfSubStr{\tempgraph}{south-carolina}{\def\result{true}}{}%
    \IfSubStr{\tempgraph}{utha-AM2}{\def\result{true}}{}%
    \IfSubStr{\tempgraph}{washington-AM1}{\def\result{true}}{}%
    \IfSubStr{\tempgraph}{west-virginia-AM2}{\def\result{true}}{}%
    \IfSubStr{\tempgraph}{wisconsin}{\def\result{true}}{}%
    \IfSubStr{\tempgraph}{youtube}{\def\result{true}}{}%
    \IfSubStr{\tempgraph}{wiki-Vote}{\def\result{true}}{}%
}

\newcommand{\sampleGraphs}[1]{%
    \def\result{false}%
    \edef\tempgraph{#1}%
    \IfSubStr{\tempgraph}{body}{\def\result{true}}{}%
    \IfSubStr{\tempgraph}{pwt}{\def\result{true}}{}%
    \IfSubStr{\tempgraph}{rotor}{\def\result{true}}{}%
    \IfSubStr{\tempgraph}{ocean}{\def\result{true}}{}%
    \IfSubStr{\tempgraph}{buddha}{\def\result{true}}{}%
    \IfSubStr{\tempgraph}{ecat}{\def\result{true}}{}%
    \IfSubStr{\tempgraph}{turtle}{\def\result{true}}{}%
    \IfSubStr{\tempgraph}{dragonsub}{\def\result{true}}{}%
    \IfSubStr{\tempgraph}{columbia-AM3}{\def\result{true}}{}%
    \IfSubStr{\tempgraph}{hawaii-AM3}{\def\result{true}}{}%
    \IfSubStr{\tempgraph}{kentucky-AM3}{\def\result{true}}{}%
    \IfSubStr{\tempgraph}{rhode-island-AM3}{\def\result{true}}{}%
    \IfSubStr{\tempgraph}{as-skitter}{\def\result{true}}{}%
    \IfSubStr{\tempgraph}{LiveJournal}{\def\result{true}}{}%
    \IfSubStr{\tempgraph}{soc-pokec}{\def\result{true}}{}%
    \IfSubStr{\tempgraph}{roadNet-CA}{\def\result{true}}{}%
    \IfSubStr{\tempgraph}{ca2010}{\def\result{true}}{}%
    \IfSubStr{\tempgraph}{il2010}{\def\result{true}}{}%
    \IfSubStr{\tempgraph}{fl2010}{\def\result{true}}{}%
    \IfSubStr{\tempgraph}{ga2010}{\def\result{true}}{}%
}

\usepackage{pgfplots}
\usepackage{pgfplotstable}
\pgfplotsset{compat=1.18}
\pgfdeclarelayer{foreground}
\usepgfplotslibrary{colorbrewer}
\pgfsetlayers{main,foreground}

\usepackage{tikz} 
\usetikzlibrary{arrows.meta, positioning, tikzmark, calc}

\tikzset{
	graph/.style={draw, fill=remainingGColor, rounded corners=4.5mm, inner sep=2.5mm, align=center, minimum width=2.5cm},
	node/.style={circle, draw, fill=black, inner sep=0pt, minimum size=6pt},
	lnode/.style={node, label=left:#1},
	rnode/.style={node, label=right:#1},
	tnode/.style={node, label=above:#1},
	bnode/.style={node, label=below:#1},
	nodeR/.style={ node, draw=lightgray, fill=lightgray},
	nodeE/.style={ node, draw=redgray, fill=redgray, fill opacity=0.6, draw opacity=0.6},
	nodeI/.style={ node, draw=greengray, fill=greengray, fill opacity=0.6, draw opacity=0.6},
	weight/.style={node, rectangle, fill opacity=0, draw opacity=0, text opacity = 1,text=weightColor, yshift=0.4cm, font=\small},
	lweight/.style={node, rectangle, fill opacity=0, draw opacity=0, text opacity = 1,text=weightColor, xshift=-0.5cm, font=\small},
	rweight/.style={node, rectangle, fill opacity=0, draw opacity=0, text opacity = 1,text=weightColor, xshift=0.5cm, font=\small},
	weightR/.style={node, rectangle, fill opacity=0, draw opacity=0, text opacity = 1,text=redgray, yshift=0.5cm, font=\small},
	weightText/.style={node, rectangle, fill opacity=1, text opacity = 1,text=weightColor, fill=white, draw=white , font=\small},
	edge/.style={draw=black, fill=black, thick},
	edgeR/.style={ edge, draw=lightgray, fill=lightgray, dashed} 
}

\usetikzlibrary{pgfplots.colorbrewer}
\pgfplotsset{
    cycle list/.define={my marks}{
        every mark/.append style={solid,fill=\pgfkeysvalueof{/pgfplots/mark list fill}},mark=*\\
        every mark/.append style={solid,fill=\pgfkeysvalueof{/pgfplots/mark list fill}},mark=square*\\
        every mark/.append style={solid,fill=\pgfkeysvalueof{/pgfplots/mark list fill}},mark=triangle*\\
        every mark/.append style={solid,fill=\pgfkeysvalueof{/pgfplots/mark list fill}},mark=diamond*\\
        every mark/.append style={solid,fill=\pgfkeysvalueof{/pgfplots/mark list fill}},mark=halfcircle*\\
        every mark/.append style={solid,fill=\pgfkeysvalueof{/pgfplots/mark list fill}, mark options={rotate=90}, scale=1.25},mark=halfsquare right*\\
        every mark/.append style={solid,fill=\pgfkeysvalueof{/pgfplots/mark list fill}},mark=triangle*\\
        every mark/.append style={solid,fill=\pgfkeysvalueof{/pgfplots/mark list fill}},mark=halfdiamond*\\
        every mark/.append style={solid,fill=\pgfkeysvalueof{/pgfplots/mark list fill}},mark=10-pointed star\\        
    },
    perf_legend/.style={
        hide axis,
        cycle list/Dark2,
        cycle multiindex* list={
                my marks\nextlist
                Dark2\nextlist
                linestyles\nextlist
            },
        mark list fill={.!75!white},
        legend style={
                font=\footnotesize,
                at={(0.75,0.5)},
                anchor=center,
                draw=none,
                /tikz/every even column/.append style={column sep=5mm},
                legend columns=4
            }
    },
    perf/.style={
        ymin=0, ymax=1, xmax=1, xmin=0.8,
        axis lines=middle, 
        height=7cm, width=8.5cm,
        axis x line*=bottom,
        axis y line*=left,
        cycle list/Dark2,
        cycle multiindex* list={
                my marks\nextlist
                Dark2\nextlist
                linestyles\nextlist
            },
        mark list fill={.!75!white},
        tick label style={font=\footnotesize},
        x dir=reverse,
        x axis line style={-},
        y axis line style={-}
    },
	perf_max/.style={
        perf,
        ymin=0, ymax=1, 
		xmax=1,
        x axis line style={Stealth-},
        y axis line style={-}
    },
    perf_min/.style={
        ylabel={Fraction of instances},
        axis lines=middle, 
        axis x line*=bottom,
        axis y line*=left,
        cycle list/Dark2,
        cycle multiindex* list={
                my marks\nextlist
                Dark2\nextlist
                linestyles\nextlist
            },
        mark list fill={.!75!white},
        tick label style={font=\footnotesize},
        ymin=0, ymax=1,
		xmin=1,
        height=7cm, width=8.5cm,
        x axis line style={-Stealth},
        y axis line style={-},
        ylabel near ticks,
        xlabel={$\tau$},
    },
    perf_left/.style={
        perf,
        height=7cm, width=5.5cm,
        ylabel near ticks,
        ylabel={Fraction of instances}
    },
    perf_right/.style={
        perf,
        xmax=0.8, xmin=0.5,
        height=7cm, width=3cm,
        ytick={},
        yticklabels={},
        xlabel={$\tau$},
        x axis line style={-,Stealth-}
    }
}

\newcommand{\remainingG}[2]{
	\begin{pgfonlayer}{foreground}
		\node[graph] at (#1) (box) {#2};
	\end{pgfonlayer}
}
\definecolor{remainingGColor}{RGB}{226,243,253}
\definecolor{weightColor}{RGB}{204,0,0}
\definecolor{darkgreen}{rgb}{0.0, 0.5, 0.0}
\definecolor{lightgreen}{rgb}{0.235, 0.7, 0.235}
\definecolor{green}{RGB}{0, 150, 130}
\definecolor{green70}{RGB}{76, 181, 167}
\definecolor{blue}{RGB}{70, 100, 170}
\definecolor{blue70}{RGB}{125,146,195}
\definecolor{blue50}{RGB}{162,177,212}
\definecolor{blue30}{RGB}{199,208,229}
\definecolor{blue15}{RGB}{227,232,242}
\definecolor{lightgray}{rgb}{0.86,0.86,0.86}
\definecolor{greengray}{RGB}{193,228,224}
\definecolor{redgray}{RGB}{233,183,183}
\definecolor{red}{RGB}{204,0,0}
\definecolor{limegreen}{rgb}{0.6, 0.8, 0.0}
\definecolor{lipicsYellow}{rgb}{0.99,0.78,0.07}

\newtcbtheorem{obs}{Observation}{
    colback=lightgray,  
    colframe=lightgray, 
    boxsep=0pt,       
    left=1mm,         
    right=1mm,        
    top=5pt,          
    bottom=5pt,      
    fonttitle=\bfseries,
    coltitle=black,
    coltext=black, 
    terminator sign={},
    attach title to upper,
    after title={\hspace{1em}} 
}{obs}

\counterwithin*{theorem}{subsection}

\newtheoremstyle{custom}
  {3pt}{3pt} 
  {\itshape}
  {}
  {\bfseries}
  {}
  { }
  {\thmname{#1}\thmnumber{ #2}\thmnote{ \textnormal{(#3)}}}
\makeatother

\theoremstyle{custom}
\newtheorem{internalreduction}{Reduction}[subsection]

\newenvironment{reduction}[1][\unskip]{
    \begin{internalreduction}[#1]
    \noindent\trivlist\item\vspace{-.5em}\ignorespaces
}{
    \endtrivlist
    \end{internalreduction}
}
\newenvironment{figreduction}[2][\unskip]{
\begin{internalreduction}[#1. Figure~\ref{fig:#2}]\label{red:#2}
    \noindent\trivlist\item\vspace{-.5em}\ignorespaces
}{
    \endtrivlist
    \end{internalreduction}
}

\title{Accelerating Reductions Using Graph Neural Networks and a New Concurrent Local Search for the Maximum Weight Independent Set Problem} 

\titlerunning{Accelerating Reductions Using GNNs and a New Concurrent LS for MWIS}

\author{Ernestine Gro{\ss}mann}{Heidelberg University,  Faculty of Mathematics and Computer Science, Germany}{e.grossmann@informatik.uni-heidelberg.de}{https://orcid.org/0000-0002-9678-0253}{Supported by DFG grant SCHU 2567/3-1.}
\author{Kenneth Langedal\footnote{Corresponding author}}{University of Bergen, Department of Informatics, Norway}{kenneth.langedal@uib.no}{https://orcid.org/0009-0001-6838-4640}{Supported by the Research Council of Norway under contract 303404 and Meltzer Research Fund, project number 104066111.}
\author{Christian Schulz}{Heidelberg University, Faculty of Mathematics and Computer Science, Germany}{christian.schulz@informatik.uni-heidelberg.de}{https://orcid.org/0000-0002-2823-3506}{}

\authorrunning{E. Gro{\ss}mann, K. Langedal, C. Schulz} 
\Copyright{Ernestine Gro{\ss}mann, Kenneth Langedal, and Christian Schulz} 



\category{}

\relatedversion{} 

\funding{\,}

\acknowledgements{The inspiration for this work and collaboration started at the Dagstuhl Seminar 23491 on Scalable Graph Mining and Learning~\cite{koutra2024scalable}. We also thank Fredrik Manne for his helpful feedback throughout the writing process. }

\nolinenumbers

\EventEditors{John Q. Open and Joan R. Access}
\EventNoEds{2}
\EventLongTitle{42nd Conference on Very Important Topics (CVIT 2016)}
\EventShortTitle{CVIT 2016}
\EventAcronym{CVIT}
\EventYear{2016}
\EventDate{December 24--27, 2016}
\EventLocation{Little Whinging, United Kingdom}
\EventLogo{}
\SeriesVolume{42}
\ArticleNo{23}

\begin{document}

\maketitle

\begin{abstract}
    The \textsc{Maximum Weight Independent Set} problem is a fundamental \NP-hard problem in combinatorial optimization with several real-world applications. Given an undirected vertex-weighted graph, the problem is to find a subset of the vertices with the highest possible weight under the constraint that no two vertices in the set can share an edge. An important part of solving this problem in both theory and practice is data reduction rules, which several state-of-the-art algorithms rely on. However, the most complicated rules are often not used in applications since the time needed to check them exhaustively becomes infeasible. In this work, we introduce three main results. First, we introduce several new data reduction rules and evaluate their effectiveness on real-world data. Second, we use a machine learning screening algorithm to speed up the reduction phase, thereby enabling more complicated rules to be applied. Our screening algorithm consults a Graph Neural Network oracle to decide if the probability of successfully reducing the graph is sufficiently large. For this task, we provide a dataset of labeled vertices for use in supervised learning. We also present the first results for this dataset using established Graph Neural Network architectures. Third, we present a new concurrent metaheuristic called \textsc{Concurrent Difference-Core Heuristic}. On the reduced instances, we use our new metaheuristic combined with iterated local search, called \pils{} (\textsc{Concurrent Hybrid Iterated Local Search}). For this iterated local search, we provide a new implementation specifically designed to handle large graphs of varying densities. \pils{} outperforms the current state-of-the-art on all commonly used benchmark instances, especially the largest ones.
\end{abstract}

\newpage

\section{Introduction}

Consider an undirected vertex-weighted graph $G=(V,\, E,\, \w)$, where $V$ is the set of vertices, $E$ is the set of edges, and $\w : V \rightarrow \mathbb{R}^{+}$ is a function that maps each vertex to a positive weight. An \textit{independent set} $\I \subseteq V$ is a subset of the vertices such that no two members of the independent set share an edge, \ie for all $u,\, v \in \I$ it holds that $\{u,\, v\} \notin E$. The \textsc{Maximum Weight Independent Set} (MWIS) problem asks for an independent set $\I$ of maximum weight, where the weight of $\I$ is defined as the sum $\sum_{v \in \I} \w(v)$ of the vertices. MWIS is a generalization of the classical \NP-hard problem \textsc{Maximum Independent Set} (MIS), where all weights equal one. For MIS and MWIS, the complement of a solution forms a minimum (weighted) vertex cover. A \textit{vertex cover} $C \subseteq V$ is a set of vertices that cover all the edges, \ie for all $\{u,\, v\} \in E$ we have $u \in C$ or $v \in C$. For all practical purposes, the problem of finding a maximum independent set is equivalent to that of finding a minimum vertex cover. The decision version of the \textsc{Minimum Vertex Cover} (MVC) problem was one of Karp's original 21 \NP-complete problems~\cite{karp2010reducibility}.

The MWIS problem and other closely related problems have several practical applications ranging from matching molecular structures to wireless networks~\cite{butenko2003maximum}. Recently, Dong \etal~\cite{dong2021new} introduced a new collection of instances based on real-life long-haul vehicle routing problems at Amazon. The problem they want to solve is to find a subset of vehicle routes such that no two routes share a driver or a load. Each route has a weight, and the objective is to maximize the sum of the route weights. To state this problem as a MWIS, they build a conflict graph where vertices correspond to routes and the route weights are modeled by vertex weights.
Furthermore, they connect two vertices by an edge if the corresponding routes have a conflict, \ie they share a driver or a load. For a further, more detailed overview of applications, \hbox{see Butenko~\cite{butenko2003maximum}.}
 
It is well known that data reduction rules can speed up algorithms for \NP-hard problems. They reduce an instance in such a way that an optimal solution for the reduced instance can be lifted to an optimal solution for the original instance. Several reduction rules have been developed for MIS and MWIS. Additionally, reduction rules have also improved many heuristic approaches for MWIS and associated problems. Such reduction rules are often used as a preprocessing step before running an exact algorithm or heuristic on the reduced graph.

\paragraph*{Our Results}

Our contribution is two-fold: first, we developed an advanced, exact preprocessing tool called \lnr~that employs Graph Neural Networks (GNNs) to decide where to apply data reduction rules. In addition to existing rules (see \cite{gellner2021boosting,gu2021towards,lamm2019exactly,xiao2024maximum,xiao2021efficient}), we propose new data reduction rules that can reduce the instances even further. With our GNN filtering, we can now apply reduction rules that were previously unused in other works due to the computational cost required to determine their applicability~\cite{reduction_survey}. Second, we propose a new concurrent iterated local search heuristic~\pils~to compute large-weight independent sets very fast. In particular, our heuristic works by alternating between the full graph and the \textsc{Difference-Core} (\pilscore)~which is a subgraph constructed using multiple heuristic solutions. With this heuristic and \lnr~we are able to outperform existing heuristics across a wide variety of real-world instances. Our main results can be summarized as follows.

\begin{itemize}
    \item We introduce seven new data reduction rules for the MWIS problem.
    \item We present a new dataset with labeled vertices for the problem of early reduction rule screening. The dataset contains two collections of graphs, one with unreduced instances and one after running a set of fast reductions. The second is the one we use for the early reduction rule screening, and it consists of more than one million labeled vertices. 
    \item Using \lnr, we can reduce our instances to within one percent of what is possible using the full set of reduction rules while spending less than 23\,\% of the corresponding time. These results are adjusted for fast reduction rules that would \hbox{always be applied.}
    \item For a large set of real-world test instances accumulated over several years by earlier works, \pils~combined with \lnr~finds the best solution across \emph{all} test instances. This result is kept separate from the new vehicle routing instances~\cite{dong2021new} due to the significant differences in size and difficulty.
    \item On vehicle routing instances, we compare \pils~with two recent heuristics, METAMIS~\cite{metamis} and a recent Bregman-Sinkhorn algorithm~\cite{haller2024bregman}, both designed specifically for these instances. Even though \pils~is not optimized for vehicle routing instances and, in contrast to the competitors, does not make use of the additional clique information provided for these instances, it still finds the best solution on 31/37 instances while being significantly closer to the best competitor in the cases where \pils~is not best.
    \item Running \chils~in parallel significantly improves performance to the point where \chils~computes the best solution on 35/37 instances. These results are with one-hour time limit and the same 16-core CPU used to evaluate all the heuristics. For parallel scalability, we include experiemnts on a 128-core CPU where the parallel version of \chils{} reaches speedups of up to 104 compared to the sequential version.
\end{itemize}

    \section{Preliminaries}
    In this work, a graph $G=(V,E,\w)$ is an undirected vertex-weighted graph with ${n=|V|}$ and ${m = |E|}$, where ${V =\{0,\ldots,n-1\}}$ and $\w : V \rightarrow \mathbb{R}^{+}$. The neighborhood $N(v)$ of a vertex $v \in V$ is defined as $N(v) = \{u \in V \mid \{u,v\} \in E\}$. 
    Additionally, we define $N[v]=N(v) \cup \{v\}$. The same sets are defined for the neighborhood $N(U)$ of a set of vertices ${U \subset V}$, \ie ${N(U) = \cup_{v \in U} N(v)\setminus U}$ and $N[U] = N(U) \cup U$. The degree $\mathrm{deg}(v)$ of a vertex is defined as the number of its neighbors $\mathrm{deg}(v)=|N(v)|$. The complement graph is defined as ${\overline{G}=(V,\overline{E})}$, where ${\overline{E}=\{\{u,v\} \mid \{u,v\} \notin E\}}$ is the set of edges not present in $G$.
    Furthermore, for a graph $G=(V,E)$ we define an induced subgraph $G[S]$  on a subset of vertices $S \subset V$ by $G[S] = (S,\{\{u,v\} \in E \mid u,v\in S\})$. For $S \subset V$, we use the notation $G-S$ instead of $G[V\setminus S]$. Similarly for a single vertex $v\in V$, we abbreviate $G[V\setminus\{v\}]$ to $G-v$.
    A set $I\subseteq V$ is called an \textit{independent set} (IS) if for all vertices $u,v \in I$ there is no edge $\{u,v\} \in E$. For a given IS~$\I$, a vertex $v \notin \I$ is called free if $\I \cup \{v\}$ is still an independent set. An IS is called \textit{maximal} if there are no free vertices. 
    The \textsc{Maximum Independent Set} problem (MIS) is that of finding an IS with maximum cardinality.
    Similarly, the \textsc{Maximum Weight Independent Set} problem (MWIS) is that of finding an IS with maximum weight. The weight of an independent set $\I$ is defined as $\omega(\I) = \sum_{v \in \I}\omega(v)$ and $\alpha_\omega(G)$ denotes the weight of an MWIS of $G$. The complement of a maximal independent set is a \textit{vertex cover}, \ie a subset ${C \subseteq V}$ such that every edge $e \in E$ is covered by at least one vertex $v \in C$.
    An edge is \textit{covered} if it is incident to one vertex in the set~$C$. The \textsc{Minimum Vertex Cover} problem, defined as computing a vertex cover with minimum cardinality, is thereby dual to the MIS problem. Another closely related concept is cliques. A \textit{clique} is a set $Q \subseteq V$ such that all vertices are pairwise adjacent. A clique in the complement graph $\overline{G}$ corresponds to an independent set in the original graph~$G$. A vertex is called \textit{simplicial} when its neighborhood forms a clique.

\section{Related Work}\label{sec:related_work}

We give an overview of latest work on both exact and heuristic procedures for MWIS. 
For a full overview of the related work on MWIS, MWVC, and \textsc{Maximum Weighted Clique} solvers, as well as all known reduction rules for the MWIS and MWVC problems, we defer to the recent survey~\cite{reduction_survey}.
For more general details on data reduction techniques, we refer the reader to the survey~\cite{Abu-Khzam2022}.
\subsection{Exact Methods}
Exact algorithms usually compute optimal solutions by systematically exploring the solution space. \emph{Branch-and-bound} is a frequently used paradigm in exact algorithms for combinatorial optimization problems~\cite{ostergaard2002fast,warren2006combinatorial}. In the case of MWIS, branch-and-bound algorithms compute optimal solutions by case distinctions where vertices are either included in the current solution or excluded from it, branching into two subproblems.
Over the years, branch-and-bound methods have been improved by new branching schemes or better pruning methods using upper and lower bounds to exclude specific subtrees~\cite{babel1994fast,balas1986finding, li2017minimization}. In particular, Warren and Hicks~\cite{warren2006combinatorial} proposed three branch-and-bound algorithms that combine the use of weighted clique covers and a branching scheme first introduced by Balas and Yu~\cite{balas1986finding}. Their first approach extends the algorithm by Babel~\cite{babel1994fast} by using more intricate data structures to improve performance. The second one is an adaptation of the algorithm of Balas and Yu, which uses a weighted clique heuristic that yields structurally similar results to the algorithm of Balas and Yu. The last algorithm is a hybrid version that combines both algorithms and can compute optimal solutions on graphs with hundreds of vertices.

In recent years, reduction rules have frequently been added to branch-and-bound methods yielding so-called \emph{branch-and-reduce} algorithms~\cite{akiba-tcs-2016}. Branch-and-reduce algorithms often outperform branch-and-bound algorithms by applying reduction rules between each branching step. For the unweighted case, many branch-and-reduce algorithms have been developed in the past. The currently best exact solver~\cite{hespe2020wegotyoucovered}, which won the PACE challenge 2019~\cite{hespe2020wegotyoucovered, bogdan-pace}, uses a portfolio of branch-and-reduce/bound solvers for the complementary problems. For non-portfolio solvers, Plachetta~\etal~\cite{plachetta2021sat} improved on the branch-and-reduce approach by using SAT solvers for additional pruning. Recently, novel branching strategies have been presented by Hespe~\etal~\cite{hespe2021targeted} and later enhanced by Langedal~\etal~\cite{langedal2024targeted} to improve both branch-and-bound and branch-and-reduce approaches~further.

The first branch-and-reduce algorithm for the weighted case was presented by Lamm \etal~\cite{lamm2019exactly}. The authors first introduce two meta-reductions called neighborhood removal and neighborhood folding, from which they derive a new set of weighted reduction rules. On this foundation, a branch-and-reduce algorithm was developed using pruning with weighted clique covers similar to the approach by Warren and Hicks~\cite{warren2006combinatorial} for upper bounds and an adapted version of the \textsc{ARW} local search~\cite{andrade-2012} \hbox{for lower bounds.} The exact algorithm by Lamm~\etal~was then extended by Gellner~\etal~\cite{gellner2021boosting} to utilize different variants of a transformation called \emph{struction}, originally introduced by Ebenegger~\etal~\cite{ebenegger1984pseudo} and later improved by Alexe~\etal~\cite{alexe2003struction}. In contrast to previous reduction rules, struction transformations do not necessarily decrease the graph size but rather transform the graph, which later can lead to even further reduction. These variants were integrated into the framework of Lamm~\etal~in the preprocessing as well as in the reduction step. The experimental evaluation shows that this algorithm can solve a large set of real-world instances and outperforms the original branch-and-reduce algorithm by Lamm~\etal~as well as different state-of-the-art heuristic approaches such as the heuristic {\hils} by Nogueira~\cite{hybrid-ils-2018} and two other local search algorithms \textsc{DynWVC1} and \textsc{DynWVC2} by Cai~\etal~\cite{cai-dynwvc}. Furthermore, a new branch-and-reduce algorithm was recently introduced using two \hbox{new reduction rules~\cite{DBLP:conf/icde/ZhengGPY20}.}

Recently, Xiao~\etal~\cite{xiao2024maximum} also presented a branch-and-bound algorithm using reduction rules working especially well on sparse graphs. They undertake a detailed analysis of the running time bound on special graphs in their theoretical work. With the measure-and-conquer technique, they can show that the running time of their algorithm is~$\mathcal{O}^*(1.1443^{(0.624x-0.872)n})$ where $x$ is the average degree of the graph. This improves previous time bounds using polynomial space complexity for graphs of average degree \hbox{up to three.}

Figiel~\etal~\cite{DBLP:conf/esa/FigielFNN22}~introduced a new idea to the state-of-the-art way of applying reductions. They propose not only to perform reductions but also the possibility of undoing them during the reduction process. They showed for the unweighted MVC problem that this can lead to new possibilities for applying further reductions and lead to \hbox{smaller reduced graphs.}

Finally, there are exact procedures that are either based on other extensions of the branch-and-bound paradigm~\cite{rebennack2011branch,warrier2007branch,warrier2005branch}, or on the reformulation into other \NP-hard problems, for which a variety of solvers already exist. For instance, Xu~\etal~\cite{xu2016new} developed an algorithm called \textsc{SBMS} that calculates an optimal solution for a given MWVC instance by solving a series of SAT instances. Also, for the MWVC problem, a new exact algorithm using the branch-and-bound idea combined with data reduction rules was recently presented~\cite{DBLP:journals/corr/abs-1903-05948}. We also note that several recent works on the complementary maximum weighted clique problem can handle large real-world networks~\cite{fang2016exact,held2012maximum,jiang2017exact}. However, using these solvers for the MWIS problem requires computing the complement graph. Since large real-world networks are often very sparse, processing their complements quickly becomes infeasible due to \hbox{their memory requirement.}

\subsection{Heuristic Methods}

Local search is a widely used heuristic approach for MWIS. This starts from any feasible solution and then tries to improve it by simple insertion, removal, or swap operations. Local search generally offers no theoretical guarantees for the solution quality. However, in practice, it often finds high-quality solutions significantly faster than exact procedures. For unweighted graphs, the iterated local search (\textsc{ARW}) by Andrade~\etal~\cite{andrade-2012} is a very successful heuristic. It is based on so-called $(1,2)$-swaps that removes one vertex from the solution and adds two new vertices, thus improving the current solution by one. Their heuristic uses special data structures that find such a $(1,2)$-swap in linear time in the number of edges or prove that none exists. Their heuristic is able to find (near-)optimal solutions for small to medium-size instances in milliseconds but struggles on large sparse instances with millions of vertices.

The hybrid iterated local search (\hils) by Nogueira~\etal~\cite{hybrid-ils-2018} adapts the \textsc{ARW} algorithm for weighted graphs. In addition to weighted $(1,2)$-swaps, it also uses $(\omega,1)$-swaps that adds one vertex $v$ into the current solution and excludes its neighbors. These two types of neighborhoods are explored separately using variable neighborhood descent (VND). When it was introduced, their algorithm found optimal solutions on well-known benchmark instances within milliseconds and outperformed other state-of-the-art local search heuristics.

Two other local search heuristics, \textsc{DynWVC1} and \textsc{DynWVC2}, for the complementary minimum weight vertex cover (MWVC) problem were presented by Cai~\etal~\cite{cai-dynwvc}. Their heuristic extend the existing \textsc{FastWVC} heuristic~\cite{li2017efficient} by dynamic selection strategies for vertices to be removed from the current solution. In practice, \textsc{DynWVC1} outperforms previous MWVC heuristics on map labeling instances and large-scale networks. \textsc{DynWVC2} provides further improvements on large-scale networks but performs worse on \hbox{map labeling instances.}

Li~\etal~\cite{DBLP:journals/jors/LiHCGWY20} presented a local search heuristic \textsc{NuMWVC} for the MWVC problem. Their heuristic applies reduction rules during the construction phase of the initial solution. Furthermore, they adapt the configuration checking approach~\cite{cai2011local} to the MWVC problem, which is used to reduce cycling, \ie returning to a solution that has been visited recently. Finally, they develop a technique called self-adaptive vertex-removing, which dynamically adjusts the number of removed vertices per iteration. Experiments showed that \textsc{NuMWVC} outperformed state-of-the-art approaches on graphs of up to millions of vertices.

A hybrid method was recently introduced by Langedal~\etal~\cite{langedal2022efficient} to solve the MWVC problem. For this approach, they combined elements from exact methods with local search, data reductions, and graph neural networks. In their experiments, they achieve clear improvements compared to \textsc{DynWVC2}, {\hils}, and \textsc{NuMWVC} in both solution quality and running time.

With \textsc{EvoMIS}, Lamm~\etal~\cite{lamm2015graph} presented an evolutionary approach to tackle the maximum independent set problem. The key feature of their heuristic was to use graph partitioning to come up with natural combine operations, where whole blocks of solutions can be exchanged easily. Local search algorithms were added to these combine operations to improve the solutions further.

Combining the branch-and-reduce approach with the evolutionary algorithm \textsc{EvoMIS}, a reduction evolution algorithm \textsc{ReduMIS} was presented by Lamm~\etal~\cite{redumis-2017}. In their experiments, \textsc{ReduMIS} outperformed the local search \textsc{ARW} as well as the pure evolutionary approach \textsc{EvoMIS}. Another reduction-based heuristic {\htwis} was presented recently by Gu~\etal~\cite{gu2021towards}. They repeatedly apply reductions exhaustively and then choose one vertex by a tie-breaking policy to add to the solution. Once this vertex and its neighbors has been removed from the graph, the reduction rules can be applied again. Their experiments show a significant improvement in running time. Recently, Gro{\ss}mann~\etal~\cite{mmwis} introduced a heuristic combining data reduction rules with an evolutionary approach. Here, the authors make use of exact data reductions and heuristic reductions derived from the population to reduce the graph iteratively. With this technique, they are able to achieve near-optimal solutions for a wide set of instances.

A new metaheuristic was introduced by Dong~\etal~\cite{metamis}, particularly for the vehicle routing instances introduced by Dong~\etal~\cite{dong2021new}. With their metaheuristic \textsc{METAMIS}, they developed a local search heuristic using a new variant of path-relinking to escape local optima. In their experiments, they outperform the {\hils} heuristic on a wide range of instances, both in terms of time and solution quality. The vehicle routing instances come equipped with initial warm start solutions and clique information derived from the application. Using this clique information, Haller and Savchynskyy~\cite{haller2024bregman} proposed a Bregman-Sinkhorn algorithm that addresses a family of clique cover LP relaxations. In addition to solving the relaxed dual problem, they utilize a simple and efficient primal heuristic to obtain feasible integer solutions for the initial non-relaxed problem. In their experiments, they outperform \textsc{METAMIS} on time and solution quality, but only in the cold-start configuration where \textsc{METAMIS} does not use the precomputed solutions.

\section{Learn and Reduce - GNN Guided Preprocessing}

In this section, we introduce~\lnr, our new exact preprocessing technique that reduces the input graph efficiently using data reductions and Graph Neural Networks~(GNNs).

In Section \ref{sec:reduction-rules}, we introduce several new data reduction rules for the MWIS problem. Some new reductions and existing ones need to solve additional MWIS problems on subgraphs, meaning preprocessing using these rules could take exponential time. For practical use, the more costly reductions are limited in some way, especially for large instances. This is done, for example, by limiting the degree of the vertex to apply the reduction on, the subgraph size, or the time spent on each vertex. In Section~\ref{sec:gnn-models}, we introduce a new GNN application in the form of early data-reduction screening. We also provide a new and openly available supervised-learning dataset with graphs and labeled vertices for each reduction rule. At inference, for each expensive data reduction rule, a pre-trained GNN model decides what vertices should be checked for the applicability of the rule.

\subsection{New Data Reduction Rules}
\label{sec:reduction-rules}
This section introduces our new data reduction rules for the MWIS problem following the scheme used in~\cite{reduction_survey}. For every reduction rule, we first describe the pattern that can be reduced by the corresponding rule. Afterwards, the information for the construction of the \textit{reduced graph}, called $G'$ is given. Then, the \textit{offset} describes the difference between the weight of an MWIS on the reduced graph $\aw(G')$ and the weight of an MWIS on the original graph $\aw(G)$. Lastly, the information on how the solution on the reduced instance, called $\I'$, can be reconstructed to a solution on the original graph, called $\I$, is provided.
In addition to including or excluding vertices directly, some reduction rules combine multiple vertices into potentially new vertices. We call this combine procedure \textit{folding}. Including or excluding the folded vertices in $\I$ only depends on whether the vertices they are folded into are \hbox{included in $\I'$.}
Since our new reduction rules are extensions of or proven by other reduction rules, these are presented first. 
\begin{reduction}[Simplicial Vertex by Lamm~\etal~\cite{lamm2019exactly}] \label{red:clique}
    Let $v\in V$ be simplicial with maximum weight in its neighborhood, \ie $\w(v)\geq \max\{\w(u) \mid u \in N(v)\}$, then include $v$.\\
    \reductiondetails{$G'=G-N[v]$}{$\aw(G) = \aw(G')+ \w(v)$}{$\I = \I' \cup \{v\}$}
\end{reduction}

\begin{reduction}[Domination by Lamm~\etal~\cite{lamm2019exactly}]\label{red:dom}
	Let $u,v\in V$ be adjacent vertices such that $N[u]\subseteq N[v]$. If $\w(v)\leq \w(u)$, then exclude $v$.\\
    \reductiondetails{$G'=G-v$}{$\aw(G) = \aw(G')$}{$\I = \I'$} 
\end{reduction}
\begin{reduction}[Twin by Lamm~\etal~\cite{lamm2019exactly}]\label{red:twin}
	Let $u, v\in V$ have equal, independent neighborhoods $N(u) = N(v) = \{p,q,r\}$.
	\begin{itemize}
		\item If $\w(\{u,v\}) \geq \w(\{p,\,q,\,r\})$, then include $u$ and $v$.\\
        \reductiondetails{$G' = G-N[\{u,v\}]$}{$\aw(G)=\aw(G') +\w(u) + \w(v)$}{$\I = \I' \cup \{u,\,v\}$        \\[-.5em]}
		\item If $\w(\{u,v\}) < \w(\{p,q,r\})$ but $\w(\{u,v\}) > \w(\{p,q,r\}) - \min\{\w(x) \mid x\in\{p,q,r\}\}$, then fold $u,\, v,\, p,\, q,\, r$ into a new vertex $v'$. \\
        \reductiondetails{ $G' = G[(V\cup \{v'\})\setminus(N[v]\cup \{u\})]$ with $N(v')=N(\{p,\,q,\,r\})$ and \vfill \break $\w(v') = \w(\{p,\,q,\,r\}) - \w(\{u,\,v\})$.}{$\aw(G) = \aw(G') + \w(\{u,v\})$}{If $v'\in \I'$, then $\I = (\I'\setminus \{v'\})\cup \{p,\,q,\,r\}$, else $\I = \I' \cup \{u,\,v\}$}
	\end{itemize}
\end{reduction}

\begin{reduction}[Heavy Set by Xiao~\etal~\cite{xiao2021efficient}]\label{red:heavy_set}
	Let $u$ and $v$ be non-adjacent vertices having at least one
	common neighbor $x$. If for every independent set $\tI$ in the induced subgraph $G[N(\{u,\, v\})]$, $\w(N(\tI) \cap \{u,\,v\})\geq \w(\tI)$, then include $u$ and $v$.\\    
    \reductiondetails{$G'=G-N[\{u,\, v\}]$}{$\aw(G) = \aw(G')+ \w(u) + \w(v)$}{$\I = \I' \cup \{u,\,v\}$}
\end{reduction}

\paragraph*{Domination Based Reductions}
We derive the following two reduction rules by extending Reduction~\ref{red:dom} first introduced by Lamm~\etal~\cite{lamm2019exactly}.
The idea of the original domination rule is to find two adjacent vertices $u,v\in V$ where any independent set including $v$ can be transformed into an equal or higher weight independent set by replacing $v$ with $u$.

In the extended domination rule, see Reduction~\ref{red:e_dom}, instead of removing vertices, we remove edges and reduce vertex weights. With this reduction rule, we allow two adjacent vertices $u$ and $v$ to be both in the solution on the reduced instance. If the edge removal leads to a solution including both vertices, we remove the previously dominated, lower weight vertex from the solution in the restoring process. With this reduction rule, we sparsify the graph and potentially make other rules applicable.

\begin{reduction}[Extended Domination]\label{red:e_dom}
Let $u, v\in V$ be adjacent vertices such that $N[u] \subseteq N[v]$. If $\w(v) > \w(u)$, then we remove the edge between $u$ and $v$.\\
\reductiondetails{$G'=(V, E \setminus \{u, v\})$ and $\w(v) = \w(v) - \w(u)$}{$\aw(G) = \aw(G')$}{If $v \in \I'$, then $\I = \I'\setminus \{u\}$, else $\I = \I'$}
\end{reduction}
\begin{proof}
    Let $u, v\in V$ be adjacent vertices such that $N[u] \subseteq N[v]$ and $\w(v) > \w(u)$. Further, let $G'$ be the reduced graph after applying the reduction to $u$ and $v$. We use $u'$ and $v'$ to refer to $u$ and $v$ in $G'$. First, we consider the case that $v' \in \I'$ and show that $\I=\I'\setminus \{u'\}$ is an MWIS in $G$. Since $N[u]\subseteq N[v]$ and $v' \in \I'$ it holds that $u' \in \I'$, because $\I'$ is maximal. Now, assume there is an MWIS $\tI$ in $G$ with $\w(\tI) > \w(\I)$. Then, the MWIS $\tI$ is also an independent set in the reduced instance and $\w(\tI) > \w(\I' \setminus \{u'\}) + \w(u) = \w(\I' \setminus \{u', v'\}) + \w(v) - \w(u) + \w(u) = \w(\I' \setminus \{u', v'\}) + \w(v') + \w(u') = \w(\I')$. This contradicts that $\I'$ is an MWIS in $G'$. Furthermore, since $\I'$ is an independent set and we only removed $u$ to obtain $\I$, the independent set property still holds and therefore $\I$ is an MWIS in $G$. 

    We prove the second case $v' \notin \I'$ in the same way. Here, $\I'$ is also an independent set in $G$. Since $v' \notin \I'$ it holds that the weight of $\I'$ is the same in $G$ and $G'$. Therefore, there can not exist an independent set of higher weight than $\I'$ in $G$, since this would also exist in the reduced instance contradicting $\I'$ being an MWIS in $G'$.
\end{proof}

Additionally, we add a reduction rule which reverses Reduction~\ref{red:e_dom}. This way, we can add back edges that were previously removed and introduce new edges later in \hbox{the reduction process.}

\begin{reduction}[Extended Domination Reversed]\label{red:e_dom_rev}
Let $u, v\in V$ be non-adjacent vertices such that $N(u) \subseteq N(v)$. If $\w(u) + \w(v) < \w(N(v))$, then we can add an edge between $u$ and $v$.\\
\reductiondetails{$G'=(V, E \cup \{u, v\})$ and $\w(v) = \w(v) + \w(u)$}{$\aw(G) = \aw(G')$}{If $v \in \I'$, then $\I = \I' \cup \{u\}$, else $\I = \I'$}
\end{reduction}

\begin{proof}
    Let $u, v\in V$ be non-adjacent vertices such that $N(u) \subseteq N(v)$ and $\w(u) + \w(v) < \w(N(v))$. Further, let $G'$ be the reduced graph after applying the reduction to $u$ and $v$. We use $u'$ and $v'$ to refer to $u$ and $v$ in $G'$. Note that since $\w(v') = \w(v) +\w(u)$ the weight $\w(\I) = \w(\I')$. First, we consider the case that $v' \in \I'$. We show that $\I = \I' \cup \{u'\}$ is an MWIS in $G$. Since $v' \in \I'$ and $u'$ and $v'$ are adjacent in $G'$ it holds that $u' \notin \I'$. Now, assume there is a MWIS $\tI$ in $G$ with $\w(\tI) > \w(\I)$. We only consider the case where $u,v \in \tI$, since otherwise, either $\tI$ is not maximal, or $\I$, $\I'$ and $\tI$ are all equal. With $u,v \in \tI$ we construct an independent set $\tI' = \tI \setminus \{u, v\} \cup \{v'\}$ in $G'$. It holds that $\w(\tI') = \w(\tI)$ since $\w(v') = \w(v) + \w(u)$. Since $\I'$ is an MWIS in $G'$, it follows that $\w(\I') \geq \w(\tI')$. Now, $\w(\I) < \w(\tI) = \w(\tI') \leq \w(\I')$ which contradicts that $\w(\I) = \w(\I')$. For the case of $v \notin \I'$ it holds that the weights and sets are directly equivalent between $G$ and $G'$, and therefore $\I=\I'$ is an MWIS in $G$.
\end{proof}

\begin{remark}
    Initially, the concept of adding edges as described in Reduction~\ref{red:e_dom_rev} may seem counter-intuitive since the goal is to reduce the graph size. However, this approach is particularly valuable for more complex reduction rules that involve solving independent sets within neighborhoods. By adding edges, additional vertices are incorporated into the direct neighborhood. This expansion can potentially enable further applications of other reduction rules, such as Reduction~\ref{red:heavy_set} or~\ref{red:heavy_set3}.
\end{remark}

\paragraph*{Twin Based Reductions}
    Reduction~\ref{red:twin} can be generalized to also be applied if the neighborhood is larger than three vertices or when it is not an independent set. The main idea is that for any two vertices that have an equal neighborhood any maximal solution always includes either both or none of these vertices. Therefore, we can fold the two twin vertices. In some special cases we can also fold these twin vertices \hbox{with their neighborhood.}
\begin{reduction}[Extended Twin]\label{red:e_twin}
	Let $u, v\in V$ be non-adjacent and with equal neighborhoods $N(u) = N(v)$. Let further $\I_{N(v)}$ be an MWIS on $G[N(v)]$.
	\begin{itemize}
		\item If $\w(u)+\w(v) \geq \w(\I_{N(v)})$, then include $u$ and $v$.\\
        \reductiondetails{$G' = G-N[\{u,v\}]$}{$\aw(G) = \aw(G') + \w(u) + \w(v)$}{$\I = \I' \cup \{u, v\}$\\[-.5em]}
		\item If $\w(u) + \w(v) < \w(\I_{N(v)})$ but $\I_{N(v)}$ is the only independent set in $N(v)$ with this property, then fold $u$, $v$, and $N(v)$ into $v'$.\\
        \reductiondetails{$G' = G[(V \cup \{v'\})\setminus N[\{u,v\}]]$ with $N(v') = N(N[\{u,v\}])$ \\&and ${\w(v') = \w(\I_{N(v)}) - \w(u) - \w(v)}$}{$\aw(G) = \aw(G') + \w(u) + \w(v)$}{If $v' \in \I'$, then $\I = (\I' \setminus \{v'\})\cup \I_{N(v)}$, else $\I = \I' \cup \{u, v\}$\\[-.5em]}
		\item Otherwise, fold $u$ and $v$ into a new vertex $v'$.\\
        \reductiondetails{$G' = G[(V\cup \{v'\}) \setminus \{u,v\} ]$ with $N(v') = N(\{u,\,v\})$ and $\w(v') = \w(v) + \w(u)$}{$\aw(G) = \aw(G')$}{If $v' \in \I'$, then $\I = \I' \setminus \{v'\} \cup \{u,\,v\}$, else $\I = \I'$}
    \end{itemize}
\end{reduction}
\begin{proof}
     Let $u, v \in V$ be twin vertices, \ie non-adjacent vertices such that $N(u) = N(v)$. Since they share the same neighborhood, there is no maximal solution that only includes one of these vertices. Therefore, we can fold the two vertices into a new vertex $v'$ with weight $\w(v')=\w(v)+\w(u)$ as is done in the third case. 
     For the first case, if $\w(u)+\w(v) \geq \w(\I_{N(v)})$ and there is a solution $\I$ not including $v'$, then we can always construct an equal or better solution $\I'=\I\setminus N(v) \cup \{u,v\}$ and therefore we can include the vertices $u$ and $v$.
     
     In the second case, we first assume that $v' \in \I$ and show that $\w(\I_{N(v)}) + \aw(G-N[N[v]]) \geq \w(v)+\w(u) + \aw(G- N[v])$. This implies that the set $\I_{N(v)}$ is contained in some MWIS of $G$. Since $v'\in\I$ we know that
     \begin{align*}
     \w(u) +\w(v) + \aw(G') &= \w(u)+\w(v) +\w(v')+ \aw(G'- N_{G'}[v'])\\
     &=\w(u)+\w(v) + \w(\I_{N(v)}) - \w(u)-\w(v)+ \aw(G'-N_{G'}[v']) \\
     &=  \w(\I_{N(v)}) + \aw(G'- N_{G'}[v']).
     \end{align*}
     Since $\I'$ is an MWIS of $G'$, we have $\w(u) + \w(v) +\aw(G') \geq \w(u) + \w(v) + \aw(G'-v') = \w(u) + \w(v) + \aw(G- N[v]).$ Now suppose that $v'\notin \I$. For this case we show $\w(u) +\w(v) + \aw(G- N[v]) \geq \w(\I_{N(v)})+\aw(G- N[N[v]])$, which implies that $u$ and $v$ are in some MWIS of $G$. Since now $v\notin \I'$, we get $\w(u) +\w(v) + \aw(G') = \w(u)+\w(v) + \aw(G'-v') = \w(u) +\w(v) + \aw(G-N[v])$. As $\I'$ is an MWIS of $G'$, we also get 
     \begin{align*}
         \w(u) + \w(v) + \aw(G') &\geq \w(u) + \w(v) + \w(v') + \aw(G'-N_{G'}[v'])\\
         &= \w(u) + \w(v) + \w(\I_{N(v)}) - \w(u)-\w(v) + \aw(G-N[N[v]])\\
         &= \w(\I_{N(v)}) + \aw(G-N[N[v]]).
     \end{align*}
     Since $\I_{N(v)}$ is the only independent set in $G(N[v])$ with higher weight than $\w(u)+\w(v)$ $v$ is in some MWIS of $G$. We get $\aw(G) = \w(u) + \w(v) + \aw(G- N[v]) = \w(u) + \w(v) +\aw(G')$.
\end{proof}

Note that the third case of Reduction~\ref{red:e_twin} is already used in~\cite{kamis} but not introduced in this way by Lamm~\etal~\cite{lamm2019exactly}. Now we extend the idea of Reduction~\ref{red:e_twin}. For the following reduction rule we no longer require the neighborhoods of the two vertices $u$ and $v$ to be equal, but assume that $N(u) \subseteq N(v)$. If $\w(u) + \w(v) \geq \w(N(v))$ the vertex $u$ is always in an MWIS. The idea is that we can always replace vertices $x \in N(v)$ in a solution with $u$ and $v$ and thereby get an equal or better solution.
\pagebreak
\begin{figreduction}[Almost Twin]{a_twin}
    Let $u, v \in V$ be non-adjacent vertices such that $N(u) \subseteq N(v)$. If $\w(u) + \w(v) \geq \w(N(v))$, then include $u$.\\
\reductiondetails{$G' = G - N[u]$}{$\aw(G) = \aw(G') + \w(u)$}{$\I = \I' \cup \{u\}$}
\end{figreduction}

\begin{proof}
     Let $u, v \in V$ be non-adjacent vertices such that $N(u)\subseteq N(v)$ and $\w(u) + \w(v) \geq \w(N(v))$. Assume there is an MWIS $\I$ of $G$ not containing $u$. Then, there is a vertex $x\in N(u)$ such that $x\in \I$ which again results in $v \notin \I$. We can therefore construct a new solution $\I' = \I\setminus N(v)\cup \{u,v\}$ with $\w(\I') \geq \w(\I)$. It follows that there always exists an MWIS \hbox{that includes $u$.}
\end{proof}

\begin{remark}
    The Reduction~\ref{red:a_twin} can be extended as well. Instead of requiring that $\w(u) + \w(v) \geq \w(N(v))$ the same reduction can also be applied if $\w(u) + \w(v) \geq \aw(G[N(v)])$, since an MWIS for $N(v)$ is a subset of $N(v)$.
 \end{remark}

\begin{figure}[t]
    \centering
    	\begin{tikzpicture}
		\node[rnode] (z) at (5,4) {};
		\node[rnode=$v$] (v) [right=0.75cm of z] {};
		\node[rnode] (x) [below=0.6cm of z] {};
		\node[rnode] (y) [above=0.6cm of z] {};
		\node[rnode=$w$] (w) [below=0.6cm of v] {};
		\node[lnode=$u$] (u) [left=.75cm of z] {};
		\remainingG{5,2}{$G-(N[v]\cup \{u\})$}
		\draw[edge] (w) -- (5.95,2);
		\draw[edge] (w) -- (5.5,2);
		\draw[edge] (x) -- (5.5,2);
		\draw[edge] (x) -- (5.175,2);
		\draw[edge] (x) -- (4.85,2);
		\draw[edge] (w) -- (v);
		\draw[edge] (v) -- (x);
		\draw[edge] (v) -- (y);
		\draw[edge] (v) -- (z);
		\draw[edge] (u) -- (x);
		\draw[edge] (u) -- (y);
		\draw[edge] (u) -- (z);
		\draw[edge] (y) -- (z);
		\draw[edge] (x) -- (4.4,2);
		\node[weight] (vweight) at ($(v) + (.2, -.1)$) {$4$};
		\node[weight] (uweight) at ($(u) + (-.2, -.1)$) {$8$};
  		\node[weight] (xweight) at ($(x) + (0, -.1)$) {$2$};
  		\node[weight] (yweight) at ($(y) + (0, -.1)$) {$4$};
  		\node[weight] (zweight) at ($(z) + (0.25, -.1)$) {$3$};
  		\node[weight] (wweight) at ($(w) + (.25, -.1)$) {$2$};
		\draw[black,-{Stealth[scale=1.5]}] ($(v) + (1.25,0)$) -- ($(v) + (2.75,0)$) ;	
	    \begin{scope}[xshift=6cm]
     
		\node[nodeE] (z) at (5,4) {};
		\node[rnode=$v$] (v) [right=0.75cm of z] {};
		\node[nodeE] (x) [below=0.55cm of z] {};
		\node[nodeE] (y) [above=0.6cm of z] {};
		\node[rnode=$w$] (w) [below=0.6cm of v] {};
		\node[nodeI] (u) [left=.75cm of z] {};
		\remainingG{5,2}{$G-(N[v]\cup \{u\})$}
		\draw[edge] (w) -- (5.95,2);
		\draw[edge] (w) -- (5.5,2);
		\draw[edgeR] (x) -- (5.5,2);
		\draw[edgeR] (x) -- (5.175,2);
		\draw[edgeR] (x) -- (4.85,2);
		\draw[edge] (w) -- (v);
		\draw[edgeR] (v) -- (x);
		\draw[edgeR] (v) -- (y);
		\draw[edgeR] (v) -- (z);
		\draw[edgeR] (u) -- (x);
		\draw[edgeR] (u) -- (y);
		\draw[edgeR] (u) -- (z);
		\draw[edgeR] (y) -- (z);
		\draw[edgeR] (x) -- (4.4,2);
		\node[weight] (vweight) at ($(v) + (.2, -.1)$) {$4$};
  		\node[weight] (wweight) at ($(w) + (.25, -.1)$) {$2$};
	    \end{scope}
	\end{tikzpicture}
	\caption{Illustration of the Almost Twin (Reduction~\ref{red:a_twin}). In the original graph on the left, $N(u) \subseteq N(v)$ with $\w(u) + \w(v) \geq \w(N(v))$. By applying Reduction~\ref{red:a_twin}, $u$ is included and its neighbors excluded from $\I$, resulting in the reduced graph on the right.}\label{fig:a_twin}
 \end{figure}
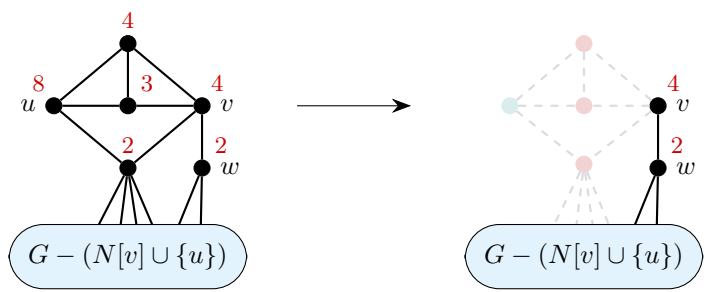
\begin{remark}
Note that applying Reduction~\ref{red:e_dom}, creates exactly the pattern needed for Reduction~\ref{red:a_twin}. If the weight constraint is satisfied, $u$ and $v$ are almost twins and Reduction~\ref{red:a_twin} \hbox{can be applied.}
\end{remark}

\paragraph*{Simplicial Vertex Based Reduction}
Recall that a simplicial vertex $v$ is a vertex where the neighborhood $N(v)$ forms a clique. The idea of Reduction~\ref{red:clique} is that a simplicial vertex $v$ can be included in the solution if it is of maximum weight in its neighborhood. This rule was first introduced by Lamm~\etal~\cite{lamm2019exactly}. In the following reduction we extend this idea further. 
Here, we consider a vertex $v$ which is almost simplicial, meaning there is one vertex $u\in N(v)$ such that if it is removed $N(v)\setminus \{u\}$ forms a clique. We call the pattern of the two vertices $u,v\in V$ such that $u\in N(v)$ and $N(v)\setminus \{u\}$ forms a clique for a $u$-$v$-\textit{funnel}.

The main idea of this reduction is that, under certain weight constraints, if $u$ and $v$ form a $u$-$v$-\textit{funnel} either $u$ or $v$ is in an MWIS. In this situation, only the three following solution patterns in $N[v]$ can occur. First, the vertex $v$ is in an MWIS. Second, the vertex $u$ and one other neighbor $x\in N(v)\setminus N[u]$ are in the solution. For this case we need to add additional edges between the vertex $u$ and the remaining vertices $x \in N(v) \setminus N[u]$. Third, only the vertex $u$ is part of the solution. Note that the third case can only occur when $\w(u)>\w(v)$. These three cases lead to the weighted version of the funnel reduction, see Reduction~\ref{red:wFunnel}. The unweighted version of this reduction was presented by Xiao~\etal~\cite{xiao2013confining}.

\begin{figreduction}[Weighted Funnel]{wFunnel}
Assume $u,v \in V$ forms a $u$-$v$-funnel and that $\w(v) \geq \max\{\w(x) \mid x \in N(v)\setminus \{u\}\}$. Furthermore, let $N'(v) = \{x \in N(v) \setminus N[u] \mid \w(x) + \w(u) > \w(v)\}$.
\begin{itemize}
    \item If $\w(v) \geq \w(u)$, fold $v$ and $u$ into its neighborhood.\\
    \reductiondetails{$G' = G - (N[v] \setminus N'(v))$ and for all $x \in N'(v)$, let $N(x) = N(x) \cup N(u)$ and \hbox{$\w(x) = \w(x) + \w(u) - \w(v)$}}{$\aw(G) = \aw(G') + \w(v)$}{If $\I' \cap N(v) = \emptyset$, then $\I = \I' \cup \{v\}$, else $\I = \I' \cup \{u\}$\\[-.5em]}
    \item If $\w(v) < \w(u)$, fold $v$ into its neighborhood.\\
    \reductiondetails{$G' = G - (N[v] \setminus (N'(v)\cup \{u\}))$ with $\w(u) = \w(u) - \w(v)$ and for all $x \in N'(v)$, let $N(x) = N(x) \cup N(u)$}{$\aw(G) = \aw(G') + \w(v)$}{If $\I' \cap N(v) = \emptyset$, then $\I = \I' \cup \{v\}$}
\end{itemize}
\end{figreduction}

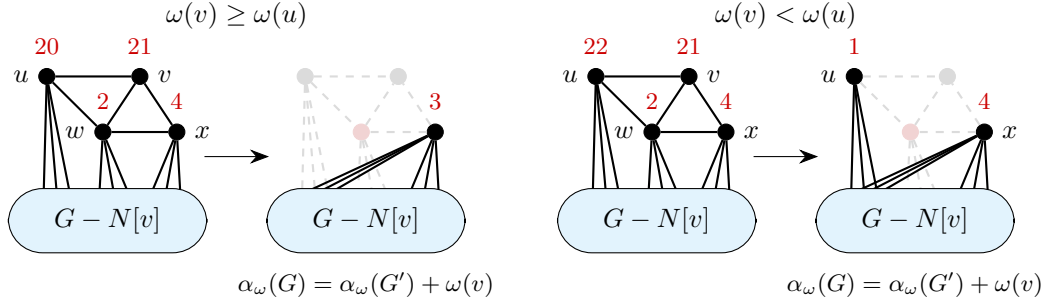
\begin{figure}[t]
    \centering
    	\begin{tikzpicture}[scale=0.85]

		\node[rnode=$v$] (v) at (5.5,4.25) {};
		\node[below=0.5cm of v] (vtmp) {};
		\node[rnode=$x$] (x) [right=0.25cm of vtmp] {};
		\node[lnode=$w$] (w) [left=0.25cm of vtmp] {};
		\node[lnode=$u$] (u) [left=1cm of v] {};
		\remainingG{5,2}{$\quad G-N[v]\quad$}
		\draw[edge] (x) -- (6.15,2);
		\draw[edge] (x) -- (5.825,2);
		\draw[edge] (x) -- (5.5,2);
		\draw[edge] (w) -- (5.5,2);
		\draw[edge] (w) -- (5.175,2);
		\draw[edge] (w) -- (4.85,2);
		\draw[edge] (x) -- (v);
		\draw[edge] (w) -- (v);
		\draw[edge] (x) -- (w);
		\draw[edge] (u) -- (v);
		\draw[edge] (u) -- (w);
		\draw[edge] (u) -- (4.5,2);
		\draw[edge] (u) -- (4.25,2);
		\draw[edge] (u) -- (4,2);
		\node[weight] (vweight) at (v) {$21$};
		\node[weight] (uweight) at (u) {$20$};
		\node[weight] (wweight) at (w) {$2$};
		\node[weight] (xweight) at (x) {$4$};
		\draw[black,-{Stealth[scale=1.5]}] (6.5,3) -- (7.5,3) ;	
        \node at (7,5.2) {$\w(v) \geq \w(u)$};
	    \begin{scope}[xshift=4cm]
			\node[nodeR] (v) at (5.5,4.25) {};
			\node[below=0.5cm of v] (vtmp) {};
			\node[node=$x$] (x) [right=0.25cm of vtmp] {};
			\node[nodeE] (w) [left=0.25cm of vtmp] {};
			\node[nodeR=$u$] (u) [left=1cm of v] {};
    		\remainingG{5,2}{$\quad G-N[v]\quad$}
			\node (offset) at (5,1) {\small$\aw(G) = \aw(G')+\w(v)$};
			\draw[edge] (x) -- (6.15,2);
			\draw[edge] (x) -- (5.825,2);
			\draw[edge] (x) -- (5.5,2);
			\draw[edgeR] (w) -- (5.5,2);
			\draw[edgeR] (w) -- (5.175,2);
			\draw[edgeR] (w) -- (4.85,2);
			\draw[edgeR] (x) -- (v);
			\draw[edgeR] (w) -- (v);
			\draw[edgeR] (x) -- (w);
			\draw[edgeR] (u) -- (v);
			\draw[edgeR] (u) -- (w);
			\draw[edgeR] (u) -- (4.5,2);
			\draw[edgeR] (u) -- (4.25,2);
			\draw[edgeR] (u) -- (4,2);
		      \draw[edge] (x) -- (3.65,2.25);
		      \draw[edge] (x) -- (3.6,2.1);
		      \draw[edge] (x) -- (3.7,2);
			\node[weight] (xweight) at (x) {$3$};
	    \end{scope}
     \begin{scope}[xshift=8.5cm]
         
		\node[rnode=$v$] (v) at (5.5,4.25) {};
		\node[below=0.5cm of v] (vtmp)  {};
		\node[rnode=$x$] (x) [right=0.25cm of vtmp] {};
		\node[lnode=$w$] (w) [left=0.25cm of vtmp] {};
		\node[lnode=$u$] (u) [left=1cm of v] {};
		\remainingG{5,2}{$\quad G-N[v]\quad$}
		\draw[edge] (x) -- (6.15,2);
		\draw[edge] (x) -- (5.825,2);
		\draw[edge] (x) -- (5.5,2);
		\draw[edge] (w) -- (5.5,2);
		\draw[edge] (w) -- (5.175,2);
		\draw[edge] (w) -- (4.85,2);
		\draw[edge] (x) -- (v);
		\draw[edge] (w) -- (v);
		\draw[edge] (x) -- (w);
		\draw[edge] (u) -- (v);
		\draw[edge] (u) -- (w);
		\draw[edge] (u) -- (4.5,2);
		\draw[edge] (u) -- (4.25,2);
		\draw[edge] (u) -- (4,2);
		\node[weight] (vweight) at (v) {$21$};
		\node[weight] (uweight) at (u) {$22$};
		\node[weight] (wweight) at (w) {$2$};
		\node[weight] (xweight) at (x) {$4$};
		\draw[black,-{Stealth[scale=1.5]}] (6.5,3) -- (7.5,3);
        \node at (7,5.2) {$\w(v) < \w(u)$};

	    \begin{scope}[xshift=4cm]
			\node[nodeR] (v) at (5.5,4.25) {};
			\node[below=0.5cm of v] (vtmp)  {};
			\node[rnode=$x$] (x) [right=0.25cm of vtmp] {};
			\node[nodeE] (w) [left=0.25cm of vtmp] {};
			\node[lnode=$u$] (u) [left=1cm of v] {};
    		\remainingG{5,2}{$\quad G-N[v]\quad$}
			\node (offset) at (5,1) {\small$\aw(G) = \aw(G')+\w(v)$};
			\draw[edge] (x) -- (6.15,2);
			\draw[edge] (x) -- (5.825,2);
			\draw[edge] (x) -- (5.5,2);
			\draw[edgeR] (w) -- (5.5,2);
			\draw[edgeR] (w) -- (5.175,2);
			\draw[edgeR] (w) -- (4.85,2);
			\draw[edgeR] (x) -- (v);
			\draw[edgeR] (w) -- (v);
			\draw[edgeR] (x) -- (w);
			\draw[edgeR] (u) -- (v);
			\draw[edgeR] (u) -- (w);
			\draw[edge] (u) -- (4.5,2);
			\draw[edge] (u) -- (4.25,2);
			\draw[edge] (u) -- (4,2);
			\draw[edge] (x) -- (3.65,2.25);
			\draw[edge] (x) -- (3.6,2.1);
			\draw[edge] (x) -- (3.7,2);
			\node[weight] (xweight) at (x) {$4$};
			\node[weight] (uweight) at (u) {$1$};
	    \end{scope}
     \end{scope}
	\end{tikzpicture}    
	\caption{Illustration for the two cases of the Weighted Funnel, Reduction~\ref{red:wFunnel}.}\label{fig:wFunnel}
\end{figure}

Before giving the proof for the Reduction~\ref{red:wFunnel}, we introduce the following lemma using the transformation of the extended weighted struction reduction introduced and proven by Gellner~\etal~\cite{gellner2021boosting}. The intuition is that for a vertex $v$ in the graph, we can either choose $v$ to be in the solution or some independent set in $G[N(v)]$. These independent sets are encoded as vertices connected to form a clique such that only one independent set in $G[N[v]]$ can be chosen. If the weight of an independent set in the neighborhood of $v$ is less than the weight of $v$, it is not added in the transformation, since these vertices could always be swapped with $v$ to create an equal or higher weight independent set.
\begin{lemma}\label{lem:struction}
    Let $v \in V$ and $C$ be the set of independent sets in $G[N(v)]$ with higher weight than $v$. Then, we can construct an instance $G'$ by removing $N[v]$ and for each solution $c \in C$, add a vertex $v_c$ with weight $\w(c)-\w(v)$. All these solution vertices are connected to form a clique. Furthermore, all $v_c$ are connected to vertices in $N(c)\setminus N(v)$ that are adjacent to a vertex in $c$, but not adjacent to $v$. It then holds that $\aw(G) = \aw(G') + \w(v)$.
\end{lemma}
\begin{proof}
    Proven by Gellner~\etal~\cite{gellner2021boosting} in the extended weighted struction reduction.
\end{proof}

\begin{proof}[Proof for Reduction~\ref{red:wFunnel}]
    Assume $u,v \in V$ forms a $u$-$v$-\textit{funnel} and $\w(v) \geq \max\{\w(x) \mid x \in N(v) \setminus \{u\}\}$. Let $N'(v)$ be as it is defined in the reduction rule. As a first step, we prove that an MWIS contains either $u$ or $v$. Therefore, assuming $u$ is not in an MWIS, we can exclude $u$, resulting in $v$ being a simplicial vertex of maximum weight. Now, we can include $v$ (see Reduction~\ref{red:clique}). In the other case, vertex $u$ is in an MWIS. Using this, we can exclude their common neighborhood $N(v) \cap N(u)$.
    Since $N'(v)$ still forms a clique, only the following three cases for an MWIS $\I$ in $G$ have to be considered. First, $v$ is in $\I$, second only $u$ is part of $\I$ and third, $u$ and \emph{one} of the remaining neighbors $x\in N'(v)$ is in the solution.
    
    We now apply the graph transformation of Lemma~\ref{lem:struction} to the solutions described. Note that when the weight of the independent set is less than $\w(v)$, no vertex is added in the transformation. In our case, this means that we can remove the remaining neighbors $x\in N'(v)$ if $\w(x) + \w(u) < \w(v)$. 

    Furthermore, every solution including a vertex $x \in N'(v)$ has to also contain $u$ by the weight assumption $\w(v) \geq \w(x)$. The vertices $x \in N'(v)$ now represent the solutions containing $x$ and $u$. Therefore, all these remaining neighbors have to be connected to the neighborhood of $u$. By Lemma~\ref{lem:struction} it holds that $\aw(G) = \aw(G') + \w(v)$.

    Next, we consider the different weight relations between $u$ and $v$. If $\w(u) < \w(v)$, there is no solution of higher weight than $\w(v)$ that only contains $u$. Therefore, $u$ is not present in the transformed graph. 
    
    Otherwise, the vertex $u$ remains in the reduced instance.
    Note that in this case, there are edges connecting $u$ to all remaining neighbors of $v$. For each remaining neighbor $x \in N'(v)$, its weight is increased by $\w(u) - \w(v)$. Afterward, we apply Reduction~\ref{red:e_dom} to all pairs of $u$ and an $x\in N'(v)$, which removes all edges connecting $u$ to $x \in N'(v)$. Reduction~\ref{red:e_dom} also reduces the weights of the vertices in $N'(v)$ by the current weight of $u$, which was reduced by the weight of $v$. This results in the original weight of the neighbors since $\w(x) = \w(x) + \w(u) - \w(v) -(\w(u) - \w(v))$.
    These steps give the reduced graph $G'$ as described by Reduction~\ref{red:wFunnel}.  
\end{proof}
\paragraph*{Extended Heavy Set Reduction}
An independent set $\I$ is called a heavy set, if for any independent set $C$ in the induced subgraph $G[N(\I)]$ it holds $\w(N(C) \cap \I) \geq \w(C)$. This concept was introduced by~\cite{xiao2021efficient}. Using this definition, for every independent set $\tI$, an equivalent or higher weight independent set $\I^*$ can be constructed by $\I^* \coloneqq \tI \setminus N(\I) \cup \I$.
Therefore, these vertices can always be included.
We now extend the Reduction~\ref{red:heavy_set} presented by Xiao~\etal~\cite{xiao2021efficient} to the case of a heavy set of three vertices in Reduction~\ref{red:heavy_set3}.

\begin{reduction}[Heavy Set 3]\label{red:heavy_set3}
	Let $u$, $v$ and $w$ be vertices forming a heavy set, then include $u$, $v$ and $w$.\\
    \reductiondetails{$G' = G-N[\{u, v, w\}]$}{$\aw(G) = \aw(G') + \w(u) + \w(v)+ \w(w)$}{$\I = \I' \cup \{u,v,w\}$}	
\end{reduction}

\begin{proof}
    The proof follows by definition of a heavy set.
\end{proof}
\paragraph*{Extended Unconfined Reduction}
We extend the unconfined vertex reduction rule for the weighted version of MIS by Xiao et al.~\cite{xiao2024maximum}. The intuition behind this rule is that a vertex $v$ can be removed if a contradiction is obtained by assuming every MWIS of $G$ includes $v$. To define the rule, we first introduce the following definitions and lemma.

Let $S$ be an independent set of $G$. As in the original unconfined reduction, a vertex $x \in N(S)$ is called a \textit{child} of $S$ if $\w(x) \geq \w(S \cap N(x))$. For each vertex $y \in N(x) \setminus N[S]$, let $\tilde{\I}_y$ be the MWIS of $G[(N(x) \setminus \{y\}) \setminus N[S]]$. Unlike the original rule, here, a child is called an \textit{extending child} if for some $y \in N(x) \setminus N[S]$ it holds that $\w(x) \geq \w(S \cap N(x)) + \w(\tilde{\I}_y)$. Such a vertex $y$ is called a satellite of $S$. Intuitively, a satellite is a vertex that must be included in every MWIS under the assumption that $S$ is contained in every MWIS. If it was not, we would have the contradiction we are looking for.

\begin{lemma}
    \label{lem:unconfined}
    Let $S$ be an independent set that is contained in every MWIS of $G$. Then, every MWIS also contains the satellites from each extending child $x$ of $S$.
\end{lemma}

\begin{proof}
    Let $S$ be an independent set that is contained in every MWIS of $G$ and $x\in N(S)$ be an extending child. Assume, towards a contradiction, that there exists a satellite $y$ that is not part of every MWIS. By definition, $\w(x)$ is now greater or equal to $\w(S \cap N(x)) + \w(\tilde{\I}_y)$. Thus, for any MWIS that includes $S$ but not $y$, we can replace $S \cap N(x)$ and $\tilde{\I}_y$ with $x$ to obtain a greater or equally large independent set, which contradicts the assumption that $S$ was contained in every MWIS of $G$.
\end{proof}

\begin{reduction}[Extended Unconfined Vertices]\label{red:e_unconfined}
    A vertex $v \in V$ can be removed if it is unconfined proven by the following procedure.
    Start with a set $S=\{v\}$. We assume $S$ is contained in every MWIS in $G$. We can search for a contradiction by repeatedly extending $S$ with satellites from one extending child until one of the following conditions hold
    \begin{enumerate}
        \item \label{redcase:unconfined1} There exists a child $x$ such that $\w(x) \geq \w(S \cap N(x)) + \aw(G[N(x) \setminus N[S]])$.
        \item \label{redcase:unconfined2} There exist no further satellites to extend $S$.
    \end{enumerate}
    In the second case, the set $S$ confines $v$, meaning every maximum weight independent set that contains $v$, also contains $S$ and we can not remove $v$.
    In the first case, $v$ is called unconfined and can be excluded.\\
        \reductiondetails{$G'=G-v$}{$\aw(G) = \aw(G')$}{$\I = \I'$} 
\end{reduction}
\begin{proof}
    Assume that $S$, initially just $\{v\}$, is contained in every MWIS and that Condition~\ref{redcase:unconfined1} holds. By Lemma \ref{lem:unconfined}, after every extension of $S$ with satellites, $S$ is still contained in every MWIS in $G$. But when Condition~\ref{redcase:unconfined1} holds, any MWIS $\I$ of $G$ with $S \subseteq \I$ can be modified using the child $x$ to obtain a new independent set $\I' = \{x\} \cup (\I \setminus N(x))$. From Condition~\ref{redcase:unconfined1}, it follows that $\w(\I') \geq \w(\I)$, breaking the assumption that $S$ is contained in every MWIS of $G$. Therefore, $v$ is removable.
\end{proof}
\begin{remark}
    This extended version of unconfined was already suggested in a remark by Xiao \etal~\cite{xiao2021efficient}, but they did not introduce it formally as is done here. It is also important to note that this version of unconfined is not practical in its most general form. Several MWIS problems need to be solved for each extending child, making it too computationally intensive for practical implementations. We implement this rule by restricting it to only search for satellites in neighborhoods that form an independent set. This way, we can detect multiple satellites from an extending child without solving any additional MWIS problems.
\end{remark}

\subsection{GNN Models}
\label{sec:gnn-models}

In the following, we introduce a new method for screening when to apply data-reductions based on Graph Neural Networks (GNN). We also provide a new and openly available supervised-learning dataset with graphs and labeled vertices for each reduction rule. For each expensive reduction rule used in the preprocessing, a pre-trained GNN model will limit the set of vertices that should be checked for using the reduction rule. GNNs are recent additions to the field of artificial intelligence that bring successful ideas from conventional deep learning to the irregular structure of graphs \cite{wu2020comprehensive}. Where traditional deep learning focuses on structured input, such as the grid of pixels in an image, GNNs accept the unstructured data of a graph.

As a first stage of this task, we evaluate the most popular GNN architectures used in combinatorial optimization, namely Graph Convolutional Network (GCN)~\cite{kipf2016semi} and Graph Sample and Aggregate (GraphSAGE)~\cite{hamilton2017inductive}. In addition to GCN and GraphSAGE, we also introduce a slightly altered GNN architecture, which we name the \lnr{}\textsc{-model} (LR). At a high level, these GNN architectures combine conventional neural networks with message passing. In this \textit{message passing}, every vertex sends messages to its neighbors and aggregates the received messages. One iteration of message passing and a neural network transformation makes up a GNN \textit{layer}. Several of these layers are stacked on top of each other to make up a GNN model. The goal is that after these message-passing layers, the final vertex embedding can be used to estimate how likely it is that a reduction rule will succeed at this vertex.

An undirected graph $G=(V, E, H)$ is given as input to the model, where $H$ is the initial feature representation for the vertices in the graph. Any number of vertex features can be used. Using the vertex weight is an obvious choice, but additional features such as vertex degrees and local clustering coefficients are common~\cite{lauri2023learning}. Note that computing complicated features for the vertices can add significant computational overhead to the model, which is important for our application since running the model should not take longer than checking the reduction rules. This feature representation will change after each layer in the model. For a vertex $u \in V$ at layer $l$, the feature representation is denoted by $H^{(l)}_u$. The length $d$ of the feature vector at layer $l$ is denoted by $d^{(l)}$. Stacking all the feature vectors at the $l$'th layer gives the matrix $H^{(l)} \in \mathbb{R} ^{|V| \times d^{(l)}}$. Independent from any input graph, every layer in the GNN model has trainable parameters $W^{(l)}$, bias $b^{(l)}$, and a non-linear activation function $\sigma$. Every model uses the $\textrm{ReLU}(x) \coloneqq \max(0,x)$ activation function for internal layers and $\textrm{Sigmoid}(x) \coloneqq \frac{1}{1 + e^{-x}}$ for the output layer. Note that these activation functions are applied element-wise when the input is a vector. With this, we define each model used for testing.

\paragraph*{Graph Convolutional Network}

The GCN architecture was the first successful extension of convolutional neural networks to work directly on graphs. At each layer in a GCN model, every node aggregates information from its immediate neighbors and combines it with its own data. After this, the information stored in each node is passed through the layer-specific neural network to create the node information for the next layer. The layer-wise propagation rule can be seen in Algorithm~\ref{alg:gcn}. Not that the GCN model assumes that there are self-edges added to the input graph. If not, a vertex will not include its own feature representation in the neighborhood aggregation.

\RestyleAlgo{ruled}
\begin{algorithm}[H]
    \DontPrintSemicolon
    \caption{GCN propagation rule. Self edges are added to the input graph and $T_u$ is a temporary variable holding the aggregated feature vectors of the neighbors of $u$.}\label{alg:gcn}
    \For{$u \in V$}{
        $T_u \gets \sum_{v \in N(u)} H^{(l)}_v / \sqrt{|N(v)|}$\;
        $H^{(l+1)}_u \gets \sigma(W^{(l)}\cdot T_u + b^{(l)})$\;
    }
\end{algorithm}

\paragraph*{Graph Sample and Aggregate}

GraphSAGE expands on GCN in two notable ways. First, it can use any differentiable aggregation function, such as aggregating the neighborhood based on mean, max, sum, or even more complicated functions, such as a long short-term memory (LSTM) machine-learning model. Second, GraphSAGE concatenates the feature vector of a vertex with the aggregated information from its neighborhood. This is in contrast to the self-edges of the GCN and allows information to skip from one layer to the next without going through the neighborhood aggregation. The general layer-wise propagation rule can be seen in Algorithm~\ref{alg:graphsage}.

\begin{algorithm}[H]
    \DontPrintSemicolon
    \caption{GraphSAGE propagation rule. $T_u$ is a temporary variable holding the aggregated feature vectors of the neighbors of $u$.}
    \label{alg:graphsage}
    \For{$u \in V$}{
        $T_u \gets \text{AGGREGATE}(\{H^{(l)}_v \: | \: v \in N(u)\})$\;
        $H^{(l+1)}_u \gets \sigma(W^{(l)}\cdot \text{CONCAT}(H^{(l)}_u, T_u) + b^{(l)})$\;
    }
\end{algorithm}

\paragraph*{Learn and Reduce}

The proposed \lnr{} (LR) architecture differs slightly from the GCN and GraphSAGE architectures. Instead of applying the weighted transformation after the aggregation, it is applied during the neighborhood aggregation. Algorithm \ref{alg:lnr} gives the exact layer-wise propagation rule. To give an intuition for why this approach could learn reduction rules better than GCN and GraphSAGE, consider the case of the extended single-edge reduction. For this reduction rule, we are looking for two adjacent vertices $u, v\in V$ such that $\w(v) \geq \w(N(v))-\w(u)$. For more information on this reduction rule, see Table~\ref{tab:full_tab_reductions} or Gu \etal~\cite{gu2021towards}. Assuming the node weights and neighborhood weights are given as input features, the LR architecture could conceivably detect this pattern during the first layer of the model since the concatenated feature representation of $u$ and $v$ would contain all the necessary information to decide if the rule can be applied. In contrast, GCN and GraphSAGE would aggregate the entire neighborhood before the first weighted transformation, potentially obscuring the necessary information required to make the correct prediction. Note that this is not a novel GNN architecture and could arguably fit into the GraphSAGE framework.

\begin{algorithm}[H]
    \DontPrintSemicolon
    \caption{\lnr{} propagation rule. The main difference with this architecture compared to GCN and GraphSAGE is the use of weighted transformation during neighborhood aggregation.}
    \label{alg:lnr}
    \For{$u \in V$}{
        $H^{(l+1)}_u \gets \text{MEAN}(\{ \sigma(W^{(l)}\cdot \text{CONCAT}(H^{(l)}_u, H^{(l)}_v) + b^{(l)}) \: | \: v \in N(u)\})$ \;
    }
\end{algorithm}

\subsection{Training Data Generation}
\label{sec:data-generation}

The supervised-learning dataset contains each instance twice, first the original instances without any reductions applied, and second the same instance after applying a set of fast reductions. These fast reductions are not considered relevant for early screening because they are so computationally cheap that we would always check them exhaustively. A comprehensive list of these reductions is given in Table~\ref{tab:full_tab_reductions}. We refer to these two versions of the same instance as \textsc{original} and \textsc{reduced}. The truth labels for each reduction rule are generated by testing the rule on each vertex without performing the reduction in the successful case. This way, each label is created using a clearly defined procedure without ambiguity. For the most costly reduction rules that rely on solving MWIS instances in subgraphs, a third class is used in the case of a timeout event. This gives the following labeling for each combination of a vertex and reduction rule.

$$
    \text{LABEL}(v) = 
    \begin{cases*}
        0 : \text{Unsuccessful reduction}\\
        1 : \text{Successful reduction}\\
        2 : \text{Timeout}\\
    \end{cases*}
$$

It is not always clear what vertex should be labeled as successful, especially for reductions where the rule starts from one vertex but ends up reducing another. It would be easiest if the starting vertex always received the successful label, given that the purpose of the trained model is to perform early screening. However, this would make the task of learning the rules unnecessarily hard for some reduction rules. The complete list of reduction rules used and how the labels are generated can be seen in Table \ref{tab:full_tab_reductions}. Note that the labeled vertices are always chosen to make the training as simple as possible. For example, the heavy set reduction starts from a source vertex and looks for two vertices in the neighborhood that can be included. If the successful labels were placed on the source vertices, the model would need to detect vertices adjacent to two vertices likely to be part of a solution. In contrast, if the labels are placed on the heavy vertices themselves, the pattern is simply to detect vertices that are likely part of a solution. We can still use the trained model for screening, but instead of deciding what vertices to start the reduction from, we only use suggested vertices from the neighborhood of the source vertex.

\subsection{The Learn and Reduce Approach}\label{sec:learn_and_reduce_approach}

\renewcommand{\floatpagefraction}{1}

In the following, we describe our \lnr~approach in more detail. In addition to our new reduction rules, we also use several data reduction rules introduced by \cite{butenko-trukhanov,gu2021towards,gellner2021boosting,xiao2021efficient}. The survey by Gro{\ss}mann~\etal~\cite{reduction_survey} presents a comprehensive collection of known reduction rules, and we enumerate them in the same way here. The complete list and a short description of each can be seen in Table \ref{tab:full_tab_reductions}. The most costly reduction rules, the ones below the bold line in the center of Table \ref{tab:full_tab_reductions}, are the ones that undergo early GNN screening before applying them. Each reduction rule will have its own pre-trained GNN model for early screening. Since we use the Sigmoid activation function, the output from the GNN will be a real number in the range from zero to one. If the GNN prediction for a vertex is greater than 0.5, that vertex will be checked using the reduction rule for which the model was trained.

The \lnr~reduction procedure maintains a queue of vertices to check for each reduction rule. The reduction rules are ordered based on complexity, and each rule is only checked when all the queues of easier reductions are empty. The ordering is the same as shown in Table~\ref{tab:full_tab_reductions}, except for Reduction~\ref{red:e_dom_rev} which is applied before Reduction~\hyperref[red:generalized_fold]{Reduction 2.5}. Even if it is more expensive, it yielded better results in our tests to place Reduction~\hyperref[red:cwis]{Reduction 6.8} before applying the struction approach introduced by Gellner~\etal~\cite{gellner2021boosting}. An outline of how the \lnr~reduction procedure works is shown in Algorithm~\ref{alg:learnAndReduce}. Each queue is initialized with all the vertices from the graph. Whenever a successful reduction occurs, the adjacent vertices that saw a change in their neighborhood are queued up again for every reduction. This way, the reduction procedure will always return to the easier reductions as soon as possible. The vertices are not ordered within the queues beyond the first-in, first-out principle. Reduced vertices may still reside in other queues but are ignored when they are popped from the queue.

\SetKwComment{Comment}{/* }{ */}
\RestyleAlgo{ruled}
\begin{algorithm}[!ht]
	\caption{The \lnr~reduction procedure.}\label{alg:learnAndReduce}
    \DontPrintSemicolon
    \KwData{Graph $G=(V,E,\w)$, set of reductions $R$, and set of GNN models $M$}
    \KwResult{Reduced instance $G$}
    $i \gets 0$\;
    $Q \gets \{ \{0, 1, ..., |V| - 1\} \cdot |R| \}$ \hfill \Comment{Initialize FIFO queues for each rule}
    \While{$i < |R|$}{
        \eIf{$Q[i] = \emptyset $}{
            $i \gets i + 1$\;
            \Comment{With the Initial/Initial-Tight configurations, the next If statement is changed accordingly}
            \If{$i < |R|$ and $Q[i] \neq \emptyset$ and $M[i] \, \text{exists}$}{
                $Q[i] \gets \emptyset$\;
                $P \gets M[i](G)$\;
                \For{$v \in V \text{ such that } P[v] > 0.5$}{
                    $push(Q[i], \, v)$\;
                }
            }
        }{
            $v = pop(Q[i])$\;
            \If{$v$ \text{is not reduced}}{
                $G' \gets R[i](G, \, v)$ \hfill \Comment{Try reduction on $v$}
                \If{$G' \neq G$}{
                    $G \gets G'$ \hfill \Comment{Successful reduction}
                    $i \gets 0$\;
                    \For{\text{each changed vertex} $u \in V$ \text{and} $q     \in Q$}{
                        \If{$u \notin q$}{
                            $push(q, \, u)$\;
                        }
                    }
                }
            }
        }
    }
\end{algorithm}

\begin{table}
\caption{The full list of reductions used in \lnr{} along with a brief description of each one. The id's are taken from~\cite{reduction_survey}. The last two columns connect the reduction rules to the supervised-learning dataset. Note that we often combine similar reduction rules for the dataset. The bold vertical line in the center of the table indicates the distinction between cheap and expensive reduction rules. }\label{tab:full_tab_reductions}
\resizebox{\textwidth}{!}{
\begin{tabular}{m{0.6cm}>{\RaggedRight}m{2.5cm}>{\RaggedRight}m{7cm}>{\RaggedRight}m{2cm}>{\RaggedRight\arraybackslash}m{1.5cm}}

\textbf{Id} & \textbf{Name} & \textbf{Brief Description} & \textbf{Data Name} &\textbf{Label} \\ \Xhline{4\arrayrulewidth}
2.2 & Neighborhood Removal~\cite{lamm2019exactly} \label{red:neighborhood} & Includes a vertex $v$ if $\w(v) \geq \w(N(v))$ & Neighborhood Removal & Included vertex \\ \midrule

1.1 & Degree One~\cite{gu2021towards}\label{red:deg-one} & Reduces or folds all degree-one vertices & Degree-One & Degree-one vertex \\ \midrule

1.2 & Triangle~\cite{gu2021towards} \label{red:triangle} & Reduces degree-two vertices with \hbox{adjacent neighbors} & \multirow{3}{2cm}{Degree-Two} & \multirow{3}{1.5cm}{Degree-two vertex} \\ \cmidrule{1-3}

1.3 & V-Shape~\cite{lamm2019exactly,gu2021towards} \label{red:vShape} & Reduces degree-two vertices with \hbox{non-adjacent neighbors} & & \\ \midrule

3.1 & Simplicial Vertex~\cite{lamm2019exactly} \label{red:clique_list} & Includes a vertex $v$ if $N(v)$ forms a clique \hbox{and ${\w(v)\geq \max_{u\in N(v)}\w(u)}$} & \multirow{3}{*}{Clique} & \multirow{3}{2cm}{Simplicial vertex} \\ \cmidrule{1-3}

3.2 & Simplicial Weight Transfer~\cite{lamm2019exactly} \label{red:clique_weight_transfer} & Folds a vertex $v$ if $N(v)$ forms a clique \hbox{and ${\w(v) \le \max_{u\in N(v)}\w(u)}$} & & \\ \midrule

4.1 & Domination~\cite{lamm2019exactly} \label{red:dom_list} & Exclude a vertex that dominates a lighter vertex & Domination & Excluded vertex \\ \midrule

4.2 & Basic Single Edge~\cite{gu2021towards} \label{red:bse} & Exclude a vertex $v$ if for some vertex $u \in N(v)$, $\w(N(u) \setminus N(v)) \leq \w(u)$ & Single Edge & Excluded vertex \\ \midrule

4.3 & Extended Single Edge~\cite{gu2021towards} \label{red:ese} & For $\{u, v\} \in E$, exclude $N(u) \cap N(v)$ if ${\w(v) \geq \w(N(v)) - \w(u)}$ & Extended Single Edge & Excluded vertices \\ \midrule

7.1 & Twin~\cite{lamm2019exactly} \label{red:twin_list} & Limited to independent neighborhood & \multirow{2}{*}{Twin} & \multirow{2}{*}{Both twins} \\ \cmidrule{1-3}

New   & Extended Twin \label{red:e_twin_list} & Fold or reduce all twins & & \\ \midrule

New   & Almost Twin \label{red:almost_twin_list} & Includes a vertex $v$ if for some non-adjacent $u$, $N(v) \subseteq N(u)$ and $\w(v) + \w(u) \geq \w(N(u))$ & Almost Twin & Included vertex \\ \midrule

New   & Weighted Funnel\label{red:wFunnel_list} & For $\{u,v\} \in E$ with $\w(v) = \max\{\w(w) : w \in N(v)\setminus\{u\}\}$ and $N(v)\setminus\{u\}$ forms a clique, we can fold or reduce parts of $N(v)$ & Funnel & Almost simplicial vertex \\ \midrule

2.3 & Clique Neighborhood Removal~\cite{lamm2019exactly} \label{red:clique_neighborhood} & Extension of \hyperref[red:neighborhood]{Reduction 2.2} for cliques \hbox{in the neighborhood of $v$} & --- & --- \\ \midrule

New & Extended Domination (Reverse) & For $\{u, v\} \in E$ with $N[u] \subseteq N[v]$ and $\w(u) < \w(v)$, remove $\{u,v\}$ and set $\w'(v) = \w(v) - \w(u)$ & --- & --- \\ \Xhline{4\arrayrulewidth}

New   & Extended Unconfined \label{red:eUnconf_list} & Exclude a vertex $v$ if a contradiction is obtained by assuming $v$ is part of every MWIS & Unconfined & Excluded vertex \\ \midrule

6.8 & Critical Weight Independent Set~\cite{butenko-trukhanov}\label{red:cwis} & Included the vertices of a (possibly empty) independent set $\I_c$ such that $\w(\I_c) - \w(N(\I_c)) = \max$ $\{\w(\I)-\w(N(\I)) : \I \text{ is an IS of }G \}$ & Critical Set & Included vertices \\ \midrule

5.4 & Extended Reduced Struction~\cite{gellner2021boosting} \label{red:struction} & A vertex $v$ and its neighborhood can be transformed into a clique where each vertex represents an independent set in $N[v]$ & --- & --- \\ \midrule

2.5 & Generalized Neighborhood Folding~\cite{lamm2019exactly,kamis} \label{red:generalized_fold} & Extension of \hyperref[red:neighborhood]{Reduction 2.2} with solving the MWIS in the neighborhood of $v$ & Generalized Fold & Included vertex \\ \midrule

2.7 & Heavy Set~\cite{xiao2021efficient} \label{red:heavy_set_list} & Finds a pair of vertices that can be included based on their combined neighborhoods & Heavy Set & Included vertices \\ \midrule

New   & Heavy Set 3 \label{red:heavy_set3_list} & Extension of \hyperref[red:heavy_set_list]{Reduction 2.6} for three vertices instead of two & Heavy Set 3 & Included vertices \\ \hline

6.6 & Cut Vertex~\cite{xiao2024maximum} \label{red:cut_vertex}& Fold or reduce an articulation point $v$ along with the smaller component in $G - v$ & Cut Vertex & Articulation points \\
\end{tabular}
}
\vspace{10pt}
\end{table}

There are multiple ways to incorporate the GNN screening. As reduction rules are applied, the graph continuously changes. Predictions made at the start of the reduction procedure might not be relevant later. Therefore, we propose three GNN screening configurations.

\paragraph*{Always}
In this configuration, the queue for each expensive reduction is only checked to see if it is non-empty. If that is the case, the queue is always replaced by new suggestions from the GNN model. This means the GNN model could be evaluated multiple times during \hbox{the reduction process.}

\paragraph*{Initial}
Instead of always evaluating the GNN model, the Initial configuration only evaluates the GNN model the first time the rule is used. Recall that each queue is initialized with all the vertices in the graph, meaning that the first time the rule is used is also when the potential screening effect is the highest. After this initial screening, the reduction rule goes back to working as without GNN screening, meaning if further reductions occur and vertices are added back to the queue, then these will be checked exhaustively.

\paragraph*{Initial Tight}
Again, the GNN model is only evaluated the first time the reduction rule is used. However, instead of returning to the normal queue-based application afterward, the queue cannot accept any additional vertices. This means that the initial suggestions from the GNN model are the only vertices that will be checked using that reduction rule. Intuitively, this configuration is the most aggressive form of early screening and will take the least time.

\hfill

For two reductions, namely \hyperref[red:cwis]{Reduction 6.8 Critical Set} and \hyperref[red:cut_vertex]{Reduction 6.6 Cut Vertex}, the queue base approach does not fit so well. The Critical Set reduction needs to construct a new bipartite graph and compute the maximum flow, and the Cut Vertex reduction needs to find all articulation points. Therefore, it is natural to apply these reductions globally. Regarding the reduction queues, these two reductions still have the same queues as the other reductions. However, as long as the queue is non-empty, the reduction is applied globally, and the entire queue is cleared afterward. This also extends to the GNN screening, but we only apply the reduction globally if the GNN suggests applying the rule on more than 1\,\% of the graph.

\section{The \pils{} Algorithm}
\label{sec:pils}

The proposed heuristic consists of two parts: (1) a simple iterated local search procedure we refer to as \baseline{}, and (2) a new, concurrent heuristic called \pils{} (Concurrent Hybrid Iterated Local Search). In the following section, we first give a high-level overview of the proposed approach and then provide a detailed description of each component of our heuristic.

\subsection{High-Level Description}

An uncomplicated heuristic implemented very efficiently is often better in practice than a complicated one that runs slowly, especially when the heuristic makes heavy use of random decisions. This can be seen from the results of programming competitions such as PACE (Parameterized Algorithms and Computational Experiments), where fast randomized local search has become the default heuristic strategy among the winning solvers~\cite{kellerhals2021pace, grossmann2022pace, bannach2023pace}. Note that the PACE competition setting restricts the time to solve the problem. We also used this approach for our baseline local search. First, randomly alter a local area in the current solution. Then, apply greedy improvements in random order until a local optimum is found. Finally, if the new local optimum is worse than the previous, backtrack the changes and repeat. Compared to more complicated heuristics, the main benefit of our \baseline{} heuristic is that it is inherently local. Regardless of the graph size, one iteration of the search will typically only touch vertices a few edges away from where the random alteration started. Backtracking to the best solution is also local, as long as queues are used to track changes. Using queues also differentiates our \baseline{} from other heuristics, such as \hils{}~\cite{hybrid-ils-2018}, that make a copy of the entire solution instead.

The high-level idea of \pils{} is a new metaheuristic we call \textsc{Concurrent Difference-Core Heuristic}. Instead of only trying to improve one solution, we maintain several and update each one concurrently. At fixed intervals, the solutions are used to create a \textsc{Difference-Core} (\pilscore) instance based on where the solutions differ. In other words, if a vertex is part of all or none of the solutions, it will not be part of the \pilscore{}. The intuition is that the intersection of the solutions is likely to be part of an optimal solution, and where they differ indicates areas where further improvements could be made. Since there is no guarantee that the intersection of the solutions is part of an optimal solution, the \textsc{Concurrent Difference-Core Heuristic} alternates between looking for improvements on the original instance and the \pilscore{}, were the \pilscore{} is always constructed according to the current solutions. While the general \textsc{Concurrent Difference-Core Heuristic} can work using any heuristic method, our \pils{} heuristic uses the \baseline{} \hbox{iterated local search.}

\subsection{Baseline Local Search}

The outline of the baseline local search is shown in Algorithm \ref{alg:ls}. Each search iteration starts by picking a random vertex $u$ from the graph. If the vertex is already in the solution $u \in S$, or the tightness of $u$ is one, \ie~$|N(u) \cap S| = 1$, then an alternating augmenting path (AAP) is used to perturb the current solution. If $u$ is not currently in the solution, it is added along with a random number of additional changes in its close proximity. The vertices that see a change in their neighborhood are considered “close proximity” to $u$. These vertices are continuously queued up in $Q$ to find greedy improvements more efficiently later. We also use the size of $Q$ to control the amount of random changes. We stop perturbing the solution once $|Q|$ exceeds $m_q$, where $m_q$ is a hyperparameter for the \baseline{} heuristic. The benefit of using $Q$ and $m_q$ is that it stabilizes the computational cost for one iteration across different graph densities. For instance, the number of changes needed to fill $Q$ in a dense area will be higher than in a sparse area.

\SetKwComment{Comment}{/* }{ */}
\RestyleAlgo{ruled}
\begin{algorithm}[ht]
    \DontPrintSemicolon
    \caption{Baseline Local Search}
    \label{alg:ls}
    \KwData{Graph $G=(V,E,w)$, independent set $S$, max queue size $m_q$, and time limit $t$}
    \KwResult{Improved independent set $S$}
    $Q \gets \emptyset$ \hfill \Comment{Always filled with vertices that observe changes to $S$}
    \While{time spent $< t$}{
        $\text{cost} \gets \w(S)$\;
        $u \gets \text{uniform random from }[0,|V|-1]$\;
        \eIf{$u \in S \; or \; |N(u) \cap S| = 1$}{
            AAP-\textit{moves}$(G, S, u)$\;
        }{
            $S = \{u\} \cap S \setminus N(u)$\;
            \While{$|Q| < m_q$}{
                $v \gets \text{random element from } Q$\;
                \eIf{$v \in S$}{
                    $S = S \setminus \{v\}$\;
                }{
                    $S = \{v\} \cap S \setminus N(v)$\;
                }
            }
        }
        $\textsc{greedy}(G, S, Q)$\;
        \If{$w(S) < \text{cost}$}{
            undo changes to $S$\;
        }
    }
    \Return{S}
\end{algorithm}

After perturbing the current solution, greedy improvements are made to find a new local optimum. Having the queue $Q$ of candidate vertices helps speed up the search for this new local optima. We incorporate three greedy improvement operators for \textsc{greedy} in \baseline{} described \hbox{in the following.}

\paragraph*{Neighborhood Swap}
For a vertex $u \notin S$, if $\w(u) > \w(N(u) \cap S)$, then the independent set obtained by inserting the vertex $u$ and removing all neighbors of $u$ that are currently in the solution, i.e. $S=\{u\} \cup (S \setminus N(u))$, leads to an independent set of higher weight. The total sum of the weights of each neighbor currently in the solution $\Sigma_{v \in (N(u) \cap S)} \: w(v)$ is maintained at all times. This additional data structure allows the neighborhood swap to be checked in $\mathcal{O}(1)$ time.

\paragraph*{Two-One Swap}
For a vertex in the current solution $u \in S$, consider the set $T = \{v \in N(u) \mid N(v) \cap S = \{u\} \}$ with vertices in the neighborhood of $u$ whose only neighbor included in the solution is $u$. If there exists a pair of vertices $\{x,y\} \subseteq T$ such that $\{x,y\} \notin E$ and where $w(u) < w(x) + w(y)$, then swapping $u$ for $x$ and $y$ leads to a better solution. Examining all two-one swaps can be done in $\mathcal{O}(|E|)$~\cite{andrade-2012}.

\paragraph*{Alternating Augmenting Path}
An alternating augmenting path (AAP) is a sequence of vertices that alternate between being in and out of the solution $S$. AAP moves are used both for perturbation and finding greedy improvements. When used for perturbation, the AAP is always extended in a random direction. After no more vertices can be added to the path, the entire AAP is applied to the solution unless a strictly improving prefix of the path exists. When searching for greedy improvements, the path always starts from a one-tight vertex and extends in the \hbox{best direction possible.}

We use a slightly modified version of the AAP introduced by Dong \etal~\cite{metamis}. Let $U$ be the set of vertices on the AAP in $S$, and $\overline{U}$ be the vertices on the AAP not in $S$. A valid AAP has the property that the set $S'$ obtained by applying the AAP, i.e. $S'=(S \setminus U) \cup \overline{U}$ is also an independent set. This means $\overline{U}$ must also be an independent set where $N(\overline{U}) \cap S = U$. The main purpose of APPs is to find improving paths where $\w(\overline{U}) > \w(U)$, and swap the vertices in $U$ for those in $\overline{U}$ to obtain a heavier independent set. As a secondary use case, they can also perturb the current solution. The benefit of using AAPs instead of random perturbation is that the cardinality of the new solution can only be one worse than the original. Constructing an APP starts with either a single vertex $u \in S$ or a pair of vertices $v \notin S, u \in S$, such that $N(v) \cap S = \{u\}$. The AAP is then extended from the last vertex $u$ on the path two vertices $x \in N(u), y \in N(x)$ at a time, such that the following three \hbox{conditions are met:}
\begin{enumerate}
    \item $x$ is not adjacent to any vertices in $\overline{U}$
    \item $x$ is not currently in $\overline{U}$
    \item $x$ is adjacent to precisely two vertices in the solution $N(x) \cap S= \{u,\, y\}$
\end{enumerate} 
Notice that $y$ is allowed to already be on the path. In this scenario, after applying the AAP to $S$, $x$ would not have neighbors currently in the solution. This is the main difference to the definition of AAP given in~\cite{metamis}.

\subsection{CHILS}

\begin{figure}
    \centering
    	\begin{tikzpicture}
    \node[node, fill=red] (z) at (5,4) {};
    \node[below=1.2cm of z] (sol)  {solution 1};
    \node[node, fill=green] (v) [right=0.7cm of z] {};
    \node[node, fill=red] (x) [below=0.6cm of z] {};
    \node[node, fill=red] (y) [above=0.6cm of z] {};
    \node[node, fill=red] (w) [below=0.6cm of v] {};
    \node[node, fill=green] (u) [left=.7cm of z] {};
    \node[node, fill=red] (a) [left=0.3 of u] {};
    \node[node, fill=red] (b) [above=0.4 of v] {};
    \node[node, fill=green] (c) [above=0.3 of y] {};
    \node[node, fill=green] (d) [below=0.6 of u] {};
    \draw[edge] (w) -- (v);
    \draw[edge] (v) -- (x);
    \draw[edge] (v) -- (y);
    \draw[edge] (v) -- (z);
    \draw[edge] (u) -- (x);
    \draw[edge] (u) -- (y);
    \draw[edge] (u) -- (z);
    \draw[edge] (y) -- (z);
    \draw[edge] (w) -- (x);
    \draw[edge] (a) -- (u);
    \draw[edge] (a) -- (c);
    \draw[edge] (b) -- (z);
    \draw[edge] (b) -- (v);
    \draw[edge] (c) -- (y);
    \draw[edge] (d) -- (x); 
    
     \begin{scope}[xshift=3.5cm]
    		\node[node, fill=green] (z) at (5,4) {};
            \node[below=1.2cm of z] (sol) {solution 2};
    		\node[node, fill=red] (v) [right=0.7cm of z] {};
    		\node[node, fill=red] (x) [below=0.6cm of z] {};
    		\node[node, fill=red] (y) [above=0.6cm of z] {};
    		\node[node, fill=green] (w) [below=0.6cm of v] {};
    		\node[node, fill=red] (u) [left=.7cm of z] {};
            \node[node, fill=red] (a) [left=0.3 of u] {};
            \node[node, fill=red] (b) [above=0.4 of v] {};
            \node[node, fill=green] (c) [above=0.3 of y] {};
            \node[node, fill=green] (d) [below=0.6 of u] {};
    		\draw[edge] (w) -- (v);
    		\draw[edge] (v) -- (x);
    		\draw[edge] (v) -- (y);
    		\draw[edge] (v) -- (z);
    		\draw[edge] (u) -- (x);
    		\draw[edge] (u) -- (y);
    		\draw[edge] (u) -- (z);
    		\draw[edge] (y) -- (z);
            \draw[edge] (w) -- (x);
            \draw[edge] (a) -- (u);
            \draw[edge] (b) -- (v);
            \draw[edge] (c) -- (y);
            \draw[edge] (b) -- (z);
            \draw[edge] (d) -- (x);   
            \draw[edge] (a) -- (c);     
    \end{scope}

     \begin{scope}[xshift=7cm]
    		\node[node, fill=red] (z) at (5,4) {};
            \node[below=1.2cm of z]  (sol) {solution 3};
    		\node[node, fill=red] (v) [right=0.7cm of z] {};
    		\node[node, fill=red] (x) [below=0.6cm of z] {};
    		\node[node, fill=red] (y) [above=0.6cm of z] {};
    		\node[node, fill=green] (w) [below=0.6cm of v] {};
    		\node[node, fill=green] (u) [left=.7cm of z] {};
            \node[node, fill=red] (a) [left=0.3 of u] {};
            \node[node, fill=green] (b) [above=0.4 of v] {};
            \node[node, fill=green] (c) [above=0.3 of y] {};
            \node[node, fill=green] (d) [below=0.6 of u] {};
    		\draw[edge] (w) -- (v);
    		\draw[edge] (v) -- (x);
    		\draw[edge] (v) -- (y);
    		\draw[edge] (v) -- (z);
    		\draw[edge] (u) -- (x);
    		\draw[edge] (u) -- (y);
    		\draw[edge] (u) -- (z);
    		\draw[edge] (y) -- (z);
            \draw[edge] (w) -- (x);
            \draw[edge] (a) -- (u);
            \draw[edge] (b) -- (v);
            \draw[edge] (c) -- (y);
            \draw[edge] (d) -- (x);    
            \draw[edge] (b) -- (z);
            \draw[edge] (a) -- (c);    
    		\draw[black,-{Stealth[scale=1.5]}] ($(v) + (.5,0)$) -- ($(v) + (1.5,0)$);	
    		\draw[black,-{Stealth[scale=1.5]}] ($(v) + (.5,.5)$) -- ($(v) + (1.5,.5)$);	
    		\draw[black,-{Stealth[scale=1.5]}] ($(v) + (.5,-.5)$) -- ($(v) + (1.5,-.5)$);	   
    \end{scope}

     \begin{scope}[xshift=11.5cm]
    		\node[node] (z) at (5,4) {};
            \node[below=1.2cm of z] (sol)  {\pilscore};
    		\node[node] (v) [right=0.7cm of z] {};
    		\node[nodeE] (x) [below=0.6cm of z] {};
    		\node[nodeE] (y) [above=0.6cm of z] {};
    		\node[node] (w) [below=0.6cm of v] {};
    		\node[node] (u) [left=.7cm of z] {};
            \node[nodeE] (a) [left=0.3 of u] {};
            \node[node] (b) [above=0.4 of v] {};
            \node[nodeI] (c) [above=0.3 of y] {};
            \node[nodeI] (d) [below=0.6 of u] {};
    		\draw[edge] (w) -- (v);
    		\draw[edgeR] (v) -- (x);
    		\draw[edgeR] (v) -- (y);
    		\draw[edge] (v) -- (z);
    		\draw[edgeR] (u) -- (x);
    		\draw[edgeR] (u) -- (y);
    		\draw[edge] (u) -- (z);
    		\draw[edgeR] (y) -- (z);
            \draw[edgeR] (w) -- (x);
            \draw[edgeR] (a) -- (u);
            \draw[edge] (b) -- (v);
            \draw[edgeR] (c) -- (y);
            \draw[edgeR] (d) -- (x);    
            \draw[edge] (b) -- (z);
            \draw[edgeR] (a) -- (c); 
    \end{scope}

	\end{tikzpicture}
    \caption{Illustration for the \pilscore{}. Vertices that are part of all or non of the solutions are not part of the \pilscore{}.}
    \label{fig:pilscore}
\end{figure}
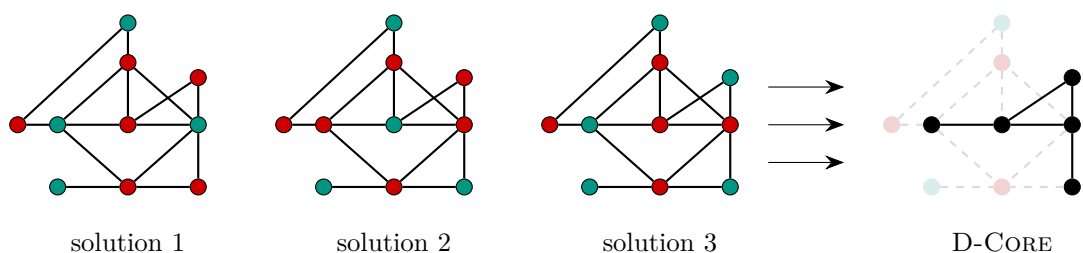

We give an overview of our concurrent heuristic \pils~in Algorithm~\ref{alg:pils}. First, we run \baseline{} $P$ times, with a time limit of $t_G$ seconds, on the full graph with different seeds and slightly modified $Q$ values to increase the randomness. We assign each of the computed solutions an \textit{id} and always keep track of the best solution found so far. This solution will only be modified if a better solution is found. Note that the id of the best solution can change during the algorithm's execution time. After the runs on the full graph are finished, \ie after $P\times t_G$ seconds, we construct the \pilscore{} instance. This is done by removing vertices that are part of all or none of the $P$ independent sets, see Figure~\ref{fig:pilscore} for an example.

\RestyleAlgo{ruled}
\begin{algorithm}[th]
    \DontPrintSemicolon
    \caption{The \pils~Algorithm}
    \label{alg:pils}
    \KwData{Graph $G=(V,E,\w)$, number of solutions $P$, max queue size $m_q$, time limit for the full graph $t_G$, time limit for the \pilscore{} $t_C$, and overall time limit $t$}
    \KwResult{Independent set $S$}
    $C = \{S_1, S_2, ..., S_P\}$ \hfill\Comment{using \textsc{greedy}($G,\,\emptyset,V)$ with different seeds}
    \While{time spent $< t$}{
        \textbf{parallel} \For{$S_i \in C$}{
            $S_i =$ \baseline{}($G, S_i, m_q + 4i, t_G$)\;
        }
        $G' \gets$ compute \pilscore{} using $C$\;
        \textbf{parallel} \For{$S_i \in C$}{
            $S' =$ \baseline{}($G', \emptyset, m_q + 4i, t_C$)\;
            \If{$\w(S') + \text{offset} \geq \w(S_i)$ OR ($i$ is odd AND $S_i$ is not best)}{
                apply $S'$ to $S_i$\;
            }
        }
        \If{$|V'| < 500$}{
            \textsc{parallel\_perturbe}($C$)\;
        }
    }
    \Return{$S \in C$ with largest weight}
\end{algorithm}

On the \pilscore{}, we again start our baseline local search $P$ times with a time limit of $t_C$ to generate new solutions. An independent set computed on the \pilscore{} can replace the previous solution with the same id. However, for solutions with an even id as well as for the best solution found so far, replacements are only made if the \pilscore{} solution is of higher weight. Letting half of the solutions always accept the \pilscore{} solution helps \hbox{diversify the search.}

Using the \pilscore{} helps concentrate the local search on the “more difficult” regions of the graph, where our $P$ solutions differ. However, since the areas where they agree are not necessarily part of an optimal solution, \pils{} alternates between using the \baseline{} on the original instance and the recomputed \pilscore{} based on the current $P$ solutions.

If the number of vertices in the \pilscore{} falls below some small constant value, it indicates that the $P$ solutions are all quite similar. Since this reduces the benefit of our approach, we perturb all solutions with an odd id, except for the best solution, and without backtracking in the case where the new local optima is worse. In Algorithm~\ref{alg:pils}, this is shown as \textsc{parallel\_perturbe} on line 14, where perturbing one solution is done as described in lines 9-16 of Algorithm~\ref{alg:ls}. As with accepting replacement solutions from the \pilscore{}, perturbing only half the solutions here helps to diversify the search and \hbox{escape local optima.}

\paragraph*{Parallel \pils{}}
Our \pils{} approach is easily parallelizable, with a natural choice for the number of solutions being exactly the number of cores available on the machine, allowing each solution to be improved simultaneously. In this configuration, the worst-case scenario will be similar to running the underlying \baseline{} sequentially, at least when technical details such as memory bandwidth and dynamic clock speeds are ignored. For larger numbers of solutions, $P$ should be divisible by the number of threads running to ensure that no threads are idle.

	\interfootnotelinepenalty=10000
\section{Experimental Evaluation}
\label{sec:experiments}

The following section introduces the experimental setup and establishes the datasets used for evaluating the proposed approaches. Then, we present GNN training results, preprocessing results, and finally state-of-the-art comparison for \baseline{} and \pils{}.

\subsection{Methodology} 
All the experiments were run on a machine with an Intel Xeon w5-3435X 16-core processor and 132\,GB of memory, running Ubuntu 22.04.4 with Linux kernel 5.15.0-113. Both \pils{} and \baseline{} were implemented in C and compiled with GCC version 11.4.0 using the -O3 flag. OpenMP is used for the parallel implementations. 
We evaluate eight instances in parallel for the sequential experiments. To ensure fairness between the algorithms, the instances start in the same order for each code, and only one program is evaluated at a time. We evaluate each program once for each instance.

We compare our algorithms \baseline{} and \pils{} to the state-of-the-art heuristic algorithms \htwis{} presented in \cite{gu2021towards} by Gu~\etal, the hybrid iterated local search \hils{} by Nogueira~\etal~\cite{hybrid-ils-2018}, the memetic algorithm \mmwiss{} by Gro{\ss}mann~\etal~\cite{mmwis}, the new metaheuristic \metamis{} by Dong~\etal~\cite{metamis}, and the Bregman-Sinkhorn algorithm \ilp{} by Haller and Savchynskyy~\cite{haller2024bregman}. The first three, \htwis{}, \hils{}, and \hils{}, were all used in their default configurations. The source code for \metamis{} is not publicly available, and therefore, we had to use the numbers reported in \cite{metamis}. They used the Amazon Web Service r3.4xlarge compute node running Intel Xeon Ivy Bridge Processors and wrote the implementation of their heuristic in Java. They also run their algorithm five times with different seeds and report the best solution found. For the Bregman-Sinkhorn algorithm \ilp{}, we use the variation that produces integer solutions only after reaching 0.1\,\% relative duality gap for the LP relaxation by recommendation from the authors~\cite{bsa_com}\footnote{The exact command we used was \lstinline|mwis_json -l temp_cont -B 50 --initial-temperature 0.01| \lstinline|-g 50 -b 100000000 -t [seconds] [instance]|}.

For the comparison of different algorithms we use performance profiles~\cite{dolan2002benchmarking}. At a high level, performance profiles show the relationship between the solution size or running time of each algorithm and the corresponding result produced by the competing algorithms. The performance profile gives each algorithm a non-decreasing, piecewise constant function. In general, the y-axis represents the fraction of instances where the objective function is better than or equal to $\tau$ times the best objective function value. For our solution quality comparison the y-axis shows the fraction of instances $\#\{\textnormal{weight} \geq \tau \ast \textnormal{best}\}/\#G$. Here, “weight” refers to the independent set weight computed by an algorithm on an instance, and “best” corresponds to the best result among all the algorithms shown in the plot. \#$G$ is the number of graphs in the dataset. The parameter $\tau$ is plotted on the x-axis. For maximization problems we have $0 < \tau \leq 1$.
When considering performance profiles comparing running time, the y-axis displays the fraction of instances where the time taken by an algorithm is less than or equal to $\tau$ times the time taken by the fastest algorithm on that instance, more formally $\#\{\textnormal{t}\leq \tau \ast \textnormal{fastest}\}/\#\textnormal{instances}$. Here, “t” represents the time taken by an algorithm on an instance, and “fastest” refers to the time taken by the fastest algorithm on that specific instance.
As we want to minimize the time, we have $\tau \geq 1$. In the plots, we refer to these y-axes by “Fraction of instances”. 
In general, algorithms are considered to perform well if a high fraction of instances are solved within a factor of $\tau$ as close to 1 as possible, indicating that many instances are solved close to or better/faster than the optimum/fastest solution found by all competing algorithms.

\subsection{Datasets}
We use two different sets of instances for our experiments. The first set of instances consists of graphs used in previous studies. We started with all the instances used by Gellner~\etal~\cite{gellner2021boosting} and Gu~\etal~\cite{gu2021towards}. Our set also consists of large social networks from the Stanford Large Network Dataset Repository (snap) \cite{snapnets}. Additionally, we added real-world graphs from OpenStreetMaps (osm) \cite{OSMWEB,barth-2016,cai-dynwvc}. Furthermore, as in Gu~\etal~\cite{gu2021towards} we took the same six graphs from the SuiteSparse Matrix Collection (ssmc) \cite{ssmcWEB,davis2011university} where weights correspond to population data. Each weight was increased by one, to avoid a large number of nodes assigned with zero weight. Additionally, we used instances from dual graphs of well-known triangle meshes (mesh) \cite{sander2008efficient}, as well as 3d meshes derived from simulations using the finite element method (fe) \cite{walshaw}. For unweighted graphs, we assigned each vertex a random weight that is uniformly distributed in the interval [1, 200]. Overall this results in 207 graphs from which we exclude instances which were already reducible by simple and fast reduction rules\footnote{We excluded instances which were fully reduced by running the following reductions, numbered as in \cite{reduction_survey}, in the given order: \hyperref[red:neighborhood]{Reduction 2.2}, \hyperref[red:deg-one]{Reduction 1.1}, \hyperref[red:triangle]{Reduction 1.2}, \hyperref[red:vShape]{Reduction 1.3}, \hyperref[red:clique]{Reduction 3.1}, \hyperref[red:dom]{Reduction 4.1}, \hyperref[red:bse]{Reduction 4.2}, \hyperref[red:ese]{Reduction 4.3}, Reduction~\ref{red:e_twin}, Reduction~\ref{red:a_twin}, Reduction~\ref{red:wFunnel}, \hyperref[red:clique_neighborhood]{Reduction 2.3}, Reduction~\ref{red:e_dom}, \hyperref[red:struction]{Reduction 5.4} (without blow up)}, resulting in a set of \numprint{83} graphs listed in Table~\ref{tab:graphs} \hbox{in the Appendix.}

Our second set of instances consists of \numprint{37} vehicle routing instances as used by Dong \etal~\cite{dong2021new}, see Table~\ref{tab:vr_graphs} in the Appendix. Initial warm-start solutions derived from the application and clique information for the graphs are also provided for these instances. The clique information is a clique cover of the graph, \ie a collection of potentially overlapping cliques that cover the entire graph. \metamis{}~\cite{metamis} and \ilp{}~\cite{haller2024bregman} use this clique information \hbox{in their algorithms.}

The early reduction screening dataset consists of the same instances as our first dataset without removing the reducible ones. There are 207 instances in the \textsc{original} instances and 71 \textsc{reduced} graphs. The difference between the number of \textsc{reduced} graphs and the 83 graphs listed in Table~\ref{tab:graphs} is due to the slight difference in reduction rules used to generate these two reduced sets. This dataset will be made publicly available.

\subsection{Training Results}
\label{sec:training-results}

We evaluate three GNN architectures on our new supervised learning task: GCN, GraphSAGE, and LR. Each architecture is configured to have a similar number of layers, parameters, and running time. We use two message-passing layers, followed by two layers of only weighted transformation. The input to the final two layers is a concatenation of all intermediate feature vectors, including the input features. The size of each intermediate feature vector is 16, and the size of the input vector is 8. The activation function ReLU is used for internal layers, and Sigmoid is used for the output layer. After each message passing layer, we perform random dropout with a probability of 0.2 during training. We use weighted binary cross entropy loss for loss function to adjust for the rarity of vertices labeled as successful. The weight is set exactly to the ratio between successful and unsuccessful labels in the training set. Finally, the Adam optimizer~\cite{kingma2014adam} is used with its default hyperparameters and 0.001 learning rate. The input features for one vertex are as follows.

\begin{enumerate}
    \item Node weight
    \item Neighborhood weight
    \item Minimum neighborhood weight
    \item Maximum neighborhood weight
    \item Node degree
    \item Average neighborhood degree
    \item Minimum neighborhood degree
    \item Maximum neighborhood degree
\end{enumerate}

It is important to note that these reduction rules can mostly be checked in polynomial time already. Therefore, it would defeat the purpose of early screening if the screening took more time than the ensuing reduction. For this reason, the models are kept as small as possible. The choice of 16 internal features is made based on the fact that BLAS kernels typically need at least 16 times 4 elements at single precision to utilize the CPU fully. In other words, going any lower than 16 would reduce the efficiency at which we can perform the necessary computations. For more information, we refer the reader to~\cite{goto2008anatomy}.

For the task of early reduction screening, we are only interested in the \textsc{reduced} graphs, i.e. the graph after applying the inexpensive reductions. There are 71 \textsc{reduced} graphs in total, with over 1 million labeled vertices combined. Each graph in the dataset is divided into three parts: (1) a training set containing 60\,\% of the vertices, (2) a validation set containing 20\,\%, and (3) a test set containing the last 20\,\%. Vertices labeled as 2 (timeout) are not used during training or validation. Each model is trained for 300 iterations, where one iteration uses the entire dataset. 

The training data can be considered a single large graph since we split the data in terms of vertices, not individual graphs. While this graph is derived from the same instances we use to evaluate \lnr{} later in the experiments, we argue that this is not an issue for the evaluation. First, the \textsc{reduced} graphs are snapshot images of how the instances looked after running a specific set of reductions. At inference, reductions are applied continuously as reducible vertices are found. This already means the inputs the GNN models will see at inference could be completely different from those seen during training. Furthermore, the reductions that lead to the \textsc{reduced} graphs are also not the same ones used in \lnr{}, and all the training results are given for the test set, which is a completely different set of vertices not seen during training.

In addition to the loss value, we also provide three additional metrics from the training procedure that indicate how well the trained model will perform at early reduction screening.
\begin{itemize}
    \item \textbf{Accuracy} The probability that a suggested vertex can be reduced. Higher is better.
    \item \textbf{Coverage} The fraction of reducible vertices suggested by the model. Higher is better.
    \item \textbf{Screening} The fraction of the total vertices in the graph suggested. Lower is better.
\end{itemize}

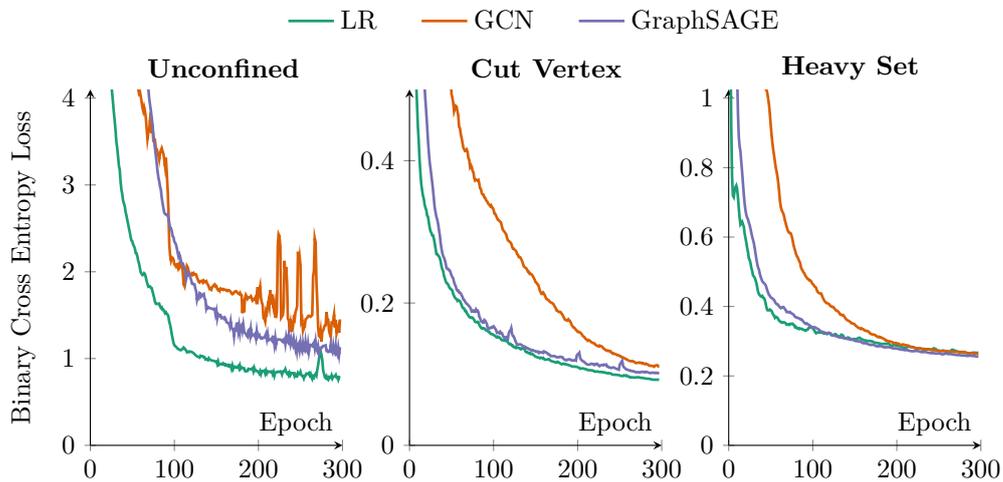
\begin{figure}[!ht]
    \centering
    \begin{tikzpicture}
        \node[align=center] at (1.75,5) {\textbf{Unconfined}};
        \begin{axis}[
            axis lines=middle,  
            axis x line*=bottom, 
            axis y line*=none, 
            cycle list/Dark2,
            ylabel={Binary Cross Entropy Loss},
            ylabel style={rotate=90, at={(axis description cs:-0.35,.5)}, anchor=north}, 
            xlabel={Epoch},
            width=0.35\linewidth,
            height=0.45\linewidth,
            legend columns=3,
            ymin=0.0, ymax=4.1,
            xmin=0, xmax=300]

            \addplot +[mark=none, line width=1pt] table [x=epoch, y=loss, col sep=semicolon] {data/gnn_training/unconfined/lr.csv};
            \addplot +[mark=none, line width=1pt] table [x=epoch, y=loss, col sep=semicolon] {data/gnn_training/unconfined/gcn.csv};
            \addplot +[mark=none, line width=1pt] table [x=epoch, y=loss, col sep=semicolon] {data/gnn_training/unconfined/sage.csv};
    
        \end{axis}
        
        \node[align=center] at (6,5) {\textbf{Cut Vertex}};
        \begin{axis}[
            axis lines=middle,  
            axis x line*=bottom, 
            axis y line*=none, 
            cycle list/Dark2,
            xlabel={Epoch},
            width=0.35\linewidth,
            height=0.45\linewidth,
            at={(0.3\linewidth,0)},
            legend style={draw=none,
                at={(0.5,1.13)}, 
                anchor=south,
                legend columns=3,
                /tikz/every even column/.append style={column sep=5mm}
            }, 
            ymin=0.0, ymax=0.5,
            xmin=0, xmax=300]

            \addplot +[mark=none, line width=1pt] table [x=epoch, y=loss, col sep=semicolon] {data/gnn_training/cut_vertex/lr.csv};
            \addplot +[mark=none,line width=1pt] table [x=epoch, y=loss, col sep=semicolon] {data/gnn_training/cut_vertex/gcn.csv};
            \addplot +[mark=none,line width=1pt] table [x=epoch, y=loss, col sep=semicolon] {data/gnn_training/cut_vertex/sage.csv};
    
            \legend{LR, GCN, GraphSAGE}
        \end{axis}
        
        \node[align=center] at (10,5) {\textbf{Heavy Set}};
        \begin{axis}[
            axis lines=middle,  
            axis x line*=bottom, 
            axis y line*=none, 
            cycle list/Dark2,
            xlabel={Epoch},
            width=0.35\linewidth,
            height=0.45\linewidth,
            at={(0.6\linewidth,0)},
            legend columns=3,
            ymin=0.0, ymax=1.025,
            xmin=0, xmax=300]

            \addplot +[mark=none, line width=1pt] table [x=epoch, y=loss, col sep=semicolon] {data/gnn_training/heavy_set/lr.csv};
            \addplot +[mark=none, line width=1pt] table [x=epoch, y=loss, col sep=semicolon] {data/gnn_training/heavy_set/gcn.csv};
            \addplot +[mark=none, line width=1pt] table [x=epoch, y=loss, col sep=semicolon] {data/gnn_training/heavy_set/sage.csv};
    
        \end{axis}
    \end{tikzpicture}
    \caption{Training loss for each architecture using the unconfined, critical set, and generalized fold reduction rules on the test set. Note that the scaling of the y-axes differ. }
    \label{fig:training-loss}
\end{figure}

\begin{figure}[!ht]
    \centering
    \begin{tikzpicture}
        
        \node[align=center] at (1.75,5) {\textbf{Accuracy}};
        \begin{axis}[
            axis lines=middle,  
            axis x line*=bottom, 
            axis y line*=none, 
            cycle list/Dark2,
            xlabel={Epoch},
            width=0.35\linewidth,
            height=0.45\linewidth,
            ymin=0.0, ymax=0.82,
            xmin=0, xmax=300]

            \addplot +[mark=none, line width=1pt] table [x=epoch, y=A, col sep=semicolon] {data/gnn_training/heavy_set/lr.csv};
            \addplot +[mark=none, line width=1pt] table [x=epoch, y=A, col sep=semicolon] {data/gnn_training/heavy_set/gcn.csv};
            \addplot +[mark=none, line width=1pt] table [x=epoch, y=A, col sep=semicolon] {data/gnn_training/heavy_set/sage.csv};
    
        \end{axis}
        
        \node[align=center] at (6,5) {\textbf{Coverage}};
        \begin{axis}[
            axis lines=middle,  
            axis x line*=bottom, 
            axis y line*=none, 
            cycle list/Dark2,
            legend style={draw=none,
                at={(0.5,1.13)}, 
                anchor=south,
                legend columns=3,
                /tikz/every even column/.append style={column sep=5mm}
            }, 
            xlabel={Epoch},
            width=0.35\linewidth,
            height=0.45\linewidth,
            at={(0.3\linewidth,0)},
            ymin=0.0, ymax=1.25,
            xmin=0, xmax=300]

            \addplot +[mark=none, line width=1pt] table [x=epoch, y=C, col sep=semicolon] {data/gnn_training/heavy_set/lr.csv};
            \addplot +[mark=none, line width=1pt] table [x=epoch, y=C, col sep=semicolon] {data/gnn_training/heavy_set/gcn.csv};
            \addplot +[mark=none, line width=1pt] table [x=epoch, y=C, col sep=semicolon] {data/gnn_training/heavy_set/sage.csv};
            \legend{LR, GCN, GraphSAGE}
        \end{axis}
        
        \node[align=center] at (10,5) {\textbf{Screening}};
        \begin{axis}[
            axis lines=middle,  
            axis x line*=bottom, 
            axis y line*=none, 
            cycle list/Dark2,
            xlabel={Epoch},
            width=0.35\linewidth,
            height=0.45\linewidth,
            at={(0.6\linewidth,0)},
            ymin=0.0, ymax=0.125,
            xmin=0, xmax=300,
            yticklabel style={/pgf/number format/fixed},]

            \addplot +[mark=none, line width=1pt] table [x=epoch, y=R, col sep=semicolon] {data/gnn_training/heavy_set/lr.csv};
            \addplot +[mark=none, line width=1pt] table [x=epoch, y=R, col sep=semicolon] {data/gnn_training/heavy_set/gcn.csv};
            \addplot +[mark=none, line width=1pt] table [x=epoch, y=R, col sep=semicolon] {data/gnn_training/heavy_set/sage.csv};
        \end{axis}
    \end{tikzpicture}
    \caption{Accuracy, coverage, and the fraction remaining after the screening for each architecture using the heavy set reduction rule. The values shown here are computed based on the test dataset. Note that the scaling of the y-axes differ.}
    \label{fig:critical_detailed}
\end{figure}
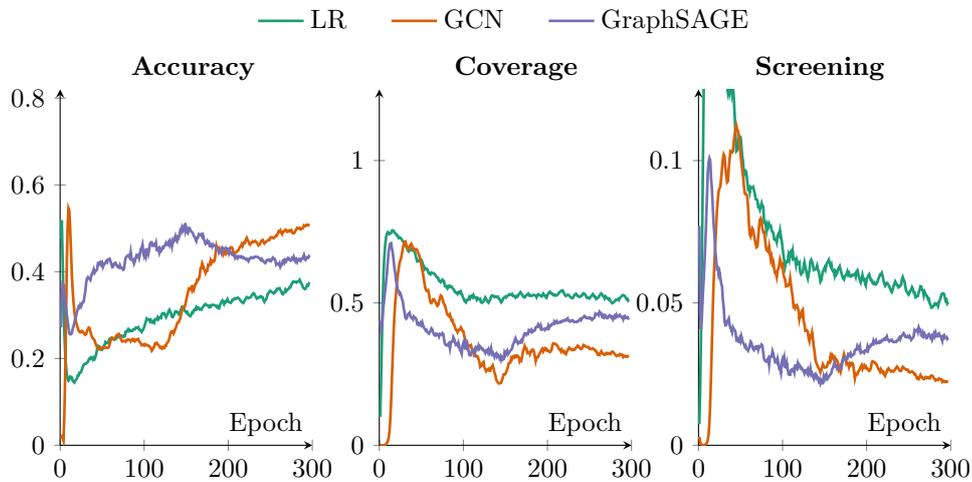

The detailed results for the reductions \hyperref[red:e_unconfined]{Unconfined}, \hyperref[red:cut_vertex]{Cut Vertex}, and \hyperref[red:heavy_set]{Heavy Set} can be seen in Figure \ref{fig:training-loss}. The proposed LR architecture clearly outperforms the GCN and GraphSAGE for the unconfined rule, but the difference is less noticeable for the critical and heavy sets. Figure \ref{fig:critical_detailed} provides accuracy, coverage, and the fraction of vertices removed by the screening for the heavy set reduction. Despite the similar-looking loss values for this reduction rule, it is clear from the other metrics that the models are discovering different things. GCN suggests the least amount of vertices. It also has the highest probability that a suggested vertex can actually be reduced. However, this is at the cost of catching less of the total reducible vertices. GraphSAGE and LR both suggest a larger fraction of the reducible vertices but at the cost of reduced suggestion accuracy and overall screening effect. Note that the models do not have any notion of accuracy or coverage during training, so reducing the loss of training data is the only thing being optimized for during training. To summarize the plots in Figure \ref{fig:critical_detailed} for LR, it suggests approximately 5\,\% of the graph for the heavy set rule, and among the suggested vertices are approximately 50\,\% of the vertices that actually can be reduced using the heavy set rule.

The final results for all the expensive reduction rules on the test set are given in Table~\ref{tab:detailed_gnn_training_all_reductions}. Training all the models took approximately 24 hours on our test machine using PyTorch Geometric~\cite{Fey/Lenssen/2019}. While PyTorch is a powerful framework for quickly developing and training models, several optimizations that could speed up the screening phase in our final program are not utilized. This includes combining the message passing and weighted transformation to reduce the number of times the data needs to pass through the cache hierarchy, loading feature vectors directly into AVX registers, and hard-coding the inner kernels for precisely the dimensions used by the models. To be clear, if the number of features had been higher, PyTorch would undoubtedly perform at a level close to the theoretical limit of the CPU. The optimizations listed above are useful because the number of features is low, and the main bottleneck is the memory bandwidth and latency. For this reason, we only continue with further experiments using one GNN architecture, and we do so using a manual implementation written in C. 
\begin{obs}{}{obs:trainings}
    Even though the differences are small, LR archives the best training loss and coverage on 5/6 instances compared to GCN and GraphSAGE. Therefore, we proceeded with the LR architecture.
\end{obs}
\begin{table}[t]
\caption{Detailed results in percentage on the test data for all reduction rules used in the \lnr~and GNN architectures. Red.~refers to the fraction of vertices in the data that can be reduced using the reduction rule. The \textbf{best} numbers are shown in bold, where Scr.,~Cov.,~Acc.,~and Loss are compared individually.}
\label{tab:detailed_gnn_training_all_reductions}
\resizebox{\textwidth}{!}{
\setlength\tabcolsep{4pt}
\begin{tabular}{lrrrrrrrrrrrrrrrrr}
   \multicolumn{2}{l}{} &\hspace{.2cm} & \multicolumn{4}{c}{GCN}   &\hspace{.2cm} & \multicolumn{4}{c}{GraphSAGE}   &\hspace{.2cm} & \multicolumn{4}{c}{LR}  \\
   \cmidrule{4-7}\cmidrule{9-12}\cmidrule{14-17}
Reduction   & \multicolumn{1}{l}{Red.} && \multicolumn{1}{l}{Scr.} & \multicolumn{1}{l}{Cov.} &\multicolumn{1}{l}{Acc.}&\multicolumn{1}{l}{Loss} && \multicolumn{1}{l}{Scr.} & \multicolumn{1}{l}{Cov.} & \multicolumn{1}{l}{Acc.} &\multicolumn{1}{l}{Loss} && \multicolumn{1}{l}{Scr.} & \multicolumn{1}{l}{Cov.} & \multicolumn{1}{l}{Acc.} &\multicolumn{1}{l}{Loss} \\

   \cmidrule{1-2}\cmidrule{4-7}\cmidrule{9-12}\cmidrule{14-17}
Critical    &6.19 &&7.03          &\textbf{44.06} & 38.81           & 0.45
                  &&\textbf{4.72} & 31.73         & 41.63           & 0.44
                  && 6.52         & 43.99         & \textbf{42.23}  & \textbf{0.43} \\
Cut         &0.53 &&0.03          & 2.89          &  55.91          & 0.11
                  &&\textbf{0.00} & 0.16          & 24.25           & 0.10
                  && 0.06         & \textbf{8.07} & \textbf{77.80}  & \textbf{0.09}\\
Gen. fold   &2.45 &&\textbf{0.14} & 2.34          & \textbf{42.17}  & 0.32
                  && 0.83         & 12.68         & 38.28           & 0.26
                  && 0.98         & \textbf{14.33}& 36.73           & \textbf{0.25} \\
Heavy set   &3.64 &&\textbf{2.25} & 31.39         & \textbf{50.65}  & \textbf{0.26}
                  && 3.77         & 44.84         & 43.39           & \textbf{0.26}
                  && 5.02         & \textbf{51.03}& 37.13           & \textbf{0.26} \\
Heavy set 3 &2.68 &&\textbf{0.47} & 7.88          & \textbf{45.45}  & \textbf{0.29}
                  && 0.51         & 5.74          & 30.60           & 0.34
                  && 0.90         & \textbf{9.35} & 28.00           & 0.33  \\
Unconfined  &7.94 && 10.15        & 25.95         & 25.37           & 1.45
                  && \textbf{5.65}& 17.20         & \textbf{27.07}  & 1.13
                  && 30.97        & \textbf{53.36} & 14.10          & \textbf{0.80}
\end{tabular}
}
\vspace{5pt}
\end{table}

\subsection{Learn and Reduce}
\label{sec:preprocessing-results}
This section presents experimental results for different configurations of our \lnr~approach. We exclude the vehicle routing instances here since most of them are very hard to reduce by any of our configurations. For the cyclic struction in \lnr~we have the two configurations \textit{Fast} and \textit{Strong} introduced by Gellner~\etal~\cite{gellner2021boosting}. We run the five different screening approaches introduced in Section~\ref{sec:learn_and_reduce_approach} for both of these configurations. In Table~\ref{tab:screening_configs}, we present the average results for all configurations running on the first set of 83 graphs. With the \textit{no gnn red} screening, we can see that these configurations are the fastest but are computing the largest reduced instances. When we use the expensive rule without screening in \textit{never}, we can see the potential gain in reduction size. In the \textit{Fast} configuration, we can further reduce the instances by \numprint{13,22}\,\% on average when using expensive reduction rules. However, this comes at the cost of more than doubling the preprocessing time on average. All of our \textit{Fast} screening methods can reduce the preprocessing time compared to \textit{never}. The fastest \textit{initial tight} is on average a factor of \numprint{1.7} times faster and manages to reduce the instances by an additional \numprint{12,59}\,\% compared to \textit{no gnn red}. For the \textit{Strong} configurations without screening, we can reduce the instances by \numprint{14,94}\,\% on average when using the expensive reduction rules. However, the configuration without any expensive rules already takes more time than all of the \textit{Fast} configurations using our screening approach. Therefore, we only focused on the different \textit{Fast} configurations in the following.

\begin{table}[t]
    \centering
    \caption{Arithmetic mean reduction results. Here, $K$ is the instance reduced by the different configurations, $\tilde{K}$ the instance reduced by \textit{Fast - no gnn red} and $G$ the original graph. We refer to the number of nodes $n$, and edges $m$ of the different graphs in the index. Furthermore, we give the offset and reduction time. In the column Time exp., we give the time which is used only for expensive reductions and the reduction screening. \vspace{-1em}}\label{tab:screening_configs}
    \begin{adjustbox}{width=\textwidth}
        \begin{tabular}{lcrrrrrrr}
            \textbf{Config} & \textbf{GNN}  & $\bf n_K/n_{\tilde{K}}$ & $\bf n_K/n_G$ & $\bf m_K/m_G$  & \textbf{Offset}       & \textbf{Time}  & \textbf{Time exp.} & $\bf \# n_K=0$   \\ \midrule
            \textit{Fast}   & no gnn red    & 100.00                  & 5.50          & 51.13          & 13\,179\,456          & \textbf{36.55} & 0.00               & 49 / 83          \\
            \textit{Fast}   & never         & 86.78                   & 4.77          & 50.85          & 13\,270\,743          & 74.74          & 38.19              & 52 / 83          \\
            \textit{Fast}   & always        & 87.11                   & 4.79          & 50.90          & 13\,269\,029          & 65.51          & 28.96              & 52 / 83          \\
            \textit{Fast}   & initial       & 87.13                   & 4.79          & 50.88          & 13\,268\,824          & 58.59          & 22.04              & 52 / 83          \\ \vspace{0.3em}
            \textit{Fast}   & initial tight & 87.41                   & 4.81          & 50.91          & 13\,267\,958          & 45.12          & 8.57               & 51 / 83          \\
            \textit{Strong} & no gnn red    & 98.74                   & 5.43          & 50.81          & 13\,188\,267          & 68.86          & 0.00               & 50 / 83          \\
            \textit{Strong} & never         & \textbf{85.06}          & \textbf{4.68} & \textbf{50.51} & \textbf{13\,280\,038} & 101.56         & 32.70              & 53 / 83          \\
            \textit{Strong} & always        & 85.56                   & 4.71          & 50.55          & 13\,278\,028          & 107.17         & 38.31              & 53 / 83          \\
            \textit{Strong} & initial       & 85.66                   & 4.71          & 50.56          & 13\,277\,789          & 88.24          & 19.38              & \textbf{54 / 83} \\
            \textit{Strong} & initial tight & 85.65                   & 4.71          & 50.58          & 13\,277\,974          & 85.20          & 16.34              & 53 / 83
        \end{tabular}
    \end{adjustbox}
\end{table}

\def\file{../data/mwis_results_no_vr/kernel_configs/chils_Fast_configurations}
\begin{figure*}
	\captionsetup[subfigure]{justification=centering}
	\centering
	\includegraphics{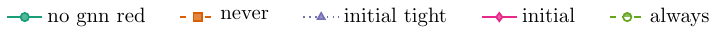}
	\begin{subfigure}[T]{0.5\textwidth}
		\centering
		\includegraphics[width=\textwidth]{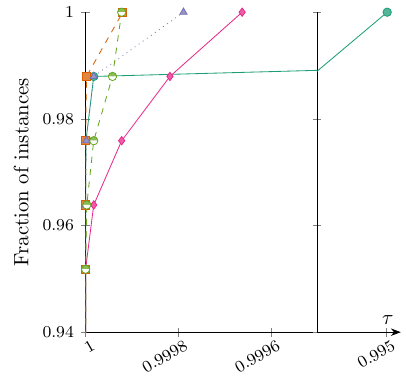}
	\end{subfigure}%
	\begin{subfigure}[T]{0.49\textwidth}
		\centering
		\includegraphics[width=\textwidth]{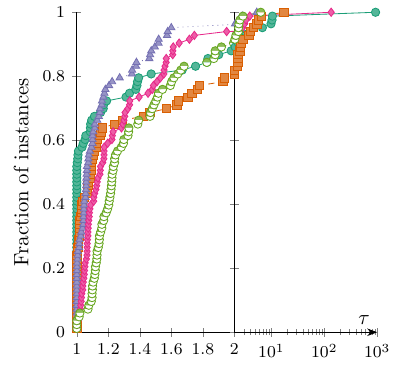}
	\end{subfigure}
	\caption{Performance profiles for different (\textit{Fast}) \lnr~configurations. After reducing the instance, we run the \pils~heuristic to evaluate the practical impact of the reduction configurations. We present solution quality, including reduction offset (left) and running time including reduction time (right). The vertical line in each plot marks the change from one \hbox{scale to another.} \vspace{-1.5em}}
	\label{fig:screening_chils}
\end{figure*}

In Figure~\ref{fig:screening_chils}, we present two performance profiles comparing the different screening approaches for the \textit{Fast} configuration with additionally running \pils{} on the reduced instances afterward. On the left, we have the solution quality achieved using one hour per instance for preprocessing and \pils{}. On the right we see the performance profile for the time needed to find the best solution within this time limit.
We can see that for all but five instances, all approaches achieve the same solution quality, and for the instances where the quality differs, the differences are minor. When considering the running times however, we can see definite differences. On more than \numprint{50}\,\% of the instances, not using the expensive reduction rules and having more time for the local search is the best strategy. However, this approach can be multiple orders of magnitude slower in other instances. Not using the screening for expensive rules, \ie \textit{never} is often more than twice as slow as other configurations.
Overall, \textit{initial tight} performed best on most instances regarding running time and is also very close to not using the screening regarding solution quality. There are only four instances where this approach took more than 1.6 times as long as the fastest solver on the respective instance. For all other configurations, this is true for more \hbox{than 12 instances.}

\begin{obs}{}{obs:configurations}
    The \textit{Strong} configuration without any expensive rules already takes more time than all of the \textit{Fast} configurations using our screening approach. Therefore, we only focus on the different \textit{Fast} configurations. The \textit{initial tight} performs best overall, especially for more difficult instances. Regarding solution quality, the \textit{initial tight} approach is very close to not using the screening. Therefore, we use the \textit{Fast - initial tight} configuration for our \lnr{} routine.
\end{obs}

\npthousandsep{\,}
\subsection{Parameter Tuning}\label{sec:parameter_tuning}
In this section, we perform experiments to find good choices for the number of concurrent solutions $P$ and the time $t$ spent improving them in the sequential version of \pils{}. We choose a subset of our two sets of instances for parameter tuning. We took six graphs from each set. First, two graphs from each graph class of the first set of instances where after \textit{fast early struction} the reduced graph has more than \numprint{500} vertices. For the vehicle routing instances, we took three graphs with more than \numprint{500000} vertices and three with less than \numprint{500000} vertices. Note that we did not include any \textit{mesh} or \textit{ssmc} graphs since these are all reduced to less than \numprint{500} vertices. The instances chosen for these experiments are marked with a $\star$ in tables~\ref{tab:graphs} and~\ref{tab:vr_graphs}.

We begin our experiments with the parameter defining the number of concurrent solutions, ${P\in\{8,16,32,48,64\}}$. The second parameter, highly correlating with $P$, is the time spent per local search run. For this experiment we set $t=t_G=t_C$ and test all ${t\in\{0.1,1,2.5,5,10,20,30\}}$ seconds. Recall that $t_G$ is the time spent on the original graph and $k_C$ is the time spent on the \pilscore{}. We present the geometric mean solution weight after one hour for the different configurations of these two parameters in Table~\ref{tab:parameter_tuning_ps}.

\begin{table}[b!]
    \centering
    \caption{Geometric mean weight computed by the sequential \pils~after one hour for different configurations for the number of solutions $P$ and time limits for the local search $t=t_G=t_C$. Note that the time $t$ is spent per solution in $P$, and not divided between them.}\label{tab:parameter_tuning_ps} 
    	\begin{tabular}{@{}lrrrr@{}}
    		& \multicolumn{1}{c}{$P=8$} & \multicolumn{1}{c}{$P=16$} & \multicolumn{1}{c}{$P=32$} & \multicolumn{1}{c}{$P=64$} \\ \cmidrule(l){2-5} 
    		$t=0.1$ & 14\,119\,825              & 14\,119\,386               & 14\,119\,516               & 14\,100\,937               \\
    		$t=1$   & 14\,124\,084              & 14\,126\,354               & 14\,127\,388               & 14\,117\,719               \\
    		$t=2.5$ & 14\,127\,282              & 14\,128\,074               & 14\,130\,823               & 14\,118\,203               \\
    		$t=5$   & 14\,127\,913              & 14\,129\,545               & 14\,129\,347               & 14\,117\,284               \\
    		$t=10$  & 14\,127\,684              & \textbf{14\,132\,120}              & 14\,124\,680               & 14\,116\,758               \\
    		$t=20$  & 14\,130\,214              & 14\,125\,181               & 14\,119\,467               & 14\,110\,491               \\
    		$t=30$  & 14\,126\,896              & 14\,118\,547               & 14\,115\,472               & 14\,112\,364              
    	\end{tabular}
\end{table}

In general, we can see that the best choice for $t$ gets lower as the number of solutions $P$ increase. Comparing the configurations with the most solutions, \ie $P=64$, the solution quality is worse than using fewer concurrent local search runs for all $t$ values tested. The best configuration we found is $P=16$ and $t=10$ seconds, resulting in a geometric mean weight of \numprint{14132120}. Note that these optimal choices are only for running \pils~sequentially.
\begin{obs}{}{obs:parameter_tuning}
    The best parameter configuration we found experimentally was ${P=16}$ and ${t=10}$ for running \pils~sequentially with a time limit of one hour per instance.
\end{obs}

\subsection{State-of-the-Art Comparison}
\label{sec:soa-results}

The state-of-the-art comparison in this section is divided into two parts. We start by discussing the results of the first set of instances. Here we compare our algorithms \baseline{} and \pils{} to state-of-the-art heuristic algorithms \htwis{} presented in \cite{gu2021towards} by Gu~\etal, the hybrid iterated local search \hils{} by Nogueira~\etal~\cite{hybrid-ils-2018} as well as the memetic algorithm \mmwiss{} by Gro{\ss}mann~\etal~\cite{mmwis}. The second part of this section is dedicated to the other algorithms and results for the vehicle routing instances. For these instances, we compare our heuristics with the results of \metamis{}~\cite{metamis}\footnote{The code is not publicly available, therefore, we could not rerun the experiments.}. Additionally, we compare against the new Bregman-Sinkhorn algorithm \ilp{} by Haller~\etal~\cite{haller2024bregman}. We do not evaluate \ilp{} on the first set of instances as \ilp{} requires clique information that is only available for the vehicle routing instances. Similarly, the vehicle routing instances differ significantly from previously established testing instances. The older state-of-the-art heuristics are not designed or implemented for these instances and perform significantly worse. In splitting the state-of-the-art comparison as described, we evaluate each heuristic in spirit with what they were designed for while demonstrating that \baseline{} and \pils{} are competitive in both categories.

Table~\ref{tab:soa_set1_sample} presents results for a subset of graphs from our first set. For this subset we choose the largest four graphs from each graph class. The results on the full set can be found in Table~\ref{tab:soa_set1} in the Appendix. We can see that both our algorithms \baseline{} and \pils{} outperform \htwis{} and \hils{} in terms of solution quality on \emph{every} instance. For large instances, such as the \textit{osm} instances, we find these best solutions even faster than any competitor find their best solution within the time limit. When comparing against \mmwiss{}, we observe that especially for the \textit{mesh} and \textit{snap} instances, our approaches find slightly worse solutions and need more time. However, \mmwiss{} already uses many of the reductions we incorporated in the \lnr{} preprocessing, which we did not use for \baseline{} or \pils{}. Comparing Table~\ref{tab:soa_set1_sample} to Table~\ref{tab:kernel_soa_set1_sample} clearly shows the importance of our preprocessing routine. In Table~\ref{tab:kernel_soa_set1_sample}, we present the results on the same instances, which have already been reduced with \lnr. From this it is clear that all evaluated algorithms generally benefit from the preprocessing. On these instances, \pils{} finds the best solutions on \emph{all} reduced instances.

\begin{figure}[!ht]
    \captionsetup[subfigure]{justification=centering}
    \centering
    \ref{named_non_vr_perf}
    \begin{subfigure}[t]{0.5\textwidth}
        \begin{tikzpicture}
            \begin{axis}[
                    perf_max,
                    perf_left, name=c2,
                    width=4.25cm,
                    xtick={1,0.9995,0.999},
                    xticklabels={$1$,0.9995,0.999},
                    xmin=0.999,
                    restrict x to domain=0.98:1]
                \addplot table[x=perf, y=fraction, col sep=comma] {data/datatool/chilsWeightPerf.tex};
                \addplot table[x=perf, y=fraction, col sep=comma] {data/datatool/baseWeightPerf.tex};
                \addplot table[x=perf, y=fraction, col sep=comma] {data/datatool/mmwisWeightPerf.tex};
                \pgfplotsset{cycle list shift=2}
                \addplot table[x=perf, y=fraction, col sep=comma] {data/datatool/htwisWeightPerf.tex};
                \addplot table[x=perf, y=fraction, col sep=comma] {data/datatool/hilsWeightPerf.tex};
            \end{axis}
            \begin{axis}[
                    perf_max, perf_right, at=(c2.south east),
                    width=4.25cm,
                    xmax=0.999,
                    xmin=0.89,
                    xtick={1, 0.96,0.93,0.9 },
                    xticklabels={1, 0.96,0.93,0.9 }
                ]
                \addplot table[x=perf, y=fraction, col sep=comma] {data/datatool/chilsWeightPerf.tex};
                \addplot table[x=perf, y=fraction, col sep=comma] {data/datatool/baseWeightPerf.tex};
                \addplot table[x=perf, y=fraction, col sep=comma] {data/datatool/mmwisWeightPerf.tex};
                \pgfplotsset{cycle list shift=2}
                \addplot table[x=perf, y=fraction, col sep=comma] {data/datatool/htwisWeightPerf.tex};
                \addplot table[x=perf, y=fraction, col sep=comma] {data/datatool/hilsWeightPerf.tex};
            \end{axis}
        \end{tikzpicture}
        \caption{Solution quality on original instances.}
        \label{fig:non-vr-performance-original-w}
    \end{subfigure}%
    \begin{subfigure}[t]{0.5\textwidth}
        \centering
        \begin{tikzpicture}
            \begin{axis}[
                    perf_min,
                    width=7.cm,
                    xmode=log,
                    log basis x=10,
                    xticklabel={$10^{\pgfmathprintnumber{\tick}}$}
                ]
                \addplot table[x=perf, y=fraction, col sep=comma] {data/datatool/chilsTimePerf.tex};
                \addplot table[x=perf, y=fraction, col sep=comma] {data/datatool/baseTimePerf.tex};
                \addplot table[x=perf, y=fraction, col sep=comma] {data/datatool/mmwisTimePerf.tex};
                \pgfplotsset{cycle list shift=2}
                \addplot table[x=perf, y=fraction, col sep=comma] {data/datatool/htwisTimePerf.tex};
                \addplot table[x=perf, y=fraction, col sep=comma] {data/datatool/hilsTimePerf.tex};
            \end{axis}
        \end{tikzpicture}
        \caption{Running time on original instances.}
        \label{fig:non-vr-performance-original-t}
    \end{subfigure}
    \begin{subfigure}[t]{0.5\textwidth}
        \centering
        \begin{tikzpicture}
            \begin{axis}[
                    perf_max,
                    perf_left,
                    width=4.25cm,
                    name=c1,
                    xtick={1, 0.995, 0.99,0.985},
                    xticklabels={$1$, $0.995$,$0.99$,$0.985$ },
                    xmin=0.99,
                    restrict x to domain=0.97:1,
                ]

                \addplot table[x=perf, y=fraction, col sep=comma] {data/datatool/kernel_data/kernel_chilsWeightPerf.tex};
                \addplot table[x=perf, y=fraction, col sep=comma] {data/datatool/kernel_data/kernel_baseWeightPerf.tex};
                \addplot table[x=perf, y=fraction, col sep=comma] {data/datatool/kernel_data/kernel_mmwisWeightPerf.tex};
                \pgfplotsset{cycle list shift=2}
                \addplot table[x=perf, y=fraction, col sep=comma] {data/datatool/kernel_data/kernel_htwisWeightPerf.tex};
                \addplot table[x=perf, y=fraction, col sep=comma] {data/datatool/kernel_data/kernel_hilsWeightPerf.tex};
            \end{axis}
            \begin{axis}[
                    perf_max, perf_right, at={(c1.south east)},
                    width=4.25cm,
                    xtick={1, 0.95,0.9,0.85 },
                    xticklabels={1, 0.95,0.9,0.85 },
                    xmax=0.99,
                    xmin=0.84,
                ]
                \addplot table[x=perf, y=fraction, col sep=comma] {data/datatool/kernel_data/kernel_chilsWeightPerf.tex};
                \addplot table[x=perf, y=fraction, col sep=comma] {data/datatool/kernel_data/kernel_baseWeightPerf.tex};
                \addplot table[x=perf, y=fraction, col sep=comma] {data/datatool/kernel_data/kernel_mmwisWeightPerf.tex};
                \pgfplotsset{cycle list shift=2}
                \addplot table[x=perf, y=fraction, col sep=comma] {data/datatool/kernel_data/kernel_htwisWeightPerf.tex};
                \addplot table[x=perf, y=fraction, col sep=comma] {data/datatool/kernel_data/kernel_hilsWeightPerf.tex};
            \end{axis}
        \end{tikzpicture}
        \caption{Solution quality on reduced instances.}
        \label{fig:non-vr-performance-reduced-w}
    \end{subfigure}%
    \begin{subfigure}[t]{0.49\textwidth}
        \centering
        \newcommand{\kernelSoaFile}{data/datatool/kernel_data/kernel_}
        \begin{tikzpicture}
            \begin{axis}[perf_min, 
                width=7.cm,
                xmode=log, log basis x=10, xticklabel={$10^{\pgfmathprintnumber{\tick}}$}, legend style = {legend columns=5,
                anchor=center,
                draw=none,
                /tikz/every even column/.append style={column sep=5mm}}, legend to name=named_non_vr_perf]

                \foreach \config in {chilsTimePerf,baseTimePerf,mmwisTimePerf}{
                        \addplot table[x=perf,y=fraction, col sep=comma] {\kernelSoaFile\config.tex};
                    }

                \pgfplotsset{cycle list shift=2}

                \foreach \config in {htwisTimePerf,hilsTimePerf}{
                        \addplot table[x=perf,y=fraction, col sep=comma] {\kernelSoaFile\config.tex};
                    }
                \legend{\chils, \baseline, \mmwiss, \htwis, \hils}

            \end{axis}
        \end{tikzpicture}
        \caption{Running time on reduced instances.}
        \label{fig:non-vr-performance-reduced-t}
    \end{subfigure}
    \caption{\label{fig:non-vr-performance-full} Performance profiles for solution quality and running time for original and reduced instances not including the vehicle routing instances. The reduction offset and reduction time are not added for the reduced instances, and fully reduced instances are not included. The vertical line in the plots on the left indicate a change from one scale to another. }
\end{figure}
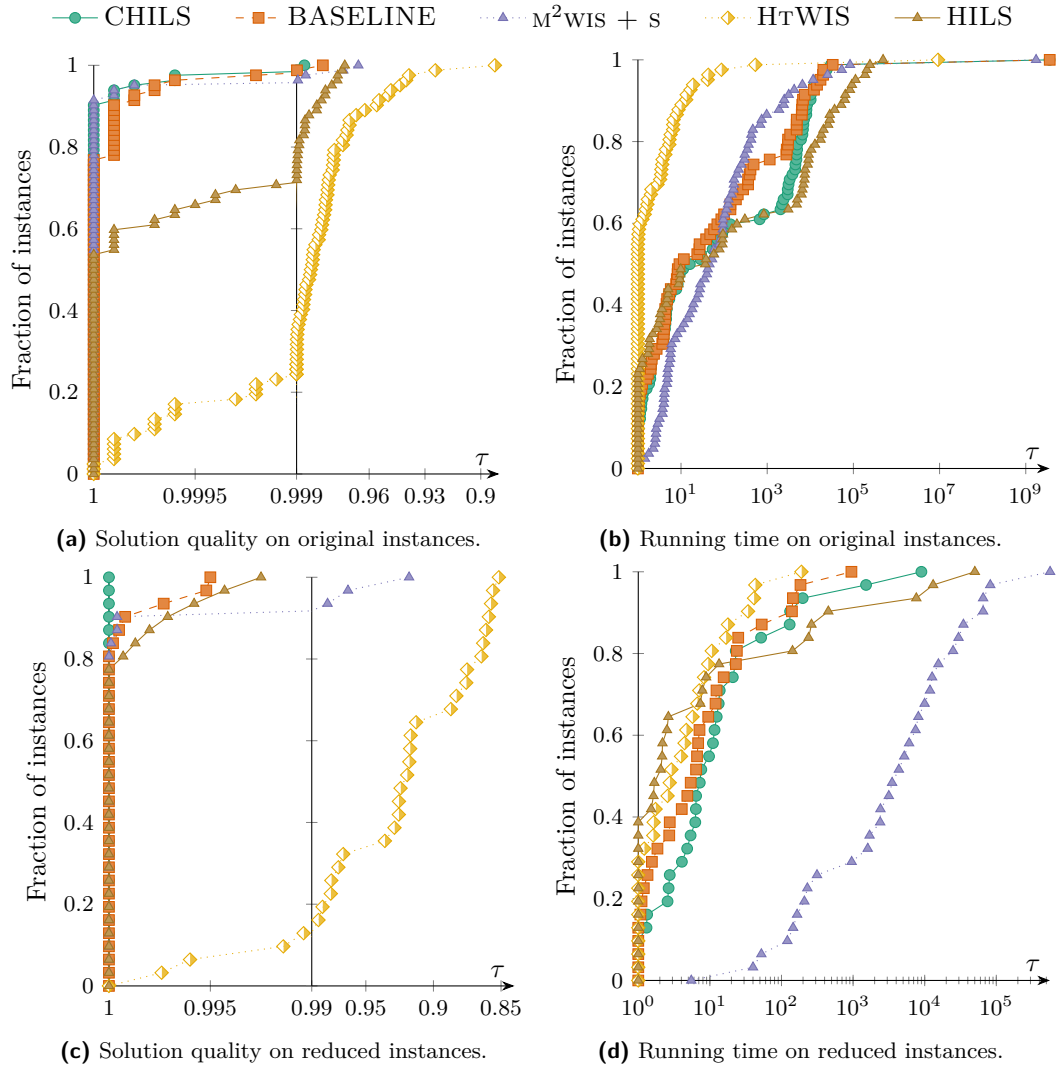

In Figure~\ref{fig:non-vr-performance-full}, we show performance profiles comparing the different algorithms' solution quality and running time on the full set of original and reduced instances. Note that most of the solutions on these instances found by \mmwiss{} are optimal~\cite{mmwis}. On approximately \numprint{60}\,\% of the original instances \htwis{} finds the solution fastest, see Figures~\ref{fig:non-vr-performance-original-t}. However, the fast running times come with a decrease in solution quality of up to 10\,\%, as can be seen in Figure~\ref{fig:non-vr-performance-original-w}. On more than \numprint{90}\,\% of the instances \pils{} finds the best solution while having comparable running times as the other schemes except for \htwis. The \mmwiss{} algorithm performs very similarly to \pils{} on the original instances. However, on the reduced instances, \mmwiss{} has the longest running times; see Figure~\ref{fig:non-vr-performance-reduced-t}. Note that the initial reductions used in \mmwiss{} are also included in the \lnr~framework. All other algorithms have similar running times on the reduced instances. From Figure~\ref{fig:non-vr-performance-reduced-w}, it follows that all algorithms except for \htwis{} find the same solution on almost 80\,\% of the reduced instances. On the other 20\,\% \pils{} finds the best solutions which are up to 0.7\,\% better than the \hils{} solutions and more than 7\,\% better than \mmwiss{}. On more than 85\,\% \htwis{} computes solutions more than 1\,\% worse than the best solutions found and for approximately 30\,\% of the instances these solutions are even more than 10\,\% worse.

\begin{obs}{}{obs:set_one_soa}
    All algorithms tested benefit from running the \lnr~preprocessing. \pils{} finds the best solution on \emph{all} reduced instances, while performing similar to \hils{} regarding running time. It is slightly slower than \htwis{}, which computes the worst solution quality on these instances. Without preprocessing, the improvement of \pils{} over \htwis{} and \hils{} regarding solution quality increases. Here, only \mmwiss{}, which already uses reductions, can find better results than \pils{} on some graphs.
\end{obs}

\begin{figure}[!ht]
    \captionsetup[subfigure]{justification=centering}
    \centering
    \ref{named}
    \begin{subfigure}[t]{0.5\textwidth}
        \begin{tikzpicture}
            \centering
            \begin{axis}[perf_left, xmin=0.975, name=p1]
                \foreach \column in {chils_cold,baseline,metamis_cold}{
                        \addplot+[] table[x={\column},y=fraction, col sep=comma] {data/t10_perf.csv};
                    }
                \pgfplotsset{cycle list shift=1}
                \foreach \column in {chils_warm,ilp,metamis_warm}{
                        \addplot+[] table[x={\column},y=fraction, col sep=comma] {data/t10_perf.csv};
                    }
            \end{axis}
            \begin{axis}[perf_right, xmax=0.975, xmin=0.4, at={(p1.south east)}]
                \foreach \column in {chils_cold,baseline,metamis_cold}{
                        \addplot+[] table[x={\column},y=fraction, col sep=comma] {data/t10_perf.csv};
                    }
                \pgfplotsset{cycle list shift=1}
                \foreach \column in {chils_warm,ilp,metamis_warm}{
                        \addplot+[] table[x={\column},y=fraction, col sep=comma] {data/t10_perf.csv};
                    }
            \end{axis}
        \end{tikzpicture}
        \caption{Time limit 6 minutes}\label{fig:vr-performance_less_time10}
    \end{subfigure}%
    \begin{subfigure}[t]{0.5\textwidth}
        \centering
        \begin{tikzpicture}
            \begin{axis}[perf_left, xmin=0.975, name=p1]
                \foreach \column in {chils_cold,baseline,metamis_cold}{
                        \addplot+[] table[x={\column},y=fraction, col sep=comma] {data/t50_perf.csv};
                    }
                \pgfplotsset{cycle list shift=1}
                \foreach \column in {chils_warm,ilp,metamis_warm}{
                        \addplot+[] table[x={\column},y=fraction, col sep=comma] {data/t50_perf.csv};
                    }
            \end{axis}
            \begin{axis}[perf_right, xmax=0.975, xmin=0.4, at={(p1.south east)}]
                \foreach \column in {chils_cold,baseline,metamis_cold}{
                        \addplot+[] table[x={\column},y=fraction, col sep=comma] {data/t50_perf.csv};
                    }
                \pgfplotsset{cycle list shift=1}
                \foreach \column in {chils_warm,ilp,metamis_warm}{
                        \addplot+[] table[x={\column},y=fraction, col sep=comma] {data/t50_perf.csv};
                    }
            \end{axis}
        \end{tikzpicture}
        \caption{Time limit 30 minutes}\label{fig:vr-performance_less_time50}
    \end{subfigure}
    \begin{subfigure}[t]{0.5\textwidth}
        \centering
        \begin{tikzpicture}
            \begin{axis}[perf_left, xmin=0.975, name=p1]
                \foreach \column in {chils_cold,baseline,metamis_cold}{
                        \addplot+[] table[x={\column},y=fraction, col sep=comma] {data/t100_perf.csv};
                    }
                \pgfplotsset{cycle list shift=1}
                \foreach \column in {chils_warm,ilp,metamis_warm}{
                        \addplot+[] table[x={\column},y=fraction, col sep=comma] {data/t100_perf.csv};
                    }
            \end{axis}
            \begin{axis}[perf_right, xmax=0.975, xmin=0.4, at={(p1.south east)}]
                \foreach \column in {chils_cold,baseline,metamis_cold}{
                        \addplot+[] table[x={\column},y=fraction, col sep=comma] {data/t100_perf.csv};
                    }
                \pgfplotsset{cycle list shift=1}
                \foreach \column in {chils_warm,ilp,metamis_warm}{
                        \addplot+[] table[x={\column},y=fraction, col sep=comma] {data/t100_perf.csv};
                    }
            \end{axis}
        \end{tikzpicture}
        \caption{Time limit 1 hour}\label{fig:vr-performance-100}
    \end{subfigure}%
    \begin{subfigure}[t]{0.5\textwidth}
        \centering
        \begin{tikzpicture}
            \begin{axis}[perf_left, xmin=0.975, name=p1]
                \foreach \column in {chils_cold,baseline,metamis_cold,chils_par,chils_warm,ilp,metamis_warm,chils_red}{
                        \addplot+[] table[x={\column},y=fraction, col sep=comma] {data/t100_best_perf.csv};
                    }
            \end{axis}
            \begin{axis}[perf_right, xmax=0.975, xmin=0.4, at={(p1.south east)},
                    legend style={
                            font=\footnotesize,
                            draw=none,
                            /tikz/every even column/.append style={column sep=5mm},
                            legend columns=4
                        },
                    legend to name=named
                ]
                \foreach \column in {chils_cold,baseline,metamis_cold,chils_par,chils_warm,ilp,metamis_warm,chils_red}{
                        \addplot+[] table[x={\column},y=fraction, col sep=comma] {data/t100_best_perf.csv};
                    }
                \legend{\chils-cold, \baseline, \metamis-cold, \chils-parallel, \chils-warm, \ilp, \metamis-warm, \chils-reduced}
            \end{axis}
        \end{tikzpicture}
        \caption{\centering Time limit 1 hour, added \\ parallel and reduced}\label{fig:vr-performance-full}
    \end{subfigure}
    \caption{Performance profiles for state-of-the-art comparison on solution quality on the vehicle routing instances with different time limits. In (d), results computed with parallel \chils{} and reduction rules are added.}\label{fig:vr-performance}
\end{figure}
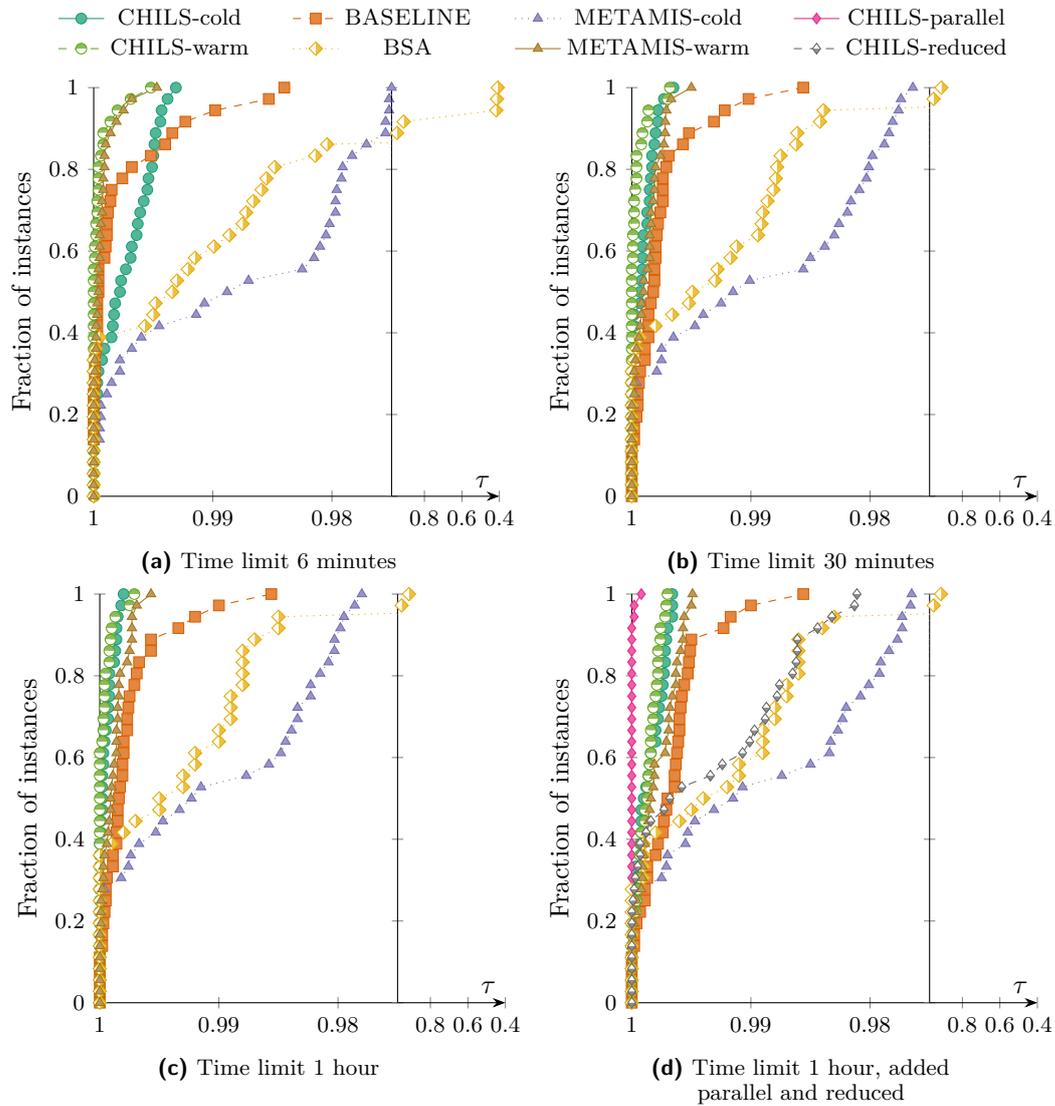

Next, we investigate the results on the second set of instances for the vehicle routing application. Note that the vertex weights for these instances are very large, which leads to only small percentage improvements despite significant differences between the solutions. Here, we only compare against the competitors \metamis{} and \ilp{}. In contrast to these competitors, \pils{} and \baseline{} are running on the full graph without exploiting the clique information that is also provided. On these instances, most of the solutions found are not optimal, evident by the significantly larger solutions we later find using our parallel \pils{}.
\renewcommand{\floatpagefraction}{.8}
We add a warm-start configuration for \pils{} where we start with the provided initial solution. The \metamis{} numbers for cold and warm-starts are taken from \cite{metamis}. In Figure~\ref{fig:vr-performance}, we present performance profiles to compare the solution quality achieved by the algorithms with different time limits. 
The first two show the state-of-the-art comparison of solution quality achieved with a 6 and 30 minutes time limit respectively. As expected, the configurations with the warm start perform best with the shortest time limit, see Figure~\ref{fig:vr-performance_less_time10}. However, if no initial solution is known, our two variants, \baseline{} and \pils{}, significantly outperform all competitor configurations. Note that \baseline{} is performing better than \pils{} on around 80\,\% of the instances here, since the configuration for \pils{} is optimized for running with a one-hour time limit, see Section~\ref{sec:parameter_tuning}. Compared to the second profile on the right, we can see that with more time, the \pils{} solutions improve most, even surpassing the \metamis{}-Warm results in some instances. Within this time limit \pils{}-Warm finds the best solutions on almost 60\,\% of the instances. Furthermore, both \pils{} configurations compute solutions that are at most 0.3\,\% worse than the best-found solution. On the other hand, on more than 38\,\% of the instances, \ilp{} and \metamis{}-Cold find solutions which are more than 1\,\% worse than the best solutions found on \hbox{the respective instances.}

In Figure~\ref{fig:vr-performance-100}, we present the state-of-the-art comparison where all algorithms have an hour time limit and run on the original instances. Here, we can see that \pils{} with and without the initial solution performs generally the best. The results shown here are very similar to the results computed within the shorter limit of 30 minutes; see \hbox{Figure~\ref{fig:vr-performance_less_time50}.}

For the profile in Figure~\ref{fig:vr-performance-full}, however, this changes. Here, we added a configuration called \pils{}-Parallel, which utilizes the warm start solution and runs in parallel with the configuration of ${P=64}$ and ${t=5}$ seconds, using 16 parallel threads, \ie without simultaneous multithreading. Additionally, we used \lnr~and run \pils{} on the reduced instances. With this configuration, \pils{}-Reduced, we can not compete with the other configurations on most instances since, for these large instances, the reduction rules do not work well; see also Table~\ref{tab:vr_graphs} in the Appendix. It takes considerable time to test even the fast reduction rules on these graphs. However, on \textit{MR-W-FN} or \textit{MW-D-01}, for example, using reduction rules definitely helps to find better solutions; see Table~\ref{tab:vr_results}. Considering all of our variants, we outperform all other schemes in all but two instances; see Figure~\ref{fig:vr-performance} and Table~\ref{tab:vr_results}. Furthermore, our approach is not depending on having good initial solutions as for example \metamis{}. Nevertheless, using a warm start solution can improve the performance.

\begin{obs}{}{obs:vehicle_routing_soa}
    On the vehicle routing instances without initial solutions, even our \baseline~approach outperforms the competitors. The \baseline~approach is especially good with a low time limit. Overall, \pils{} outperforms \emph{all} the other heuristics regarding solution quality after 6 minutes, 30 minutes, and 1 hour.
\end{obs}
\begin{landscape}
    \begin{table}[t]
        \centering
        \caption{Average solution weight $\w$ and time $t$ (in seconds). The \textbf{best} solutions among all algorithms are marked in bold. Here, we give the results of the four largest, non-reduced graphs of each class from our first set of instances. The results on all graphs are presented in Table~\ref{tab:soa_set1} in the Appendix. }\label{tab:soa_set1_sample} 
        \begin{adjustbox}{width=1.6\textwidth}
            
\begin{filecontents*}{data/datatool/soa_set1Full}
graph,weight,time,mmwisWeight,htwisWeight,baseWeight,chilsWeight,MaxWeight,mmwisWeightPerf,chilsWeightPerf,baseWeightPerf,hilsWeightPerf,htwisWeightPerf,mmwisTime,htwisTime,baseTime,chilsTime,MinTime,mmwisTimePerf,chilsTimePerf,baseTimePerf,hilsTimePerf,htwisTimePerf
fe-body-uniform,1678383,3471.02,1680182,1645650,1680182,1680182,1680182,1,1,1,0.9989,0.9794,2958.88,0.055079,1824.2176,165.0253,0.055079,53720.6558,2996.1564,33120.0203,63018.9364,1
fe-ocean-uniform,7081151,3384.46,7248581,6803672,7139076,7210830,7248581,1,0.9948,0.9849,0.9769,0.9386,3.38211,0.119556,16.7525,2863.4968,0.119556,28.2889,23951.0924,140.1226,28308.5751,1
fe-pwt-uniform,1175710,3580.37,1176685,1153600,1178479,1178708,1178708,0.9983,1,0.9998,0.9975,0.9787,3601.03,0.043087,49.2413,537.1301,0.043087,83575.7885,12466.1754,1142.8343,83096.2935,1
fe-rotor-uniform,2651471,3255.65,2596217,2591456,2663757,2665805,2665805,0.9739,1,0.9992,0.9946,0.9721,3601.17,0.309570,880.7024,3179.8309,0.30957,11632.8133,10271.767,2844.9217,10516.6844,1
fe-sphere-uniform,617139,3166.86,617816,608401,617738,617816,617816,1,1,0.9999,0.9989,0.9848,1.94976,0.015436,6.0504,125.5464,0.015436,126.3125,8133.3506,391.9668,205160.6634,1
mesh-blob-uniform,855004,252.711,855547,854484,855544,855547,855547,1,1,1,0.9994,0.9988,0.038847,0.006911,0.2647,45.3683,0.006911,5.621,6564.6506,38.3013,36566.4882,1
mesh-buddha-uniform,56086272,3600.05,57555880,57508556,57554770,57555764,57555880,1,1,1,0.9745,0.9992,4.25692,0.518566,3380.9415,2571.5242,0.518566,8.209,4958.914,6519.7902,6942.3179,1
mesh-bunny-uniform,3681110,3415.81,3686960,3682356,3686896,3686960,3686960,1,1,1,0.9984,0.9988,0.138293,0.036998,15.1256,167.0748,0.036998,3.7379,4515.7792,408.8221,92324.1797,1
mesh-dragonsub-uniform,31833085,3599.67,32213898,32163872,32212574,32213769,32213898,1,1,1,0.9882,0.9984,1.6175,0.280377,3409.6818,2857.1540,0.280377,5.769,10190.4008,12161.061,12838.6779,1
mesh-dragon-uniform,7946347,3485.87,7956530,7950526,7956496,7956526,7956530,1,1,1,0.9987,0.9992,0.268275,0.066966,32.3505,160.6388,0.066966,4.0061,2398.8113,483.0884,52054.3261,1
mesh-ecat-uniform,36122506,3600.03,36650298,36606394,36648982,36650039,36650298,1,1,1,0.9856,0.9988,2.59901,0.484511,3279.3517,2155.8664,0.484511,5.3642,4449.5716,6768.3741,7430.2338,1
mesh-fandisk-uniform,463049,3212.81,463288,462765,463280,463288,463288,1,1,1,0.9995,0.9989,0.017416,0.006573,2.3732,21.8195,0.006573,2.6496,3319.5649,361.0528,488788.9852,1
mesh-feline-uniform,2204823,3055.87,2207219,2204947,2207213,2207219,2207219,1,1,1,0.9989,0.999,0.101915,0.027325,1.6022,105.2000,0.027325,3.7297,3849.9543,58.6349,111834.2177,1
mesh-gameguy-uniform,2323182,3553.14,2325878,2324088,2325802,2325878,2325878,1,1,1,0.9988,0.9992,0.060395,0.019054,5.2518,42.0398,0.019054,3.1697,2206.3504,275.6272,186477.3801,1
mesh-gargoyle-uniform,1059102,2337.95,1059559,1058656,1059559,1059559,1059559,1,1,1,0.9996,0.9991,0.0454011,0.009649,1.3517,1.3517,0.009649,4.7053,140.0871,140.0871,242299.7202,1
mesh-turtle-uniform,14238437,3597.64,14263005,14247883,14262790,14262993,14263005,1,1,1,0.9983,0.9989,0.579221,0.138954,2753.3210,1665.3490,0.138954,4.1684,11984.8943,19814.6221,25890.8704,1
osm-alabama-AM2,174309,0.105678,174309,172797,174309,174309,174309,1,1,1,1,0.9913,0.031713,0.012066,0.1017,0.1012,0.012066,2.6283,8.3872,8.4286,8.7583,1
osm-alabama-AM3,185744,0.303255,185744,182667,185744,185744,185744,1,1,1,1,0.9834,28.9536,0.516524,14.4397,16.4507,0.303255,95.4761,54.2471,47.6157,1,1.7033
osm-california-AM2,47153,0.000287056,47153,45377,47153,47153,47153,1,1,1,1,0.9623,0.0257852,0.002292,0.0006,0.0013,0.000287056,89.8264,4.5287,2.0902,1,7.9845
osm-california-AM3,49365,0.083359,49365,48356,49365,49365,49365,1,1,1,1,0.9796,19.0425,0.044472,0.0460,0.0460,0.044472,428.1908,1.0344,1.0344,1.8744,1
osm-canada-AM3,86018,0.014745,86018,84503,86018,86018,86018,1,1,1,1,0.9824,2.06882,0.017174,0.0051,0.0049,0.0049,422.2082,1,1.0408,3.0092,3.5049
osm-colorado-AM3,54741,0.025454,54741,54435,54741,54741,54741,1,1,1,1,0.9944,0.072005,0.005051,0.0193,0.0193,0.005051,14.2556,3.821,3.821,5.0394,1
osm-district-of-columbia-AM1,196475,1.0464,196475,193364,196475,196475,196475,1,1,1,1,0.9842,1.05881,0.018109,0.0739,0.0351,0.018109,58.4687,1.9383,4.0808,57.7834,1
osm-district-of-columbia-AM2,209132,476.447,209132,198327,209132,209132,209132,1,1,1,1,0.9483,1643.2,2.402327,64.2834,160.1815,2.402327,684.0035,66.6776,26.7588,198.3273,1
osm-district-of-columbia-AM3,227652,939.718,227681,210461,227683,227683,227683,1,1,1,0.9999,0.9244,1890.32,330.184692,948.9367,258.4256,258.4256,7.3148,1,3.672,3.6363,1.2777
osm-florida-AM3,237333,1.93949,237333,234218,237333,237333,237333,1,1,1,1,0.9869,6.91146,0.196515,1.1378,1.1294,0.196515,35.1701,5.7471,5.7899,9.8694,1
osm-georgia-AM3,222652,0.0356388,222652,218573,222652,222652,222652,1,1,1,1,0.9817,5.9916,0.133135,0.1534,0.1512,0.0356388,168.1201,4.2426,4.3043,1,3.7357
osm-greenland-AM2,10718,0.130842,10718,10179,10718,10718,10718,1,1,1,1,0.9497,0.681488,0.058638,0.2353,0.2315,0.058638,11.622,3.948,4.0128,2.2314,1
osm-greenland-AM3,14012,3131.3,14011,12505,14012,14012,14012,0.9999,1,1,1,0.8924,1360.74,34.224670,180.8817,161.1345,34.22467,39.759,4.7081,5.2851,91.4925,1
osm-hawaii-AM2,125284,0.124023,125284,123173,125284,125284,125284,1,1,1,1,0.9832,0.307117,0.189414,1.0343,1.0448,0.124023,2.4763,8.4242,8.3396,1,1.5272
osm-hawaii-AM3,141041,3213.04,141070,134703,141074,141095,141095,0.9998,1,0.9999,0.9996,0.9547,2611.52,1423.033325,34.1224,160.0195,34.1224,76.5339,4.6896,1,94.1622,41.7038
osm-idaho-AM3,77145,115.883,77145,75527,77145,77145,77145,1,1,1,1,0.979,1795.12,52.148018,20.5178,12.0247,12.0247,149.2861,1,1.7063,9.6371,4.3367
osm-kansas-AM3,87976,2.25518,87976,87424,87976,87976,87976,1,1,1,1,0.9937,61.3928,3.247169,10.7643,36.0835,2.25518,27.223,16.0003,4.7731,1,1.4399
osm-kentucky-AM2,97397,0.454739,97397,97362,97397,97397,97397,1,1,1,1,0.9996,0.463053,0.934714,4.1913,4.2532,0.454739,1.0183,9.3531,9.2169,1,2.0555
osm-kentucky-AM3,100510,1292.9,99918,97906,100507,100511,100511,0.9941,1,1,1,0.9741,3684.53,3077.138184,35.3898,160.0263,35.3898,104.1128,4.5218,1,36.5331,86.9499
osm-louisiana-AM3,60024,0.0242,60024,59040,60024,60024,60024,1,1,1,1,0.9836,0.442016,0.025679,0.0268,0.0272,0.0242,18.2651,1.124,1.1074,1,1.0611
osm-maine-AM3,26734,0.000001,26734,26231,26734,26734,26734,1,1,1,1,0.9812,0.00188088,0.000532,0.0001,0.0001,0.000001,1880.88,100,100,1,532
osm-maryland-AM3,45496,0.0435009,45496,44539,45496,45496,45496,1,1,1,1,0.979,0.53903,0.068540,0.0341,0.0344,0.0341,15.8073,1.0088,1,1.2757,2.01
osm-massachusetts-AM2,140095,0.0454199,140095,139799,140095,140095,140095,1,1,1,1,0.9979,0.101676,0.024628,0.0256,0.0322,0.024628,4.1285,1.3075,1.0395,1.8442,1
osm-massachusetts-AM3,145866,0.409381,145866,144381,145866,145866,145866,1,1,1,1,0.9898,1292.74,2.220021,11.0560,11.5827,0.409381,3157.7919,28.2932,27.0066,1,5.4229
osm-mexico-AM2,94834,0.00469398,94834,94820,94834,94834,94834,1,1,1,1,0.9999,0.011462,0.006771,0.0005,0.0006,0.0005,22.924,1.2,1,9.388,13.542
osm-mexico-AM3,97663,0.0242021,97663,96700,97663,97663,97663,1,1,1,1,0.9901,8.94265,0.069216,0.0964,0.0973,0.0242021,369.4989,4.0203,3.9831,1,2.8599
osm-minnesota-AM3,32787,0.00957584,32787,32318,32787,32787,32787,1,1,1,1,0.9857,0.45735,0.036914,0.0056,0.0058,0.0056,81.6696,1.0357,1,1.71,6.5918
osm-montana-AM3,59822,0.102964,59822,59521,59753,59753,59822,1,0.9988,0.9988,1,0.995,3.36621,0.179889,0.0208,0.0205,0.0205,164.2054,1,1.0146,5.0226,8.7751
osm-nevada-AM3,52036,0.00872993,52036,51269,52036,52036,52036,1,1,1,1,0.9853,0.185488,0.014372,0.0032,0.0036,0.0032,57.965,1.125,1,2.7281,4.4912
osm-new-hampshire-AM3,116060,0.0429521,116060,115161,116060,116060,116060,1,1,1,1,0.9923,1.06467,0.015617,0.0094,0.0094,0.0094,113.2628,1,1,4.5694,1.6614
osm-new-york-AM2,14330,0.000638962,14330,14330,14330,14330,14330,1,1,1,1,1,0.0336602,0.004311,0.0008,0.0013,0.000638962,52.6795,2.0345,1.252,1,6.7469
osm-new-york-AM3,16268,0.0143349,16268,15293,16268,16268,16268,1,1,1,1,0.9401,60.5852,0.189482,0.1108,0.1116,0.0143349,4226.4125,7.7852,7.7294,1,13.2182
osm-north-carolina-AM3,49720,0.0507729,49720,49253,49720,49720,49720,1,1,1,1,0.9906,1277.21,0.521826,0.2345,0.2720,0.0507729,25155.3486,5.3572,4.6186,1,10.2776
osm-ohio-AM3,52634,0.000285149,52634,52199,52634,52634,52634,1,1,1,1,0.9917,0.708175,0.012258,0.0199,0.0198,0.000285149,2483.5262,69.4374,69.7881,1,42.9881
osm-oregon-AM3,175078,5.53855,175078,172813,175078,175078,175078,1,1,1,1,0.9871,1396.74,29.730036,64.3755,161.5105,5.53855,252.1851,29.1612,11.6232,1,5.3678
osm-pennsylvania-AM3,143870,0.0198829,143870,142472,143870,143870,143870,1,1,1,1,0.9903,1.32508,0.022415,0.0061,0.0111,0.0061,217.2262,1.8197,1,3.2595,3.6746
osm-puerto-rico-AM3,33590,0.00475907,33590,33556,33590,33590,33590,1,1,1,1,0.999,3.53969,0.031218,0.0074,0.0076,0.00475907,743.7777,1.597,1.5549,1,6.5597
osm-rhode-island-AM2,184596,0.140991,184596,179366,184596,184596,184596,1,1,1,1,0.9717,0.856335,0.463198,0.6049,0.1809,0.140991,6.0737,1.2831,4.2903,1,3.2853
osm-rhode-island-AM3,201753,1190.44,201771,190341,201771,201771,201771,1,1,1,0.9999,0.9434,1886.87,255.583374,366.2920,15.3432,15.3432,122.9776,1,23.8732,77.5875,16.6578
osm-tennessee-AM3,32276,0.000252962,32276,32232,32276,32276,32276,1,1,1,1,0.9986,0.0689001,0.003815,0.0002,0.0004,0.0002,344.5005,2,1,1.2648,19.075
osm-utah-AM2,95087,0.005795,95087,95080,95087,95087,95087,1,1,1,1,0.9999,0.00254798,0.001756,0.0035,0.0037,0.001756,1.451,2.1071,1.9932,3.3001,1
osm-utah-AM3,98847,0.0577071,98847,97754,98847,98847,98847,1,1,1,1,0.9889,1.34174,0.052415,0.1038,0.1021,0.052415,25.5984,1.9479,1.9803,1.101,1
osm-vermont-AM2,59310,1.90963,59310,57563,59284,59310,59310,1,1,0.9996,1,0.9705,1.07104,0.036269,0.0063,42.5171,0.0063,170.0063,6748.746,1,303.1159,5.757
osm-vermont-AM3,63304,667.25,63312,60518,63302,63302,63312,1,0.9998,0.9998,0.9999,0.9559,1317.01,4.890871,4.3775,4.4081,4.3775,300.8589,1.007,1,152.4272,1.1173
osm-virginia-AM2,295867,0.134254,295867,290535,295867,295867,295867,1,1,1,1,0.982,0.112363,0.028120,0.1211,0.1223,0.02812,3.9958,4.3492,4.3065,4.7743,1
osm-virginia-AM3,308305,1.33523,308305,300335,308305,308305,308305,1,1,1,1,0.9741,1297.87,1.375529,3.6287,3.6080,1.33523,972.0198,2.7022,2.7177,1,1.0302
osm-washington-AM2,305619,0.17322,305619,300195,305619,305619,305619,1,1,1,1,0.9823,0.171025,0.095822,0.1150,0.1154,0.095822,1.7848,1.2043,1.2001,1.8077,1
osm-washington-AM3,314288,546.599,314288,305019,314288,314288,314288,1,1,1,1,0.9705,1373.62,14.414819,115.0358,161.8076,14.414819,95.2922,11.2251,7.9804,37.9192,1
osm-west-virginia-AM3,47927,0.314467,47927,46344,47927,47927,47927,1,1,1,1,0.967,62.7196,0.322184,0.7402,0.7415,0.314467,199.4473,2.358,2.3538,1,1.0245
snap-as-skitter-uniform,123207288,3597.32,124157729,124141373,124141585,124111689,124157729,1,0.9996,0.9999,0.9923,0.9999,2828.19,1.169609,3545.4520,3491.3153,1.169609,2418.0645,2985.0277,3031.3139,3075.6603,1
snap-com-amazon,19266084,3590.7,19271031,19270078,19270569,19270939,19271031,1,1,1,0.9997,1,0.45293,0.177688,3571.0460,3556.4091,0.177688,2.549,20014.9087,20097.2829,20207.8925,1
snap-loc-gowalla-edges,12275221,3372.84,12276929,12276781,12275141,12275967,12276929,1,0.9999,0.9999,0.9999,1,4.98325,0.119892,1714.7353,3598.8512,0.119892,41.5645,30017.4424,14302.3329,28132.3191,1
snap-roadNet-CA-uniform,108370752,3600.18,111360828,111325524,111355248,111358047,111360828,1,1,0.9999,0.9731,0.9997,4.12136,0.722598,3496.8336,3579.3751,0.722598,5.7035,4953.4805,4839.2517,4982.2723,1
snap-roadNet-PA-uniform,60504456,3599.52,61731589,61710606,61730745,61731458,61731589,1,1,1,0.9801,0.9997,2.00781,0.414035,1264.1205,3347.5560,0.414035,4.8494,8085.2005,3053.173,8693.7578,1
snap-roadNet-TX-uniform,76796646,3600.03,78599946,78575460,78598317,78599552,78599946,1,1,1,0.9771,0.9997,2.38906,0.492106,3121.3741,3398.9519,0.492106,4.8548,6906.9507,6342.8897,7315.558,1
snap-soc-LiveJournal1-uniform,279386267,0.000001,284036239,283922214,283951703,283924551,284036239,1,0.9996,0.9997,0.9836,0.9996,1737.39,9.219021,3599.4601,3562.8544,0.000001,1737390000,3562854400,3599460100,1,9219021
snap-soc-pokec-relationships-uniform,82984397,3598.22,81112447,83920370,83971405,83982856,83982856,0.9658,1,0.9999,0.9881,0.9993,3637.19,39.062561,3369.2243,3520.9053,39.062561,93.1119,90.135,86.252,92.1143,1
snap-web-BerkStan-uniform,43694574,3599.51,43907482,43889843,43895988,43901482,43907482,1,0.9999,0.9997,0.9952,0.9996,4.2949,14.873417,2790.9626,3600.0813,4.2949,1,838.2224,649.8318,838.0894,3.463
snap-web-Google-uniform,56193137,3599.88,56326504,56323382,56320795,56321677,56326504,1,0.9999,0.9999,0.9976,0.9999,2.16259,0.744583,2539.9503,3580.0545,0.744583,2.9044,4808.1335,3411.2386,4834.7599,1
snap-web-NotreDame-uniform,26014630,3593.5,26016941,26013830,26014581,26016396,26016941,1,1,0.9999,0.9999,0.9999,0.986554,0.198811,3425.7288,3051.5395,0.198811,4.9623,15348.947,17231.0828,18074.9556,1
snap-web-Stanford-uniform,17779603,3595.78,17792930,17789430,17792392,17792474,17792930,1,1,1,0.9993,0.9998,1.30801,0.573333,2784.6966,3491.0494,0.573333,2.2814,6089.0432,4857.0318,6271.7129,1
ssmc-ca2010,16591401,3597.7,16867553,16792827,16867987,16869278,16869278,0.9999,1,0.9999,0.9835,0.9955,3601.26,0.548512,2650.7685,2613.1301,0.548512,6565.5081,4764.0345,4832.6536,6559.0179,1
ssmc-fl2010,8694303,3599.61,8743506,8719272,8742384,8743309,8743506,1,1,0.9999,0.9944,0.9972,183.766,0.396184,2941.6874,2832.0557,0.396184,463.84,7148.3344,7425.0535,9085.7026,1
ssmc-ga2010,4640168,3595.16,4644417,4639891,4644398,4644417,4644417,1,1,1,0.9991,0.999,61.2815,0.207957,579.1258,650.9023,0.207957,294.6835,3129.985,2784.8344,17287.997,1
ssmc-il2010,5965735,3599.95,5998539,5963974,5998017,5998520,5998539,1,1,0.9999,0.9945,0.9942,17.9288,0.392554,2516.4218,3314.4718,0.392554,45.6722,8443.3525,6410.3838,9170.5854,1
ssmc-nh2010,588843,2811.97,588996,587059,588996,588996,588996,1,1,1,0.9997,0.9967,0.436588,0.043277,15.5749,28.9120,0.043277,10.0882,668.0685,359.8886,64976.0843,1
ssmc-ri2010,458983,191.529,459275,457108,459270,459275,459275,1,1,1,0.9994,0.9953,0.490691,0.025208,5.2270,51.3621,0.025208,19.4657,2037.5317,207.3548,7597.9451,1
\end{filecontents*}
\DTLloaddb{DATA}{data/datatool/soa_set1Full}

\npthousandsep{\,}
\begin{tabular}{clrrrrrrrrrrrrrr}
        &&
        \multicolumn{2}{c}{\htwis} &&
        \multicolumn{2}{c}{\hils} &&
        \multicolumn{2}{c}{\mmwiss} &&
        \multicolumn{2}{c}{\baseline} &&
        \multicolumn{2}{c}{\pils} \\
        \cmidrule{3-4}\cmidrule{6-7}\cmidrule{9-10}\cmidrule{12-13}\cmidrule{15-16}\\[-1em]
        \multicolumn{2}{l}{\textbf{Instance}}& \multicolumn{1}{c}{$\w$} & \multicolumn{1}{c}{$t$}
        && \multicolumn{1}{c}{$\w$} & \multicolumn{1}{c}{$t$} 
        && \multicolumn{1}{c}{$\w$} & \multicolumn{1}{c}{$t$}
        && \multicolumn{1}{c}{$\w$} & \multicolumn{1}{c}{$t$} 
        && \multicolumn{1}{c}{$\w$} & \multicolumn{1}{c}{$t$} \\[-.2em]
        \cmidrule{3-4}\cmidrule{6-7}\cmidrule{9-10}\cmidrule{12-13}\cmidrule{15-16}\\[-1em]
        
    \DTLforeach*{DATA}{\instance=graph,\Whtwis=htwisWeight,\Thtwis=htwisTime,\Whils=weight,\Thils=time,\Wmmwis=mmwisWeight,\Tmmwis=mmwisTime,\Wchils=chilsWeight,\Tchils=chilsTime,\Wbase=baseWeight,\Tbase=baseTime,\max=MaxWeight}{%
    \sampleGraphs{\instance}%
    \ifthenelse{\equal{\result}{true}}{%
        \IfSubStr{\instance}{body}{\multirow{4}{*}{\rotatebox{90}{\textbf{fe}}}}{}%
        \IfSubStr{\instance}{ca2010}{\\[-.6em] \multirow{4}{*}{\rotatebox{90}{\textbf{ssmc}}}}{}%
        \IfSubStr{\instance}{skitter}{\\[-.6em] \multirow{4}{*}{\rotatebox{90}{\textbf{snap}}}}{}%
        \IfSubStr{\instance}{buddha}{\\[-.6em] \multirow{4}{*}{\rotatebox{90}{\textbf{mesh}}}}{}%
        \IfSubStr{\instance}{colum}{\\[-.6em] \multirow{4}{*}{\rotatebox{90}{\textbf{osm}}}}{}%
                & \niceGraphName{\instance}%
                & \ifnum \xintiiGtorEq{\Whtwis}{\max}=1 \relax \bfseries \fi 
                \numprint{\Whtwis} & \checkSmallTime{\Thtwis}%
                && \ifnum \xintiiGtorEq{\Whils}{\max}=1 \relax \bfseries \fi 
                \numprint{\Whils} & \checkSmallTime{\Thils}%
                && \ifnum \xintiiGtorEq{\Wmmwis}{\max}=1 \relax \bfseries \fi 
                \numprint{\Wmmwis} & \checkSmallTime{\Tmmwis}%
                && \ifnum \xintiiGtorEq{\Wbase}{\max}=1 \relax \bfseries \fi 
                \numprint{\Wbase} & \checkSmallTime{\Tbase}%
                && \ifnum \xintiiGtorEq{\Wchils}{\max}=1 \relax \bfseries \fi 
                \numprint{\Wchils} & \checkSmallTime{\Tchils}%
                \\
     }{}%
  }
\end{tabular}

\DTLgdeletedb{DATA}

        \end{adjustbox}
    \end{table}    
    \begin{table}[t]
    \centering
        \caption{Average solution weight $\w$ and time $t$ (in seconds). The \textbf{best} solutions among all algorithms are marked in bold. Here, we give the results of the graphs from Table~\ref{tab:soa_set1_sample} reduced with \lnr. Fully reduced graphs are omitted, reduction offset and reduction time is included in the results. All algorithms are run for one hour on the reduced instances, the \lnr~time was not limited. The results for all graphs are presented in Table~\ref{tab:kernel_soa_set1} in the Appendix. Detailed per instance results for the reduction time and offset are given in Table~\ref{tab:graphs}.}\label{tab:kernel_soa_set1_sample} 
        \begin{adjustbox}{width=1.6\textwidth}
            \begin{filecontents*}{data/datatool/kernel_data/soa_set1Kernel}
id,graph,redstyle,redconfig,timelimit,seed,gnnfilter,n,m,kn,km,offset,time,hilsWeight,mmwisWeight,htwisWeight,baseWeight,chilsWeight,mmwisTime,hilsTime,htwisTime,baseTime,chilsTime,MaxKWeight,MinKTime,MaxWeight,MinTime,mmwisWeightPerf,chilsWeightPerf,baseWeightPerf,hilsWeightPerf,htwisWeightPerf,mmwisTimePerf,chilsTimePerf,baseTimePerf,hilsTimePerf,htwisTimePerf
all-reductions-cyclicFast-early-cs-gnninitial-tight,fe-body-uniform,early-CS,all-reductions-cyclicFast,36000,0,initial-tight,45087,163734,395,3551,1672469,0.595169,1680182,1680182,1680066,1680182,1680182,60.629869,0.608746,0.596897,0.686669,0.684869,7713,0.001728,1680182,0.596897,1,1,1,1,0.985,34742.3032,51.9097,52.9514,7.8571,1
all-reductions-cyclicFast-early-cs-gnninitial-tight,fe-pwt-uniform,early-CS,all-reductions-cyclicFast,36000,0,initial-tight,36519,144794,15239,106943,813076,2.16788,1178273,1178735,1169768,1178735,1178735,3196.23788,1946.50788,2.206646,9.31458,9.83818,365659,0.038766,1178735,2.206646,1,1,1,0.9987,0.9755,82393.5923,197.8615,184.3548,50155.8066,1
all-reductions-cyclicFast-early-cs-gnninitial-tight,fe-rotor-uniform,early-CS,all-reductions-cyclicFast,36000,0,initial-tight,99617,662431,89294,690242,435335,5.05004,2653478,2617855,2599908,2664500,2666222,3606.23004,2764.03004,5.411118,347.15854,3251.99344,2230887,0.361078,2666222,5.411118,0.9783,1,0.9992,0.9943,0.9703,9973.4129,8992.3601,947.4643,7640.9529,1
all-reductions-cyclicFast-early-cs-gnninitial-tight,osm-alabama-AM3,early-CS,all-reductions-cyclicFast,36000,0,initial-tight,3504,309664,815,61637,173907,1.45062,185744,185744,184863,185744,185744,304.34962,1.4608241,1.523227,1.46952,1.46392,11837,0.0102041,185744,1.4608241,1,1,1,1,0.9256,29684.0486,1.3034,1.8522,1,7.1155
all-reductions-cyclicFast-early-cs-gnninitial-tight,osm-california-AM3,early-CS,all-reductions-cyclicFast,36000,0,initial-tight,587,27536,359,19705,35724,0.205125,49365,49365,48244,49365,49365,37.219725,0.20966592,0.225234,0.237625,0.233525,13641,0.00454092,49365,0.20966592,1,1,1,1,0.9178,8151.3438,6.2542,7.1571,1,4.4284
all-reductions-cyclicFast-early-cs-gnninitial-tight,osm-district-of-columbia-AM2,early-CS,all-reductions-cyclicFast,36000,0,initial-tight,13597,1609795,2361,274529,178323,5.76413,209132,209132,206445,209132,209132,1284.78413,5.971443,6.23557,5.89963,5.86483,30809,0.1007,209132,5.86483,1,1,1,1,0.9128,12701.291,1,1.3456,2.0587,4.6816
all-reductions-cyclicFast-early-cs-gnninitial-tight,osm-district-of-columbia-AM3,early-CS,all-reductions-cyclicFast,36000,0,initial-tight,46221,27729137,26728,13406329,133708,290.369,227586,227622,217052,227205,227656,3450.789,734.528,395.805752,350.7284,450.5029,93948,60.3594,227656,350.7284,0.9996,1,0.9952,0.9993,0.8871,52.36,2.653,1,7.3586,1.7468
all-reductions-cyclicFast-early-cs-gnninitial-tight,osm-florida-AM3,early-CS,all-reductions-cyclicFast,36000,0,initial-tight,2985,154043,658,41043,226379,0.52137,237333,237333,235834,237333,237333,14.06507,0.52934391,0.567098,0.71107,0.68807,10954,0.00797391,237333,0.52934391,1,1,1,1,0.8632,1698.5017,20.9057,23.7901,1,5.7347
all-reductions-cyclicFast-early-cs-gnninitial-tight,osm-georgia-AM3,early-CS,all-reductions-cyclicFast,36000,0,initial-tight,1680,74126,481,31388,207992,0.384872,222652,222652,220585,222652,222652,8.898342,0.38601379,0.423877,0.402672,0.397472,14660,0.00114179,222652,0.38601379,1,1,1,1,0.859,7456.2485,11.0353,15.5896,1,34.1613
all-reductions-cyclicFast-early-cs-gnninitial-tight,osm-greenland-AM3,early-CS,all-reductions-cyclicFast,36000,0,initial-tight,4986,3652361,3402,2141096,7572,27.3877,14012,14012,13532,14012,14012,1406.5577,416.0367,41.892963,28.3443,28.2494,6440,0.8617,14012,28.2494,1,1,1,1,0.9255,1600.5222,1,1.1101,451.0259,16.8333
all-reductions-cyclicFast-early-cs-gnninitial-tight,osm-hawaii-AM3,early-CS,all-reductions-cyclicFast,36000,0,initial-tight,28006,49444921,22849,38081016,96234,935.961,140963,141076,135846,140976,141095,4427.011,1068.106,1659.199831,1023.1769,1785.0451,44861,87.2159,141095,1023.1769,0.9996,1,0.9973,0.9971,0.883,40.0277,9.7354,1,1.5151,8.2925
all-reductions-cyclicFast-early-cs-gnninitial-tight,osm-idaho-AM3,early-CS,all-reductions-cyclicFast,36000,0,initial-tight,4064,3924080,2920,2413139,70960,40.0338,77145,77145,76227,77145,77145,1499.4838,52.1104,59.606775,185.0306,200.4611,6185,12.0766,77145,52.1104,1,1,1,1,0.8516,120.8494,13.2841,12.0064,1,1.6207
all-reductions-cyclicFast-early-cs-gnninitial-tight,osm-kansas-AM3,early-CS,all-reductions-cyclicFast,36000,0,initial-tight,2732,806912,1209,292997,84662,3.78747,87976,87976,87725,87976,87976,64.44307,4.051428,4.514076,5.48077,5.48297,3314,0.263958,87976,4.051428,1,1,1,1,0.9243,229.7926,6.4234,6.415,1,2.7527
all-reductions-cyclicFast-early-cs-gnninitial-tight,osm-kentucky-AM3,early-CS,all-reductions-cyclicFast,36000,0,initial-tight,19095,59533630,15943,53739580,91864,1510.41,100475,100193,99310,100511,100511,5125.45,2582.95,3421.082852,2167.947,2385.1401,8647,657.537,100511,2167.947,0.9632,1,1,0.9958,0.8611,5.4979,1.3303,1,1.6311,2.9058
all-reductions-cyclicFast-early-cs-gnninitial-tight,osm-massachusetts-AM3,early-CS,all-reductions-cyclicFast,36000,0,initial-tight,3703,551491,1626,340652,136331,2.98373,145866,145866,144574,145866,145866,1293.42373,3.092523,4.013044,3.57123,3.56883,9535,0.108793,145866,3.092523,1,1,1,1,0.8645,11861.4249,5.3781,5.4002,1,9.4612
all-reductions-cyclicFast-early-cs-gnninitial-tight,osm-mexico-AM3,early-CS,all-reductions-cyclicFast,36000,0,initial-tight,1096,47131,447,20443,86467,0.286542,97663,97663,96866,97663,97663,10.770342,0.28953105,0.306492,0.307042,0.308142,11196,0.00298905,97663,0.28953105,1,1,1,1,0.9288,3507.402,7.2264,6.8584,1,6.6744
all-reductions-cyclicFast-early-cs-gnninitial-tight,osm-montana-AM3,early-CS,all-reductions-cyclicFast,36000,0,initial-tight,837,69293,382,39859,55531,0.252079,59822,59822,59679,59822,59822,8.947519,0.252428045,0.31845,0.302479,0.297479,4291,0.000349045,59822,0.252428045,1,1,1,1,0.9667,24912.0887,130.0692,144.394,1,190.1503
all-reductions-cyclicFast-early-cs-gnninitial-tight,osm-new-york-AM3,early-CS,all-reductions-cyclicFast,36000,0,initial-tight,837,88728,526,50572,11447,0.416003,16268,16268,15666,16268,16268,62.237503,0.458698,0.493537,0.439403,0.435603,4821,0.0196,16268,0.435603,1,1,1,1,0.8751,3154.1582,1,1.1939,2.1783,3.9558
all-reductions-cyclicFast-early-cs-gnninitial-tight,osm-north-carolina-AM3,early-CS,all-reductions-cyclicFast,36000,0,initial-tight,1557,236739,997,133106,39783,0.872409,49720,49720,48897,49720,49720,1285.562409,51.583909,1.08438,1.446309,1.415909,9937,0.211971,49720,1.08438,1,1,1,1,0.9172,6060.6875,2.564,2.7074,239.2379,1
all-reductions-cyclicFast-early-cs-gnninitial-tight,osm-oregon-AM3,early-CS,all-reductions-cyclicFast,36000,0,initial-tight,5588,2912701,3159,1817045,162009,23.7533,175078,175078,173454,175078,175078,1410.4033,24.020772,35.27736,25.523,25.4579,13069,0.267472,175078,24.020772,1,1,1,1,0.8757,5184.281,6.373,6.6164,1,43.0851
all-reductions-cyclicFast-early-cs-gnninitial-tight,osm-rhode-island-AM2,early-CS,all-reductions-cyclicFast,36000,0,initial-tight,2866,295488,498,46183,174013,0.522558,184596,184596,184554,184596,184596,5.794418,0.52376797,0.572679,0.552658,0.550358,10583,0.00120997,184596,0.52376797,1,1,1,1,0.996,4357.0171,22.9758,24.8766,1,41.4233
all-reductions-cyclicFast-early-cs-gnninitial-tight,osm-rhode-island-AM3,early-CS,all-reductions-cyclicFast,36000,0,initial-tight,15124,12622219,12189,11225470,144106,153.055,201654,201767,197126,201480,201771,3171.015,184.2287,306.805549,167.5083,333.0726,57665,14.4533,201771,167.5083,0.9999,1,0.995,0.998,0.9194,208.8077,12.4551,1,2.1569,10.6377
all-reductions-cyclicFast-early-cs-gnninitial-tight,osm-vermont-AM3,early-CS,all-reductions-cyclicFast,36000,0,initial-tight,3436,1136164,1946,526463,52788,5.49778,63312,63312,61795,63312,63312,1317.27778,5.582173,7.007465,7.44468,7.41808,10524,0.084393,63312,5.582173,1,1,1,1,0.8559,15543.7062,22.7543,23.0694,1,17.8887
all-reductions-cyclicFast-early-cs-gnninitial-tight,osm-virginia-AM3,early-CS,all-reductions-cyclicFast,36000,0,initial-tight,6185,665903,2385,293068,272730,2.5103,308305,308305,306026,308305,308305,1294.6903,3.93287,3.050815,5.1284,5.1137,35575,0.540515,308305,3.050815,1,1,1,1,0.9359,2390.646,4.8165,4.8437,2.6319,1
all-reductions-cyclicFast-early-cs-gnninitial-tight,osm-washington-AM3,early-CS,all-reductions-cyclicFast,36000,0,initial-tight,10022,2346213,6671,1955568,264406,14.5531,314288,314288,307176,314288,314288,1406.8031,127.8021,22.977153,27.6024,130.21,49882,8.424053,314288,22.977153,1,1,1,1,0.8574,165.2708,13.7294,1.5491,13.4435,1
all-reductions-cyclicFast-early-cs-gnninitial-tight,osm-west-virginia-AM3,early-CS,all-reductions-cyclicFast,36000,0,initial-tight,1185,125620,864,98399,37678,0.516044,47927,47927,47074,47927,47927,68.217344,0.694478,0.698116,0.591944,0.586944,10249,0.0709,47927,0.586944,1,1,1,1,0.9168,954.8843,1,1.0705,2.5167,2.568
all-reductions-cyclicFast-early-cs-gnninitial-tight,snap-as-skitter-uniform,early-CS,all-reductions-cyclicFast,36000,0,initial-tight,1696415,11095298,1951,27277,124100687,3.43416,124157729,124157729,124156703,124157700,124157729,1300.24416,8.62787,3.458185,3.45396,33.48736,57042,0.0198,124157729,3.45396,1,1,0.9995,1,0.982,65495.4545,1517.8384,1,262.3086,1.2134
all-reductions-cyclicFast-early-cs-gnninitial-tight,snap-loc-gowalla-edges,early-CS,all-reductions-cyclicFast,36000,0,initial-tight,196591,950327,304,2938,12268682,0.327421,12276929,12276929,12276850,12276929,12276929,3.756341,0.32980495,0.328844,0.340921,0.338221,8247,0.001423,12276929,0.328844,1,1,1,1,0.9904,2409.6416,7.5896,9.487,1.6753,1
all-reductions-cyclicFast-early-cs-gnninitial-tight,snap-soc-LiveJournal1-uniform,early-CS,all-reductions-cyclicFast,36000,0,initial-tight,4847571,42851237,1715,28390,283975051,23.4735,284036239,284036239,284035711,284036239,284036239,1291.8035,23.645099,23.492865,23.5519,23.5522,61188,0.019365,284036239,23.492865,1,1,1,1,0.9914,65495.9979,4.064,4.0485,8.8613,1
all-reductions-cyclicFast-early-cs-gnninitial-tight,snap-soc-pokec-relationships-uniform,early-CS,all-reductions-cyclicFast,36000,0,initial-tight,1632803,22301964,782676,14234808,53362233,643.998,83762908,81480804,83913297,83986271,83991718,4311.478,4243.548,669.23922,4175.3388,3975.298,30629485,25.24122,83991718,669.23922,0.918,1,0.9998,0.9925,0.9974,145.2973,131.9786,139.9037,142.606,1
all-reductions-cyclicFast-early-cs-gnninitial-tight,snap-web-NotreDame-uniform,early-CS,all-reductions-cyclicFast,36000,0,initial-tight,325729,1090108,97,1357,26013618,0.622788,26016941,26016941,26016941,26016941,26016941,0.72531,0.62311511,0.623325,0.623688,0.623688,3323,0.00032711,26016941,0.62311511,1,1,1,1,1,313.4175,2.7514,2.7514,1,1.6416
all-reductions-cyclicFast-early-cs-gnninitial-tight,ssmc-fl2010,early-CS,all-reductions-cyclicFast,36000,0,initial-tight,484481,1173147,2045,8459,8707619,3.61887,8743506,8743506,8742642,8743506,8743506,1479.02887,37.95277,3.621486,3.65127,3.64927,35887,0.002616,8743506,3.621486,1,1,1,1,0.9759,563994.6483,11.6208,12.3853,13124.5795,1
\end{filecontents*}
\DTLloaddb{DATA}{data/datatool/kernel_data/soa_set1Kernel}
\npthousandsep{\,}

\begin{tabular}{clrrrrrrrrrrrrrr}
        &&
        \multicolumn{2}{c}{\htwis} &&
        \multicolumn{2}{c}{\hils} &&
        \multicolumn{2}{c}{\mmwiss} &&
        \multicolumn{2}{c}{\baseline} &&
        \multicolumn{2}{c}{\pils} \\
        \cmidrule{3-4}\cmidrule{6-7}\cmidrule{9-10}\cmidrule{12-13}\cmidrule{15-16}\\[-1em]
        \multicolumn{2}{l}{\textbf{Instance}} & \multicolumn{1}{c}{$\w$} & \multicolumn{1}{c}{$t$}
        && \multicolumn{1}{c}{$\w$} & \multicolumn{1}{c}{$t$} 
        && \multicolumn{1}{c}{$\w$} & \multicolumn{1}{c}{$t$}
        && \multicolumn{1}{c}{$\w$} & \multicolumn{1}{c}{$t$} 
        && \multicolumn{1}{c}{$\w$} & \multicolumn{1}{c}{$t$} \\[-.2em]
        \cmidrule{3-4}\cmidrule{6-7}\cmidrule{9-10}\cmidrule{12-13}\cmidrule{15-16}\\[-1em]
        
    \DTLforeach*{DATA}{\instance=graph,\Whtwis=htwisWeight,\Thtwis=htwisTime,\Whils=hilsWeight,\Thils=hilsTime,\Wmmwis=mmwisWeight,\Tmmwis=mmwisTime,\Wchils=chilsWeight,\Tchils=chilsTime,\Wbase=baseWeight,\Tbase=baseTime,\max=MaxWeight}{%
    \sampleGraphs{\instance}%
    \ifthenelse{\equal{\result}{true}}{%
        \IfSubStr{\instance}{body}{\multirow{4}{*}{\rotatebox{90}{\textbf{fe}}}}{}%
        \IfSubStr{\instance}{skitter}{\\[-.6em] \multirow{4}{*}{\rotatebox{90}{\textbf{snap}}}}{}%
        \IfSubStr{\instance}{fl2010}{\\[-.6em] \multirow{1}{*}{\rotatebox{90}{\textbf{ssmc}}}}{}%
        \IfSubStr{\instance}{colum}{\\[-.6em] \multirow{4}{*}{\rotatebox{90}{\textbf{osm}}}}{}%
                & \niceGraphName{\instance}%
                & \ifnum \xintiiGtorEq{\Whtwis}{\max}=1 \relax \bfseries \fi 
                \numprint{\Whtwis} & \checkSmallTime{\Thtwis}%
                && \ifnum \xintiiGtorEq{\Whils}{\max}=1 \relax \bfseries \fi 
                \numprint{\Whils} & \checkSmallTime{\Thils}%
                && \ifnum \xintiiGtorEq{\Wmmwis}{\max}=1 \relax \bfseries \fi 
                \numprint{\Wmmwis} & \checkSmallTime{\Tmmwis}%
                && \ifnum \xintiiGtorEq{\Wbase}{\max}=1 \relax \bfseries \fi 
                \numprint{\Wbase} & \checkSmallTime{\Tbase}%
                && \ifnum \xintiiGtorEq{\Wchils}{\max}=1 \relax \bfseries \fi 
                \numprint{\Wchils} & \checkSmallTime{\Tchils}%
                \\
     }{}%
  }
     \\[2em]
\end{tabular}

\DTLgdeletedb{DATA}

        \end{adjustbox}
    \end{table}
\end{landscape}

\renewcommand*{\arraystretch}{0.9}
\begin{landscape}
    \begin{figure}
        \footnotesize
        \setlength{\tabcolsep}{3.5pt}
        \begin{longtable}{llrlrlrlrlrlrlrlr}
            \caption{\normalsize Results for vehicle routing instances by Dong \etal~\cite{dong2021new}. The results for \metamis{} were taken from the paper by Dong \etal~\cite{metamis} since the code is not publicly available. Note that the authors of \metamis{} used the best out of four runs, while the results for \baseline{}, \chils{}, and \ilp{} were only run once. All the results show the best solution found after one hour.}\label{tab:vr_results}                                          \\
            \multicolumn{1}{c}{} & \multicolumn{1}{c}{} & \multicolumn{7}{c}{\textbf{Cold Start}} & \multicolumn{1}{c}{} & \multicolumn{3}{c}{\textbf{Warm Start}} & \multicolumn{1}{c}{} & \multicolumn{1}{c}{\textbf{Reduced}} & \multicolumn{1}{c}{} & \multicolumn{1}{c}{\textbf{Parallel}}                                                                                                                                                                                                                   \\ \cline{3-9} \cline{11-13} \cline{15-15} \cline{17-17} \\[-.5em]
            \textbf{Instance}    & \multicolumn{1}{c}{} & \multicolumn{1}{c}{\metamis}            & \multicolumn{1}{c}{} & \multicolumn{1}{c}{\ilp}                & \multicolumn{1}{c}{} & \multicolumn{1}{c}{\baseline}        & \multicolumn{1}{c}{} & \multicolumn{1}{c}{\chils}            & \multicolumn{1}{c}{} & \multicolumn{1}{c}{\metamis} & \multicolumn{1}{c}{} & \multicolumn{1}{c}{\chils} & \multicolumn{1}{c}{} & \multicolumn{1}{c}{\chils} & \multicolumn{1}{c}{} & \multicolumn{1}{c}{\chils} \\ \cline{1-1} \cline{3-3} \cline{5-5} \cline{7-7} \cline{9-9} \cline{11-11} \cline{13-13} \cline{15-15} \cline{17-17} \\[-.5em]
            CR-S-L-1             &                      & 5\,588\,489                             &                      & 5\,639\,292                             &                      & 5\,694\,508                          &                      & 5\,698\,608                           &                      & 5\,692\,891                  &                      & 5\,701\,015                &                      & 5\,610\,201                &                      & \textbf{5\,718\,230}       \\
            CR-S-L-2             &                      & 5\,691\,892                             &                      & 5\,729\,194                             &                      & 5\,785\,140                          &                      & 5\,798\,470                           &                      & 5\,784\,034                  &                      & 5\,799\,458                &                      & 5\,715\,463                &                      & \textbf{5\,813\,362}       \\
            CR-S-L-4             &                      & 5\,681\,336                             &                      & 5\,725\,525                             &                      & 5\,781\,519                          &                      & 5\,791\,447                           &                      & 5\,777\,081                  &                      & 5\,792\,758                &                      & 5\,698\,425                &                      & \textbf{5\,806\,975}       \\
            CR-S-L-6             &                      & 3\,859\,513                             &                      & 3\,900\,370                             &                      & 3\,936\,944                          &                      & 3\,938\,946                           &                      & 3\,936\,137                  &                      & 3\,946\,157                &                      & 3\,890\,733                &                      & \textbf{3\,952\,205}       \\
            CR-S-L-7             &                      & 1\,989\,879                             &                      & 2\,006\,810                             &                      & 2\,017\,624                          &                      & 2\,021\,238                           &                      & 2\,019\,428                  &                      & 2\,022\,339                &                      & 2\,006\,857                &                      & \textbf{2\,025\,681}       \\
            CR-T-C-1             &                      & 4\,654\,419                             &                      & 4\,717\,754                             &                      & 4\,742\,508                          &                      & 4\,744\,119                           &                      & 4\,743\,040                  &                      & 4\,750\,570                &                      & 4\,702\,548                &                      & \textbf{4\,760\,428}       \\
            CR-T-C-2             &                      & 4\,874\,346                             &                      & 4\,934\,488                             &                      & 4\,969\,512                          &                      & 4\,975\,069                           &                      & 4\,968\,952                  &                      & 4\,976\,613                &                      & 4\,918\,257                &                      & \textbf{4\,987\,562}       \\
            CR-T-D-4             &                      & 4\,817\,281                             &                      & 4\,875\,268                             &                      & 4\,912\,984                          &                      & 4\,919\,079                           &                      & 4\,911\,646                  &                      & 4\,922\,752                &                      & 4\,864\,858                &                      & \textbf{4\,932\,826}       \\
            CR-T-D-6             &                      & 2\,970\,011                             &                      & 3\,009\,020                             &                      & 3\,024\,044                          &                      & 3\,027\,211                           &                      & 3\,024\,523                  &                      & 3\,029\,782                &                      & 2\,999\,231                &                      & \textbf{3\,033\,189}       \\
            CR-T-D-7             &                      & 1\,440\,281                             &                      & 1\,453\,990                             &                      & 1\,460\,328                          &                      & 1\,461\,099                           &                      & 1\,460\,240                  &                      & 1\,460\,584                &                      & 1\,451\,070                &                      & \textbf{1\,462\,239}       \\
            CW-S-L-1             &                      & 1\,634\,950                             &                      & 1\,645\,459                             &                      & 1\,660\,063                          &                      & 1\,660\,472                           &                      & 1\,660\,815                  &                      & 1\,662\,580                &                      & 1\,646\,548                &                      & \textbf{1\,663\,763}       \\
            CW-S-L-2             &                      & 1\,708\,820                             &                      & 1\,713\,535                             &                      & 1\,737\,052                          &                      & 1\,738\,307                           &                      & 1\,738\,128                  &                      & 1\,739\,670                &                      & 1\,723\,914                &                      & \textbf{1\,743\,532}       \\
            CW-S-L-4             &                      & 1\,725\,591                             &                      & 1\,730\,641                             &                      & 1\,751\,204                          &                      & 1\,754\,409                           &                      & 1\,753\,803                  &                      & 1\,756\,602                &                      & 1\,735\,633                &                      & \textbf{1\,759\,452}       \\
            CW-S-L-6             &                      & 1\,158\,925                             &                      & 1\,162\,552                             &                      & 1\,175\,335                          &                      & 1\,177\,272                           &                      & 1\,177\,156                  &                      & 1\,176\,354                &                      & 1\,166\,774                &                      & \textbf{1\,178\,467}       \\
            CW-S-L-7             &                      & 587\,288                                &                      & 588\,279                                &                      & 594\,312                             &                      & 594\,598                              &                      & 593\,891                     &                      & 593\,807                   &                      & 590\,834                   &                      & \textbf{594\,757}          \\
            CW-T-C-1             &                      & 1\,317\,775                             &                      & 1\,325\,862                             &                      & 1\,336\,232                          &                      & 1\,338\,085                           &                      & 1\,338\,064                  &                      & 1\,340\,085                &                      & 1\,323\,170                &                      & \textbf{1\,341\,888}       \\
            CW-T-C-2             &                      & 931\,802                                &                      & 935\,238                                &                      & 944\,902                             &                      & 946\,217                              &                      & 945\,886                     &                      & 945\,913                   &                      & 935\,907                   &                      & \textbf{947\,644}          \\
            CW-T-D-4             &                      & 457\,185                                &                      & 459\,575                                &                      & 460\,955                             &                      & 460\,734                              &                      & 461\,056                     &                      & 461\,108                   &                      & 459\,543                   &                      & \textbf{461\,475}          \\
            CW-T-D-6             &                      & 457\,790                                &                      & 459\,238                                &                      & 461\,088                             &                      & 461\,354                              &                      & 461\,312                     &                      & 461\,101                   &                      & 460\,227                   &                      & \textbf{461\,709}          \\
            MR-D-03              &                      & 1\,754\,110\,286                        &                      & 1\,757\,141\,345                        &                      & 1\,752\,894\,190                     &                      & 1\,758\,762\,733                      &                      & 1\,757\,227\,519             &                      & 1\,758\,429\,721           &                      & 1\,758\,775\,738           &                      & \textbf{1\,759\,255\,435}  \\
            MR-D-05              &                      & 1\,786\,342\,921                        &                      & 1\,788\,812\,740                        &                      & 1\,781\,868\,583                     &                      & 1\,789\,601\,100                      &                      & 1\,787\,849\,777             &                      & 1\,789\,537\,256           &                      & 1\,789\,944\,187           &                      & \textbf{1\,790\,776\,639}  \\
            MR-D-FN              &                      & 1\,797\,573\,192                        &                      & 1\,801\,983\,754                        &                      & 1\,794\,240\,970                     &                      & 1\,801\,499\,264                      &                      & 1\,799\,452\,160             &                      & 1\,800\,375\,890           &                      & 1\,801\,727\,328           &                      & \textbf{1\,802\,925\,854}  \\
            MR-W-FN              &                      & 5\,358\,386\,615                        &                      & 5\,386\,842\,780                        &                      & 5\,385\,799\,979                     &                      & 5\,385\,799\,979                      &                      & \textbf{5\,386\,842\,781}    & \textbf{}            & \textbf{5\,386\,842\,781}  & \textbf{}            & \textbf{5\,386\,842\,781}  &                      & \textbf{5\,386\,842\,781}  \\
            MT-D-01              &                      & \textbf{238\,166\,485}                  & \textbf{}            & \textbf{238\,166\,485}                  & \textbf{}            & \textbf{238\,166\,485}               & \textbf{}            & \textbf{238\,166\,485}                & \textbf{}            & \textbf{238\,166\,485}       & \textbf{}            & \textbf{238\,166\,485}     & \textbf{}            & \textbf{238\,166\,485}     & \textbf{}            & \textbf{238\,166\,485}     \\
            MT-D-200             &                      & 287\,048\,909                           &                      & \textbf{287\,155\,108}                  &                      & 287\,042\,596                        &                      & 287\,086\,442                         &                      & 287\,048\,081                &                      & 287\,042\,596              &                      & 287\,086\,442              &                      & 287\,086\,442              \\
            MT-D-FN              &                      & \textbf{290\,866\,943}                  & \textbf{}            & \textbf{290\,866\,943}                  & \textbf{}            & \textbf{290\,866\,943}               & \textbf{}            & \textbf{290\,866\,943}                &                      & 290\,771\,450                &                      & 290\,771\,450              &                      & \textbf{290\,866\,943}     &                      & \textbf{290\,866\,943}     \\
            MT-W-01              &                      & \textbf{312\,121\,568}                  & \textbf{}            & \textbf{312\,121\,568}                  & \textbf{}            & \textbf{312\,121\,568}               & \textbf{}            & \textbf{312\,121\,568}                & \textbf{}            & \textbf{312\,121\,568}       & \textbf{}            & \textbf{312\,121\,568}     & \textbf{}            & \textbf{312\,121\,568}     &                      & \textbf{312\,121\,568}     \\
            MT-W-200             &                      & 383\,961\,323                           &                      & 384\,056\,011                           &                      & 383\,974\,084                        &                      & 384\,052\,157                         &                      & 383\,985\,408                &                      & 384\,052\,017              &                      & 384\,052\,157              &                      & \textbf{384\,056\,012}     \\
            MT-W-FN              &                      & 390\,854\,593                           &                      & 390\,869\,890                           &                      & \textbf{390\,869\,891}               & \textbf{}            & \textbf{390\,869\,891}                & \textbf{}            & \textbf{390\,869\,891}       & \textbf{}            & \textbf{390\,869\,891}     & \textbf{}            & \textbf{390\,869\,891}     &                      & \textbf{390\,869\,891}     \\
            MW-D-01              &                      & 476\,334\,711                           &                      & 476\,298\,607                           &                      & 475\,180\,516                        &                      & 476\,328\,138                         &                      & 475\,987\,082                &                      & 476\,120\,423              &                      & \textbf{476\,440\,656}     &                      & 476\,360\,408              \\
            MW-D-20              &                      & 525\,124\,575                           &                      & 526\,857\,183                           &                      & 523\,423\,307                        &                      & 526\,883\,302                         &                      & 525\,486\,034                &                      & 526\,333\,489              &                      & 526\,648\,329              &                      & \textbf{527\,498\,481}     \\
            MW-D-40              &                      & 536\,520\,199                           &                      & \textbf{539\,036\,121}                  &                      & 533\,671\,730                        &                      & 538\,302\,190                         &                      & 536\,735\,155                &                      & 537\,485\,389              &                      & 538\,409\,586              &                      & 538\,596\,069              \\
            MW-D-FN              &                      & 541\,918\,916                           &                      & 545\,554\,192                           &                      & 541\,193\,604                        &                      & 544\,451\,893                         &                      & 543\,857\,187                &                      & 544\,205\,239              &                      & 544\,246\,489              &                      & \textbf{545\,711\,017}     \\
            MW-W-01              &                      & \textbf{1\,270\,305\,952}               &                      & \textbf{1\,270\,305\,952}               &                      & \textbf{1\,270\,305\,952}            &                      & 1\,270\,235\,200                      &                      & 1\,269\,344\,846             &                      & 1\,269\,344\,846           &                      & \textbf{1\,270\,305\,952}  &                      & \textbf{1\,270\,305\,952}  \\
            MW-W-05              &                      & \textbf{1\,328\,958\,047}               &                      & 1\,326\,236\,043                        &                      & 1\,327\,478\,708                     &                      & 1\,328\,043\,785                      &                      & \textbf{1\,328\,958\,047}    &                      & \textbf{1\,328\,958\,047}  &                      & 1\,328\,043\,785           &                      & \textbf{1\,328\,958\,047}  \\
            MW-W-10              &                      & 1\,342\,899\,725                        &                      & 1\,280\,286\,209                        &                      & 1\,340\,268\,013                     &                      & 1\,342\,808\,634                      &                      & \textbf{1\,342\,915\,691}    &                      & \textbf{1\,342\,915\,691}  &                      & 1\,342\,809\,954           &                      & \textbf{1\,342\,915\,691}  \\
            MW-W-FN              &                      & \textbf{1\,350\,818\,543}               &                      & 1\,235\,306\,258                        &                      & 1\,331\,333\,002                     &                      & \textbf{1\,350\,818\,543}             &                      & \textbf{1\,350\,818\,543}    &                      & \textbf{1\,350\,818\,543}  &                      & \textbf{1\,350\,818\,543}  &                      & \textbf{1\,350\,818\,543}
        \end{longtable}
    \end{figure}
\end{landscape}

\subsection{Parallel Scalability}

For the parallel results shown in Figure~\ref{fig:vr-performance-full}, we used the same machine as the other heuristics to ensure fairness between each program with regard to solution quality. To gain more detailed insight into the parallel scalability of our approach, we utilize a larger machine with an AMD EPYC 9754 128-core processor for our scalability experiments. Furthermore, \chils{} as defined in Section~\ref{sec:pils} relies on wall time to alternate between local search on the full graph and the \chilscore{}. This introduces variance between runs and makes it difficult to conduct scalability experiments. Therefore, we change the implementation for this section to perform a fixed number of local search iterations instead, where one iteration refers to the while loop starting on Line 2 in Algorithm~\ref{alg:ls}. We also set a fixed number of \chils{} iterations, referring to the while loop starting on Line 2 in Algorithm~\ref{alg:pils}. By removing the wall-clock measures and using the same random seed, we ensure that the parallel and sequential versions perform the exact same computations and reach the same solution \hbox{in the end.}


\pgfplotsset{
    /pgfplots/colormap={mycm}{rgb255=(27, 158, 119) rgb255=(253, 198, 18) rgb255=(231, 41, 187)}
}

\pgfplotsset{select coords between index/.style 2 args={
            x filter/.code={
                    \ifnum\coordindex<#1\def\pgfmathresult{}\fi
                    \ifnum\coordindex>#2\def\pgfmathresult{}\fi
                }
        }}

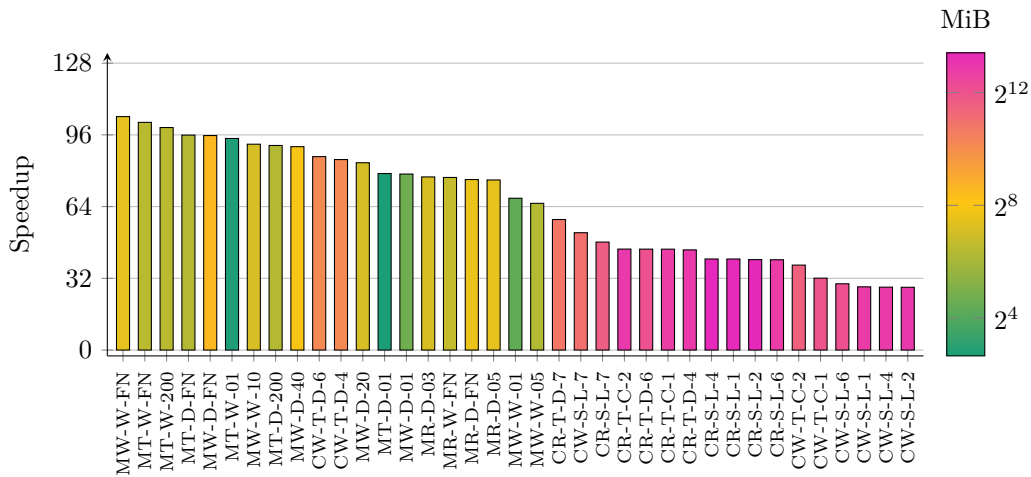
\begin{figure}[!t]
    \centering
    \begin{tikzpicture}
        \pgfplotstableread[col sep=comma]{data/parallel_speedup.csv}{\loadedtable}%
        \pgfplotstablesort[sort key = speedup, sort cmp = float >]{\sortedtable}{\loadedtable}

        \begin{axis}[
                colorbar,
                colorbar style={
                        title={MiB},
                        ytick={4,8,12},
                        yticklabels={$2^4$, $2^8$, $2^{12}$}
                    },
                scatter,
                scatter src=y,
                only marks,
                clip mode=individual,
                visualization depends on=\thisrow{speedup}\as\my,
                scatter/@pre marker code/.append code={
                        \pgfkeys{/pgf/fpu=true,/pgf/fpu/output format=fixed}
                        \fill [draw=black] (0.3,0) rectangle (-0.3,-\my);
                        \pgfplotsset{mark=none}
                    },
                point meta=\thisrow{size-log},
                ymin=0, ymax=130,
                width=.88\textwidth,
                height=.4\textwidth,
                ymajorgrids,
                ytick={0,32,64,96,128},
                xtick=data,
                xticklabels={MW-W-FN,MT-W-FN,MT-W-200,MT-D-FN,MW-D-FN,MT-W-01,MW-W-10,MT-D-200,MW-D-40,CW-T-D-6,CW-T-D-4,MW-D-20,MT-D-01,MW-D-01,MR-D-03,MR-W-FN,MR-D-FN,MR-D-05,MW-W-01,MW-W-05,CR-T-D-7,CW-S-L-7,CR-S-L-7,CR-T-C-2,CR-T-D-6,CR-T-C-1,CR-T-D-4,CR-S-L-4,CR-S-L-1,CR-S-L-2,CR-S-L-6,CW-T-C-2,CW-T-C-1,CW-S-L-6,CW-S-L-1,CW-S-L-4,CW-S-L-2},
                x tick label style={rotate=90,anchor=east,font=\scriptsize},
                axis lines=middle,
                axis x line*=bottom,
                axis y line*=left,
                x axis line style={-},
                ylabel={Speedup},
                ylabel near ticks,
                enlargelimits=0.02
            ]
            \addplot [select coords between index={0}{36}] table[col sep=comma, x expr=\coordindex, y=speedup] {\sortedtable};
        \end{axis}
    \end{tikzpicture}
    \caption{Bar chart showing the speedup on each instance using 128 threads compared to the sequential CHILS. The colors indicate the amount of allocated memory for each instance. The total amount of L3 cache for this machine is 256 MiB, and smaller instances fit entirely in the cache.}\label{fig:par}
\end{figure}

The configuration we use for these experiments consists of \numprint{1000} local search iterations, \numprint{10} \chils~iterations, and $P=128$. Depending on the instance, this ratio of local search on the full graph and \chilscore{} is reasonably close to the wall-clock version with $t_G=t_C=10$. Each run was repeated five times, and the best measurement is used. We use the best measurement because, in this configuration, there is no randomness between runs. Any difference we observe in execution time is solely due to factors outside our control, such as fluctuations in clock speed and other programs running on the machine. As such, the best measure is the closest to the true execution time. 

The speedup numbers for each instance are shown in Figure~\ref{fig:par}. The best scaling instance reaches a speedup of 104, while the worst is still 28 times faster than the sequential program. Figure~\ref{fig:par} also shows the amount of memory allocated for each instance. This includes the graph, which we store using the compressed sparse row format, and additional data structures required by the heuristic. All memory allocations are done up front before starting the heuristic, and the appropriate thread initializes any thread-local data. The CPU we use has a combined L3 cache of 256 MiB. There is a clear drop in speedup around the point when the data no longer fits in the cache. This indicates that the memory bandwidth is the main bottleneck for larger instances. In terms of load balancing, there could be variations in how much work it takes to perform \numprint{1000} local search iterations for each solution. This is because each solution has a different random seed and $m_q$ parameter. Note that this is not an issue for the version presented in Section~\ref{sec:pils}, since that version uses wall time instead of local search iterations. Table~\ref{tab:graphs} in the Appendix gives detailed execution time, speedup, and the amount of memory used \hbox{for each instance.}

\section{Conclusion and Future Work}

We have introduced \lnr{}, a preprocessing algorithm for the MWIS problem that combines Graph Neural Networks (GNNs) with a large collection of reduction rules to reduce further and faster than previously possible. In addition to using known data reduction rules, we also introduce seven new reduction rules. The GNNs are trained to predict where we can apply costly reduction rules to speed up the reduction process. Combined, this strikes a good balance between speed and quality at the preprocessing stage. In addition to exact preprocessing, we also introduce a new heuristic called \pils{} (Concurrent Hybrid Iterated Local Search) that expands on the \hils{} heuristic. This new heuristic outperforms all known heuristics across a wide range of test instances in a sequential environment. As an added benefit, \pils{} can also leverage the power offered by multicore processors. Letting \pils{} use all 16 cores available on our test machine significantly improves the solution quality on the hardest instances in our dataset.

The vehicle routing instances by Dong~\etal~\cite{dong2021new} offer a significant challenge for practical MWIS algorithms. Our result marks the third iteration of improvements to this dataset, after \metamis{}~\cite{metamis} and \textsc{BSA}~\cite{haller2024bregman}, and yet, we are still far away from optimal solutions on these instances. This is evident by the significant uplift in solution quality by running \pils{} in parallel.
Using another 128-core machine for scalability experiments, we show that \pils{} reaches speedups up to 104 for \hbox{the vehicle routing instances.}

There are several directions for future research, including finding efficient data reductions that work on these hard instances or using the clique information in the \pils{} heuristic. If the use of clique information leads to improvements, then a natural continuation would be to compute clique covers for other hard instances.

We have also introduced a supervised learning dataset for MWIS reductions. The models we use in \lnr{} are trained on this dataset. In this work, we only considered the most common GNN architectures and only the application of reduction rule screening. Besides trying more complicated architectures, there are other directions for further work starting from this dataset. One promising direction is to use GNN models to reduce the graph directly. Even though this would no longer be exact preprocessing it could lead to \hbox{a powerful heuristic.}

\pils{} is based on the proposed metaheuristic \textsc{Concurrent Difference-Core Heuristic}. This metaheuristic could lead to improvements in heuristics for other problems as well. As a metaheuristic, it only requires that solutions can be compared to find the \textsc{Difference-Core}; otherwise, any heuristic method can be used internally. Examples of problems to try include \textsc{Vertex Cover}, \textsc{Dominating set}, \textsc{Graph Coloring}, and \textsc{Connectivity Augmentation}. As an added benefit, the \textsc{Concurrent Difference-Core Heuristic} is trivially parallelizable, which could enable improvements in the parallel setting too.

\newpage

\bibliography{bib}
\newpage

\appendix

\npthousandsep{\,}

\DTLloaddb{DATA}{data/datatool/soa_set1Full}
\begin{landscape}
    \setlength{\tabcolsep}{4.5pt}
\section{State-of-the-Art Results}
\footnotesize{
    \renewcommand{\arraystretch}{1.14}
\begin{longtable}{rlrrrrrrrrrrrrrr}\caption{Average solution weight $\w$ and time $t$ (in seconds) required to compute it for the \emph{non reduced} instances from set one.
The \textbf{best} solutions among all algorithms are marked bold.}\label{tab:soa_set1} \\
        &&
        \multicolumn{2}{c}{\htwis} &&
        \multicolumn{2}{c}{\hils} &&
        \multicolumn{2}{c}{\mmwiss} &&
        \multicolumn{2}{c}{\baseline} &&
        \multicolumn{2}{c}{\pils} \\
        \cmidrule{3-4}\cmidrule{6-7}\cmidrule{9-10}\cmidrule{12-13}\cmidrule{15-16}
        \multicolumn{2}{c}{\bfseries Instance}
        & \multicolumn{1}{c}{$\w$} & \multicolumn{1}{c}{$t$} &
        & \multicolumn{1}{c}{$\w$} & \multicolumn{1}{c}{$t$} &
        & \multicolumn{1}{c}{$\w$} & \multicolumn{1}{c}{$t$} &
        & \multicolumn{1}{c}{$\w$} & \multicolumn{1}{c}{$t$} &
        & \multicolumn{1}{c}{$\w$} & \multicolumn{1}{c}{$t$} \\[-.2em]
        \cmidrule{3-4}\cmidrule{6-7}\cmidrule{9-10}\cmidrule{12-13}\cmidrule{15-16}
    \endfirsthead
    
    \multicolumn{16}{c}%
    {\tablename\ \thetable\ -- \textit{Continued from previous page}} \\    
        &&
        \multicolumn{2}{c}{\htwis} &&
        \multicolumn{2}{c}{\hils} &&
        \multicolumn{2}{c}{\mmwiss} &&
        \multicolumn{2}{c}{\baseline} &&
        \multicolumn{2}{c}{\pils} \\        
        \cmidrule{3-4}\cmidrule{6-7}\cmidrule{9-10}\cmidrule{12-13}\cmidrule{15-16}\\[-1em]
        \multicolumn{2}{c}{\bfseries Instance}
        & \multicolumn{1}{c}{$\w$} & \multicolumn{1}{c}{$t$}&
        & \multicolumn{1}{c}{$\w$} & \multicolumn{1}{c}{$t$} &
        & \multicolumn{1}{c}{$\w$} & \multicolumn{1}{c}{$t$}&
        & \multicolumn{1}{c}{$\w$} & \multicolumn{1}{c}{$t$} &
        & \multicolumn{1}{c}{$\w$} & \multicolumn{1}{c}{$t$} \\[-.2em]
        \cmidrule{3-4}\cmidrule{6-7}\cmidrule{9-10}\cmidrule{12-13}\cmidrule{15-16}
        
    \endhead

    \midrule \multicolumn{16}{r}{\textit{Continued on next page}} \\
    \endfoot
    
    \endlastfoot

    \DTLforeach*{DATA}{\instance=graph,\Whtwis=htwisWeight,\Thtwis=htwisTime,\Whils=weight,\Thils=time,\Wmmwis=mmwisWeight,\Tmmwis=mmwisTime,\Wchils=chilsWeight,\Tchils=chilsTime,\Wbase=baseWeight,\Tbase=baseTime,\max=MaxWeight}{%
    \IfSubStr{\instance}{alabama-AM2}{\\[-.6em]\textbf{osm}}{}%
    \IfSubStr{\instance}{ca2010}{\\[-.6em]\textbf{ssmc}}{}%
    \IfSubStr{\instance}{skitter}{\\[-.6em]\textbf{snap}}{}%
    \IfSubStr{\instance}{body}{\\[-.6em]\textbf{fe}}{}%
    \IfSubStr{\instance}{blob}{\\[-.6em]\textbf{mesh}}{}%
            &\niceGraphName{\instance}%
            & \ifnum \xintiiGtorEq{\Whtwis}{\max}=1 \relax \bfseries \fi 
            \numprint{\Whtwis} & \checkSmallTime{\Thtwis}%
            && \ifnum \xintiiGtorEq{\Whils}{\max}=1 \relax \bfseries \fi 
            \numprint{\Whils} & \checkSmallTime{\Thils}%
            && \ifnum \xintiiGtorEq{\Wmmwis}{\max}=1 \relax \bfseries \fi 
            \numprint{\Wmmwis} & \checkSmallTime{\Tmmwis}%
            && \ifnum \xintiiGtorEq{\Wbase}{\max}=1 \relax \bfseries \fi 
            \numprint{\Wbase} & \checkSmallTime{\Tbase}%
            && \ifnum \xintiiGtorEq{\Wchils}{\max}=1 \relax \bfseries \fi 
            \numprint{\Wchils} & \checkSmallTime{\Tchils}%
            \\
    }
    
\end{longtable}%
}
\end{landscape}
\DTLgdeletedb{DATA}

\begin{landscape}
\begin{table}[ht]
    \centering    
        \caption{Average solution weight $\w$ and time $t$ (in seconds) required to compute it for the \emph{reduced} instances from set one.
        The \textbf{best} solutions among all algorithms are marked bold.}\label{tab:kernel_soa_set1} 
        
\DTLloaddb{DATA}{data/datatool/kernel_data/soa_set1Kernel}
\npthousandsep{\,}
\setlength{\tabcolsep}{3pt} 

\footnotesize{
\begin{tabular}{clrrrrrrrrrrrrrr}
        &&
        \multicolumn{2}{c}{\htwis} &&
        \multicolumn{2}{c}{\hils} &&
        \multicolumn{2}{c}{\mmwiss} &&
        \multicolumn{2}{c}{\baseline} &&
        \multicolumn{2}{c}{\pils} \\
        \cmidrule{3-4}\cmidrule{6-7}\cmidrule{9-10}\cmidrule{12-13}\cmidrule{15-16}\\[-1em]
        \multicolumn{2}{l}{\textbf{Instance}} & \multicolumn{1}{c}{$\w$} & \multicolumn{1}{c}{$t$}
        && \multicolumn{1}{c}{$\w$} & \multicolumn{1}{c}{$t$} 
        && \multicolumn{1}{c}{$\w$} & \multicolumn{1}{c}{$t$}
        && \multicolumn{1}{c}{$\w$} & \multicolumn{1}{c}{$t$} 
        && \multicolumn{1}{c}{$\w$} & \multicolumn{1}{c}{$t$} \\[-.2em]
        \cmidrule{3-4}\cmidrule{6-7}\cmidrule{9-10}\cmidrule{12-13}\cmidrule{15-16}\\[-1em]
        
    \DTLforeach*{DATA}{\instance=graph,\Whtwis=htwisWeight,\Thtwis=htwisTime,\Whils=hilsWeight,\Thils=hilsTime,\Wmmwis=mmwisWeight,\Tmmwis=mmwisTime,\Wchils=chilsWeight,\Tchils=chilsTime,\Wbase=baseWeight,\Tbase=baseTime,\max=MaxWeight}{%
        \IfSubStr{\instance}{body}{\multirow{3}{*}{\rotatebox{90}{\textbf{fe}}}}{}%
        \IfSubStr{\instance}{skitter}{\\[-.6em] \multirow{5}{*}{\rotatebox{90}{\textbf{snap}}}}{}%
        \IfSubStr{\instance}{fl2010}{\\[-.6em] \multirow{1}{*}{\rotatebox{90}{\textbf{ssmc}}}}{}%
        \IfSubStr{\instance}{alabama}{\\[-.6em] \multirow{23}{*}{\rotatebox{90}{\textbf{osm}}}}{}%
                & \niceGraphName{\instance}%
                & \ifnum \xintiiGtorEq{\Whtwis}{\max}=1 \relax \bfseries \fi 
                \numprint{\Whtwis} & \checkSmallTime{\Thtwis}%
                && \ifnum \xintiiGtorEq{\Whils}{\max}=1 \relax \bfseries \fi 
                \numprint{\Whils} & \checkSmallTime{\Thils}%
                && \ifnum \xintiiGtorEq{\Wmmwis}{\max}=1 \relax \bfseries \fi 
                \numprint{\Wmmwis} & \checkSmallTime{\Tmmwis}%
                && \ifnum \xintiiGtorEq{\Wbase}{\max}=1 \relax \bfseries \fi 
                \numprint{\Wbase} & \checkSmallTime{\Tbase}%
                && \ifnum \xintiiGtorEq{\Wchils}{\max}=1 \relax \bfseries \fi 
                \numprint{\Wchils} & \checkSmallTime{\Tchils}%
                \\
  }
\end{tabular}
}
\DTLgdeletedb{DATA}

\end{table}
\end{landscape}

\npthousandsep{\,}
\begin{filecontents*}{data/graphs/graph_reduction_data.csv}
graph,n,m,kn,km,offset,time
fe-body-uniform,45087,163734,395,3551,1672469,0.595169
fe-ocean-uniform,143437,409593,0,0,7248581,3.55729
fe-pwt-uniform,36519,144794,15239,106943,813076,2.16788
fe-rotor-uniform,99617,662431,89294,690242,435335,5.05004
fe-sphere-uniform,16386,49152,0,0,617816,0.517799
mesh-blob-uniform,16068,24102,0,0,855547,0.0322032
mesh-buddha-uniform,1087716,1631574,0,0,57555880,2.55494
mesh-bunny-uniform,68790,103017,0,0,3686960,0.131429
mesh-dragonsub-uniform,600000,900000,0,0,32213898,1.37619
mesh-dragon-uniform,150000,225000,0,0,7956530,0.198425
mesh-ecat-uniform,684496,1026744,0,0,36650298,2.12462
mesh-fandisk-uniform,8634,12818,0,0,463288,0.0187721
mesh-feline-uniform,41262,61893,0,0,2207219,0.104564
mesh-gameguy-uniform,42623,63850,0,0,2325878,0.0692439
mesh-gargoyle-uniform,20000,30000,0,0,1059559,0.0408518
mesh-turtle-uniform,267534,401178,0,0,14263005,0.423871
osm-alabama-AM2,1164,19386,0,0,174309,0.012171
osm-alabama-AM3,3504,309664,815,61637,173907,1.45062
osm-california-AM2,231,3074,0,0,47153,0.0190418
osm-california-AM3,587,27536,359,19705,35724,0.205125
osm-canada-AM3,943,20241,0,0,86018,0.232315
osm-colorado-AM3,538,8365,0,0,54741,0.0113389
osm-district-of-columbia-AM1,2500,24651,0,0,196475,0.218856
osm-district-of-columbia-AM2,13597,1609795,2361,274529,178323,5.76413
osm-district-of-columbia-AM3,46221,27729137,26728,13406329,133708,290.369
osm-florida-AM3,2985,154043,658,41043,226379,0.52137
osm-georgia-AM3,1680,74126,481,31388,207992,0.384872
osm-greenland-AM2,686,50218,0,0,10718,0.08231
osm-greenland-AM3,4986,3652361,3402,2141096,7572,27.3877
osm-hawaii-AM2,2875,265158,0,0,125284,0.15088
osm-hawaii-AM3,28006,49444921,22849,38081016,96234,935.961
osm-idaho-AM3,4064,3924080,2920,2413139,70960,40.0338
osm-kansas-AM3,2732,806912,1209,292997,84662,3.78747
osm-kentucky-AM2,2453,643428,0,0,97397,0.274048
osm-kentucky-AM3,19095,59533630,15943,53739580,91864,1510.41
osm-louisiana-AM3,1162,37077,0,0,60024,0.0765102
osm-maine-AM3,143,850,0,0,26734,0.000829935
osm-maryland-AM3,1018,95415,0,0,45496,0.0996959
osm-massachusetts-AM2,1339,35449,0,0,140095,0.030324
osm-massachusetts-AM3,3703,551491,1626,340652,136331,2.98373
osm-mexico-AM2,516,9411,0,0,94834,0.00759721
osm-mexico-AM3,1096,47131,447,20443,86467,0.286542
osm-minnesota-AM3,683,34188,0,0,32787,0.071728
osm-montana-AM3,837,69293,382,39859,55531,0.252079
osm-nevada-AM3,569,15016,0,0,52036,0.0165861
osm-new-hampshire-AM3,1107,18021,0,0,116060,0.121337
osm-new-york-AM2,224,6399,0,0,14330,0.0143211
osm-new-york-AM3,837,88728,526,50572,11447,0.416003
osm-north-carolina-AM3,1557,236739,997,133106,39783,0.872409
osm-ohio-AM3,482,11376,0,0,52634,0.0524528
osm-oregon-AM3,5588,2912701,3159,1817045,162009,23.7533
osm-pennsylvania-AM3,1148,26464,0,0,143870,0.119685
osm-puerto-rico-AM3,494,26926,0,0,33590,0.0973129
osm-rhode-island-AM2,2866,295488,498,46183,174013,0.522558
osm-rhode-island-AM3,15124,12622219,12189,11225470,144106,153.055
osm-tennessee-AM3,212,3215,0,0,32276,0.0134401
osm-utah-AM2,589,4692,0,0,95087,0.00160289
osm-utah-AM3,1339,42872,0,0,98847,0.283235
osm-vermont-AM2,766,37607,0,0,59310,0.113445
osm-vermont-AM3,3436,1136164,1946,526463,52788,5.49778
osm-virginia-AM2,2279,60040,0,0,295867,0.056628
osm-virginia-AM3,6185,665903,2385,293068,272730,2.5103
osm-washington-AM2,3025,152449,0,0,305619,0.100536
osm-washington-AM3,10022,2346213,6671,1955568,264406,14.5531
osm-west-virginia-AM3,1185,125620,864,98399,37678,0.516044
snap-as-skitter-uniform,1696415,11095298,1951,27277,124100687,3.43416
snap-com-amazon,334863,925869,0,0,19271031,0.322631
snap-loc-gowalla-edges,196591,950327,304,2938,12268682,0.327421
snap-roadNet-CA-uniform,1965206,2766607,0,0,111360828,2.33154
snap-roadNet-PA-uniform,1088092,1541898,0,0,61731589,1.18514
snap-roadNet-TX-uniform,1379917,1921660,0,0,78599946,1.45765
snap-soc-LiveJournal1-uniform,4847571,42851237,1715,28390,283975051,23.4735
snap-soc-pokec-relationships-uniform,1632803,22301964,782676,14234808,53362233,643.998
snap-web-BerkStan-uniform,685230,6649470,0,0,43907482,3.86781
snap-web-Google-uniform,875713,4322051,0,0,56326504,2.02972
snap-web-NotreDame-uniform,325729,1090108,97,1357,26013618,0.622788
snap-web-Stanford-uniform,281903,1992636,0,0,17792930,1.39353
ssmc-ca2010,710145,1744683,0,0,16869550,7.34339
ssmc-fl2010,484481,1173147,2045,8459,8707619,3.61887
ssmc-ga2010,291086,709028,0,0,4644417,0.891715
ssmc-il2010,451554,1082232,0,0,5998539,5.63136
ssmc-nh2010,48837,117275,0,0,588996,0.16351
ssmc-ri2010,25181,62875,0,0,459275,0.18774
\end{filecontents*}
\DTLloaddb{GRAPHS}{data/graphs/graph_reduction_data.csv}

\setlength{\tabcolsep}{4pt} 
\begin{longtable}{llrrrrrr}
\caption{Detailed graph properties for the first set consisting of \textit{fe}, \textit{mesh}, \textit{osm}, \textit{snap} and \textit{ssmc} instances. Graphs marked with a $\star$ are part of the parameter tuning set. We also present the number of vertices $n_K$ and edges $m_K$ of the instances reduced by \lnr~with the configuration \textit{Fast - initial tight} as well as the computed offset and running time $t_{red}$ in seconds.} \label{tab:graphs} \\
\multicolumn{2}{c}{\textbf{Instance}} & $n$ & $m$ &$n_K$&$m_K$& offset & $t_{red}$\\ 
\midrule
\endfirsthead

\multicolumn{8}{c}%
{\tablename\ \thetable\ -- \textit{Continued from previous page}} \\
\multicolumn{2}{c}{\textbf{Instance}} & $n$ & $m$ &$n_K$&$m_K$& offset & $t_{red}$\\ 
\midrule
\endhead

\midrule \multicolumn{8}{r}{\textit{Continued on next page}} \\
\endfoot

\endlastfoot

\DTLforeach*{GRAPHS}{\g=graph,\n=n,\m=m,\kn=kn, \km=km, \offset=offset,\time=time}{%
                
        \IfSubStr{\g}{skitter}{\\[-.7em] \textbf{snap}}{}%
        \IfSubStr{\g}{body}{ \textbf{fe}}{}%
        \IfSubStr{\g}{alabama-AM2}{\\[-.7em] \textbf{osm}}{}%
        \IfSubStr{\g}{blob}{\\[-.7em] \textbf{mesh}}{}%
        \IfSubStr{\g}{ca2010}{\\[-.7em] \textbf{ssmc}}{}%
    &  \CheckMarkGraphsParameter{\g}%
            \ifthenelse{\equal{\result}{true}}{%
            \textit{\niceGraphName{\g}}$^{\star}$}{\textit{\niceGraphName{\g}}}%
    & \numprint{\n} & \numprint{\m} & \numprint{\kn} & \numprint{\km} & \numprint{\offset} & \nptwo{\time}  \\ 
}

\end{longtable}
\DTLgdeletedb{GRAPHS}

\newpage
\begin{table}[ht]
    \centering    
    \caption{
    Detailed graph properties for the second set of vehicle routing instances. We also present the number of vertices $n_K$ and edges $m_K$ of the instances reduced by \lnr~with the configuration \textit{Fast - initial tight} as well as the computed offset and running time $t_{red}$ in seconds.
    }\label{tab:vr_graphs} 
    
\begin{filecontents*}{data/mwis_results_vr/results_chils_kernel_vr_exhaustive_initial_tight.csv}
id,graph,n,m,kn,km,offset,reductiontime,kernelWeight,kernelTime,fullWeight,fullTime
gnninitial-tight-t36001,CR-S-L-1,863368,331203970,855833,328396194,26070,3671.65,,,,
gnninitial-tight-t36000,CR-S-L-2,880974,342158741,873941,339503913,20282,5234.35,5695181,"3396,8905",5715463,"8631,2405"
gnninitial-tight-t36000,CR-S-L-4,881910,344057350,874668,341151292,19400,5315.89,5679025,"3578,8289",5698425,"8894,7189"
gnninitial-tight-t36000,CR-S-L-6,578244,219717582,572888,217572336,25208,3056.68,3865525,"3505,078",3890733,"6561,758"
gnninitial-tight-t36002,CR-S-L-7,270067,94109215,266591,92817339,26955,858.967,,,,
gnninitial-tight-t36000,CR-T-C-1,602472,194753152,594883,192451235,59577,2668.29,4642971,"3514,9923",4702548,"6183,2823"
gnninitial-tight-t36000,CR-T-C-2,652497,215694927,645058,213228808,32598,2870.86,4885659,"3569,3698",4918257,"6440,2298"
gnninitial-tight-t36000,CR-T-D-4,651861,220480534,644845,218142885,22898,2342.34,4841960,"3476,9374",4864858,"5819,2774"
gnninitial-tight-t36000,CR-T-D-6,381380,115082762,376351,113599041,31498,1393.38,2967733,"3520,6574",2999231,"4914,0374"
gnninitial-tight-t36003,CR-T-D-7,163809,43028583,160630,42184292,29046,341.149,,,,
gnninitial-tight-t36000,CW-S-L-1,411950,283860106,409559,282656866,11718,4465.88,1634830,"3582,7398",1646548,"8048,6198"
gnninitial-tight-t36000,CW-S-L-2,443404,315569883,441093,314281707,8310,5816.47,1715604,"3504,212",1723914,"9320,682"
gnninitial-tight-t36000,CW-S-L-4,430379,303042962,427984,301676399,6009,5262.94,1729624,"3574,368",1735633,"8837,308"
gnninitial-tight-t36000,CW-S-L-6,267698,171132761,265820,170150154,9603,2512.19,1157171,"3550,9348",1166774,"6063,1248"
gnninitial-tight-t36000,CW-S-L-7,127871,78459291,126866,77941263,4283,859.786,586551,"3476,472",590834,"4336,258"
gnninitial-tight-t36000,CW-T-C-1,266403,144634578,264389,143674135,10530,2064.78,1312640,"3571,0232",1323170,"5635,8032"
gnninitial-tight-t36000,CW-T-C-2,194413,111098006,192871,110337336,10128,1387.53,925779,"3529,6326",935907,"4917,1626"
gnninitial-tight-t36000,CW-T-D-4,83091,37881529,82110,37453654,3108,356.973,456435,"3477,3139",459543,"3834,2869"
gnninitial-tight-t36000,CW-T-D-6,83758,38781839,82723,38299758,5338,400.99,454889,"2836,081",460227,"3237,071"
gnninitial-tight-t36000,MR-D-03,21499,130508,17390,131591,832349816,1.06007,926425922,"2408,1939",1758775738,"2409,25397"
gnninitial-tight-t36000,MR-D-05,27621,236044,24250,236758,676962870,1.8606,1112981317,"3546,2987",1789944187,"3548,1593"
gnninitial-tight-t36000,MR-D-FN,30467,296369,27301,297905,639582861,2.4617,1162144467,"3385,0036",1801727328,"3387,4653"
gnninitial-tight-t36000,MR-W-FN,15639,126800,14657,128925,691250460,0.752507,4695592321,"62,9543",5386842781,"63,706807"
gnninitial-tight-t36000,MT-D-01,979,3125,0,0,238166485,0.0167542,,,238166485,"0,0167542"
gnninitial-tight-t36000,MT-D-200,10880,505359,10676,493364,6914835,2.19186,280171607,"257,3448",287086442,"259,53666"
gnninitial-tight-t36000,MT-D-FN,10880,604041,10682,595799,45784382,2.32549,245082561,"524,2754",290866943,"526,60089"
gnninitial-tight-t36000,MT-W-01,1006,2411,0,0,312121568,0.00254798,,,312121568,"0,00254798"
gnninitial-tight-t36000,MT-W-200,12320,515871,11472,474079,85332696,2.32752,298719461,"743,9686",384052157,"746,29612"
gnninitial-tight-t36000,MT-W-FN,12320,553895,11456,503617,128611366,2.4324,262258525,"167,4501",390869891,"169,8825"
gnninitial-tight-t36000,MW-D-01,3988,13556,2539,13853,307196278,0.177796,169244378,"1205,3992",476440656,"1205,576996"
gnninitial-tight-t36000,MW-D-20,20054,606318,19473,600915,104526987,3.14987,422121342,"3518,6864",526648329,"3521,83627"
gnninitial-tight-t36000,MW-D-40,33563,1879303,33062,1870537,74805076,8.33306,463604510,"3470,809",538409586,"3479,14206"
gnninitial-tight-t36000,MW-D-FN,47504,4017196,47070,4005378,60228177,18.7965,484018312,"3566,9589",544246489,"3585,7554"
gnninitial-tight-t36000,MW-W-01,3079,22664,2844,22765,169735181,0.162872,1100570771,"50,009",1270305952,"50,171872"
gnninitial-tight-t36000,MW-W-05,10790,485261,10667,481661,25259548,1.68331,1302784237,"285,4951",1328043785,"287,17841"
gnninitial-tight-t36000,MW-W-10,18023,1451813,17921,1445959,16146685,4.80567,1326663269,"1636,0561",1342809954,"1640,86177"
gnninitial-tight-t36000,MW-W-FN,22316,2275623,22215,2269047,15670946,9.86309,1335147597,"371,4639",1350818543,"381,32699"
\end{filecontents*}
\DTLloaddb{VRG}{data/mwis_results_vr/results_chils_kernel_vr_exhaustive_initial_tight.csv}
\npthousandsep{\,}
    \begin{tabular}{rlrrrrrr}
    \multicolumn{2}{c}{\textbf{Instance}} & $n$ & $m$ &$n_K$&$m_K$& offset & $t_{red}$\\ 
        \midrule
        \DTLforeach*{VRG}{\g=graph,\n=n,\m=m, \kn=kn, \km=km,\offset=offset,\time=reductiontime}{%
           \CheckMarkGraphsParameter{\g}%
                    \ifthenelse{\equal{\result}{true}}{%
                    \textbf{$\star$}}{}
    & \textit{\niceGraphName{\g}} & \numprint{\n} & \numprint{\m} & \numprint{\kn} & \numprint{\km} & \numprint{\offset} & \checkSmallTime{\time}  \\ 
        }\\[-.8em]
    \end{tabular}
\DTLgdeletedb{VRG}

\end{table}

\end{document}